\documentclass[11pt,a4paper]{article}

\usepackage[english]{babel}
\usepackage[margin=1in]{geometry}
\usepackage[latin1]{inputenc}
\usepackage{amsmath}
\usepackage{amsthm}
\usepackage{amsfonts}
\usepackage{amssymb}
\usepackage[makeroom]{cancel}
\usepackage{graphicx}
\usepackage{epsfig}
\usepackage{fullpage}
\usepackage{color}
\usepackage{bm}
\usepackage{booktabs}
\usepackage{colortbl}
\usepackage[cmtip,all]{xy}
\usepackage{hyperref}
\usepackage{caption,subcaption}
\usepackage{stmaryrd}
\usepackage{float}
\usepackage{wrapfig}
\usepackage{caption}
\usepackage{cite}
\usepackage[title,titletoc]{appendix}
\usepackage[colorinlistoftodos,prependcaption,textsize=tiny]{todonotes}
\usepackage{cleveref}
\usepackage{cases}
\usepackage{enumitem}
\usepackage{ulem}
\usepackage{tikz}
\usepackage{ifthen}
\DeclareMathSizes{10}{10}{10}{10}
\title{A posteriori error estimates for the Richards equation \footnotetext{This project has received funding by the European Research Council (ERC) under the European Union's Horizon 2020 research and innovation program (grant agreement No 647134)}}
\usepackage[noblocks]{authblk}
\author[2,1]{K. Mitra}
\author[1,3]{M. Vohral\'ik}
\affil[1]{Inria, 2 rue Simone Iff, 75589 Paris, France}
\affil[2]{Radboud University, Heyendaalseweg 135, 6525 AJ Nijmegen, The Netherlands}
\affil[3]{CERMICS, Ecole des Ponts, 77455 Marne-la-Vall\'ee, France}
\numberwithin{equation}{section}
\makeatletter
\DeclareMathSizes{\@xpt}{\@xpt}{6}{5}
\makeatother

\newcounter{assumption}% for numbering assumprions

 \stepcounter{assumption}
 
 \newcounter{Lscassumption}
 \stepcounter{Lscassumption}
  \newcounter{properties}
 \stepcounter{properties}
  \newcounter{fproperties}
  \stepcounter{fproperties}
 
\newtheorem{theorem}{Theorem}[section]

\newtheorem{proposition}{Proposition}[section]

\newtheorem{remark}[theorem]{Remark}
\theoremstyle{definition}
\newtheorem{definition}[theorem]{Definition}

\newtheorem{algo}{Algorithm}[section]
\def \a  {\alpha}

\def \b  {\beta}
\def \g  {\gamma}

\def \f  {\varphi}
\def \vr  {\varrho}
\def \k  {\kappa}
\def \l  {\lambda}
\def \om {\omega}
\def \Om {\Omega}

\def \s  {\sigma}
\def \t  {\tau}

\def \z  {\zeta}

\def \vg {\bm{g}}
\def \vs {\varsigma}
\def \del {\nabla}

\def \p  {\partial}
\def \R  {\mathbb{R}}
\def \N  {\mathbb{N}}
\def \K {\mathbf{\bar{K}}}
\def \Kn {{H^{-1}_{\!_\K}}}
\def \Kh {{H^{1}_{\!_\K}}}

\def \J {\mathcal{J}}

\def \Sf {\theta}
\def \Gf {\mathcal{G}}

\def \Lam  {{\Lambda}}

\def \W {{\bf \mathcal{Y}}}
\def \X {{\bf \mathcal{X}}}

\def \dst {{\mathrm{dist}}}
\DeclareMathOperator*{\argmin}{argmin}
\newcommand{\dd}{\mathrm{d}}

\newcommand{\eval}[2]{\left. #1\right|_{#2}}

\newcommand{\pinh}[2]{#1^{\bar{i}}_{#2,h}}

\newcommand{\calK}{\mathcal{K}}
\newcommand{\calT}{\mathcal{T}}
\newcommand{\calE}{\mathcal{E}}
\newcommand{\calP}{\mathcal{P}}

\newcommand{\calR}{\mathcal{R}}

\newcommand{\calV}{\mathcal{V}}
\newcommand{\RTN}{\mathbf{RTN}} %Local Raviart-Thomas-Nedelec space
\newcommand{\Hdiv}{\bm{H}(\mathrm{div},\Om)}
\newcommand{\Hdivoma}{\bm{H}(\mathrm{div},\oma)}
\newcommand{\CR}{\calT_{n}}
\newcommand{\Ta}{\calT^{\ta}_n}

\newcommand{\Vn}{V_{n,h}}

\newcommand{\elCR}{K}

\newcommand{\res}{\mathcal{R}}
\newcommand{\ta}{\mathbf{a}}
\newcommand{\psia}{\psi_{\ta}}

\newcommand{\oma}{{\om_{\ta}}}

\newcommand{\PiL}{\Pi_{n,h}}

\newcommand{\PiH}{\Pi_{n,h}^{{\mathrm{RT}}}}

\newcommand{\calVh}{\mathcal{V}_n}
\newcommand{\calVhint}{\mathcal{V}^{\mathrm{int}}_n}
\newcommand{\calVhext}{\mathcal{V}^{\mathrm{ext}}_n}
\newcommand{\Va}{\bm{V}^{\ta}_{n,h}}
\newcommand{\Qa}{Q_{n,h}^{\ta}}

\newcommand{\stha}{\bm{\sigma}_{n,h}^{\ta}}

\newcommand{\sth}{\bm{\sigma}_{n,h}}
\newcommand{\sith}{\bm{\sigma}_{n,h}}

\newcommand{\tautha}{\bm{\tau}^{\ta}_{n,h}}

\newcommand{\gtautha}{g^{\ta}_{n,h}}

\newcommand{\Snt}{\mathcal{G}_{n,h}}
\newcommand{\Fnt}{\bm{F}_{n,h}}
\newcommand{\etaEq}{\eta^{\mathrm{F}}_{n,h,K}}
\newcommand{\etaJh}{\eta^{\mathrm{J},H^1}_{n,h,K}}
\newcommand{\etaJl}{\eta^{\mathrm{J},L^2}_{n,h,K}}
\newcommand{\etaOscT}{\eta^{\mathrm{qd},t}_{n,h,K}}

\newcommand{\etaOscS}{\eta^{\mathrm{qd},\mathcal{G}}_{n,h,K}}
\newcommand{\etaOscF}{\eta^{\mathrm{qd},\bm{F}}_{n,h,K}}
\newcommand{\etaOscSa}{\eta^{\mathrm{qd},\mathcal{G}}_{n,h,\oma}}
\newcommand{\etaOscFa}{\eta^{\mathrm{qd},\bm{F}}_{n,h,\oma}}
\newcommand{\etaOscTa}{\eta^{\mathrm{qd},t}_{n,h,\oma}}
\newcommand{\etaOscTom}{\eta^{\mathrm{qd},t}_{n,h,\Om}}

\newcommand{\etaOsc}{\eta^{\mathrm{osc}}_{n,\om}}

\newcommand{\etaDeg}{\eta^{\mathrm{deg}}}

\newcommand{\etaOscInitL}{\eta^{\mathrm{ini},L^2}}
\newcommand{\etaOscInitH}{\eta^{\mathrm{ini},H^{-1}}}
\newcommand{\norm}[1]{{\left\vert\kern-0.25ex\left\vert\kern-0.25ex\left\vert #1 
    \right\vert\kern-0.25ex\right\vert\kern-0.25ex\right\vert}}

\providecommand{\keywords}[1]
{
  \small	
  \textbf{\textit{Keywords---}} #1
}

\begin{document}
\maketitle

\vspace{-2em}
\begin{abstract}
The Richards equation is commonly used to model the flow of water and air through  soil, and it serves as a gateway equation for multiphase flows through porous media. It is a nonlinear advection--reaction--diffusion equation that exhibits both parabolic--hyperbolic and parabolic--elliptic kinds of degeneracies. In this study, we provide reliable, fully computable, and locally space--time efficient a posteriori error bounds for numerical approximations of the fully degenerate Richards equation. For showing global reliability, a nonlocal-in-time error estimate is derived individually for the time-integrated $H^1(H^{-1})$, $L^2(L^2)$, and the $L^2(H^1)$ errors. A maximum principle and a degeneracy estimator are employed for the last one. Global and local space--time efficiency error bounds are then obtained in a standard $H^1(H^{-1})\cap L^2(H^1)$ norm. The reliability and efficiency norms employed coincide when there is no nonlinearity.  Moreover, error contributors such as flux nonconformity, time discretization, quadrature, linearization, and  data oscillation are identified and separated. The estimates are also valid in a setting where iterative linearization with inexact solvers is considered. Numerical tests are conducted for nondegenerate and degenerate cases having exact solutions, as well as for a realistic case. It is shown that the estimators correctly identify the errors up to a factor of the order of unity.

\keywords{Richards equation, a-posteriori error estimates, nonlinear degenerate problems, flow through porous media, finite element method}
\end{abstract}

%\tableofcontents
\section{Introduction}
The Richards equation models flow of water through porous medium (e.g., soil) partially filled with air \cite{helmig1997multiphase,bear1972dynamics}. For a domain $\Om\subset \R^d$, $d\in \N$, and final time $T>0$, with water saturation $s$ and pressure $p$ being the primary unknowns, it equates
\begin{subequations}\label{eq:FullSystem}
\begin{align}
\p_t s -\del\cdot[\K(\bm{x}) \, \k(s)\, (\del p + \vg)]= f(s,\bm{x},t) \text{ in } \Om\times [0,T].\label{eq:Richards1}
\end{align}
Here,  space and time variables are denoted by $\bm{x}$ and $t$, respectively. The source term $f(s,\bm{x},t)$ represents contribution due to reaction/absorption. The gravity is represented by the constant vector $-\vg$.  The absolute permeability tensor $\K(\bm{x})$ and the relative permeability function $\k:[0,1]\to [0,1]$ are properties of the medium. Initial condition is provided for the saturation $s$, and homogeneous Dirichtlet boundary condition is provided for the pressure $p$, i.e.,
\begin{equation}\label{eq:IC}
s(\bm{x},0)=s_0(\bm{x}) \text{ for } \bm{x}\in \Om \text{ and } p=0  \text{ on } \p\Om\times (0,T].
\end{equation}
Dirichlet--Neumann mixed boundary conditions are also considered in the numerical \Cref{sec:Numerical}.
To close \eqref{eq:Richards1}--\eqref{eq:IC},  it is usually assumed that saturation and pressure are related algebraically (commonly referred to as the capillary pressure relationship \cite{helmig1997multiphase}), i.e.,  for a function $S:\R\to [0,1]$ one has
\begin{align}
s=S(p).\label{eq:CapPressure}
\end{align}
\end{subequations}
Here, we assume that the saturation $s$ is bounded in the closed interval $[0,1]$.
Equation \eqref{eq:Richards1} is obtained by combining the constitutive relation for the flux, stated by the Darcy law
$$
\bm{\sigma}:=-\K(\bm{x}) \k(s) (\del p + \vg),
$$
with the mass balance equation $\p_t s+ \del\cdot \bm{\s}=f(s,\bm{x},t)$. 
The Richards equation is important in modelling groundwater flow and various chemical and biological processes. It is a nonlinear advection--reaction--diffusion equation which degenerates into an elliptic equation if $S'(p)=0$ at some point of the domain. On the other hand, if $\k(s)=0$, then the equation becomes a first order ordinary differential equation (hyperbolic) with the loss of regularity  of the solution. Nonlinearity and degeneracy are the two main challenges in analysing the system \eqref{eq:FullSystem}. 

Existence of solutions for the Richards equation was shown in \cite{alt1983quasilinear,alt1984nonstationary}. However, in the degenerate case when $\k(s)=0$, only the existence of a weak limit can  be shown \cite{alt1984nonstationary}. Consequently, the pair $(s,p)$ might not satisfy \eqref{eq:FullSystem} in a weak sense. We give appropriate details in \Cref{sec:Formulations}. Uniqueness of solutions is proven  in \cite{otto1996l1} using the $L^1$-contraction method. 

Different spatial discretization methods have been designed for the Richards equation. Some notable examples are \cite{eymard1999finite} for finite volumes, \cite{nochetto1988approximation} for finite elements, \cite{arbogast1996nonlinear,radu2014convergence} for  mixed finite elements,  \cite{li2007local} for  the discontinuous Galerkin method, and \cite{klausen2008convergence} for multi-point flux approximations. Iterative linearization methods such as the Newton, Picard, J{\"a}ger--Ka{\v{c}}ur, and the L-schemes have been investigated in \cite{bergamaschi1999mixed}, \cite{celia_Picard}, \cite{Jager1991}, and \cite{list2016study,MITRA2018}, respectively, see also the references therein. An improvement of the Newton method was proposed in \cite{brenner2017improving} by parametrizing both the saturation and the pressure as functions of a separate primary variable. A comprehensive review of numerical methods for the Richards equation can be found in \cite{zha2019review}.

The theory of a posteriori estimates for elliptic differential equations is well studied, see, e.g. \cite{Ainsw_Oden_a_post_FE_00,Repin_book_08,Verf_13}.
A posteriori upper error bounds for the heat equation in the $L^2(H^1)\cap L^\infty(L^2)$ norm were derived in \cite{Pic_adpt_par_98}. In  \cite{Ver_a_post_heat_03}, global efficiency in space on every time step together with reliability are proven for the $L^2(H^1)\cap L^\infty(L^2)\cap H^1(H^{-1})$ norm. In \cite{Ern_Sme_Voh_heat_HO_Y_17}, a local efficiency estimate in space and in time is established for the norm further enriched by time jumps. 
 A general framework for obtaining rigorous a posteriori estimates for  nonlinear problems has been laid out in \cite{Verf_a_post_NL_par_I_98,Verf_a_post_NL_par_II_98}. 
However, the Lipschitz continuity and invertibility  of the operators associated with the differential equations are assumed, which limits the scope of the estimates.
A more specific result for the $p$-Laplacian problem is given in \cite{kreuzer2013reliable}. Using a formulation relying on the $N$-functions, the coercivity and Lipschitz-continuity of the flux function are shown with respect to the gradient. This makes it possible to derive a posteriori estimates for the problem.
 Estimators for nonlinear advection--diffusion equations were proposed in \cite{Dol_Ern_Voh_a_post_unst_13}.
Both upper and lower bounds (reliability and efficiency) were established, robust with respect to the nonlinearities and advection dominance, but for a weaker space--time mesh-dependent norm. Moreover, it was also assumed that the solutions belong to $H^1(L^2)$, which may not be the case for degenerate problems and/or if the initial condition is discontinuous. 
 Using entropy methods,  error estimates in the $L^1$-norm were derived in \cite{Ohl_a_post_FV_cell_CRD_01} for singularly perturbed nonlinear advection--diffusion problems. Degenerate parabolic equations were considered in \cite{Noch_Sch_Ver_apost_deg_par_99}. An $L^\infty(H^{-1})$ estimate was derived using dual equations of the diffusion problem. For problems having parabolic--hyperbolic degeneracy,  a  posteriori upper bounds on the $L^2(H^{-1})\cap L^\infty(H^{-1})$ norm combined with the time-integrated $L^2(L^2)$ norm of error  were derived  using Green's function  in \cite{Di_Pi_Voh_Yous_a_post_Stef_15} for a Stefan problem and on the $L^2(H^{-1})\cap L^2(H^1)\cap L^2(L^2)$ norm in \cite{cances2014posteriori} for two-phase flow through porous media. 
 For the Richards equation, a posteriori error upper bounds in the $L^2(H^1)\cap H^1(L^2)$ norm were derived in \cite{Ber_El_Al_Mghaz_a_post_Rich_14}. A regularization term was introduced to avoid degeneracy and to obtain $H^1(L^2)$ estimates. 
 
 In the present paper, we provide a posteriori error estimates for the Richards equation \eqref{eq:Richards1}. The main improvements in this study are: (a) \textit{Rigorous derivation of the upper as well as lower bounds of error} by the equivalence of the dual norm of the residual with an error metric that reduces to the $L^2(H^1)\cap H^1(H^{-1})\cap L^{\infty}(L^2)$ norm in the linear case. (b) Equivalence of the dual norm of the residual with \textit{fully computable} and \textit{locally space--time efficient estimates}. (c) \textit{ No higher-order regularity assumptions} such as the pressure in $L^2(H^2)\cap H^1(L^2)$ or that the initial condition is in $H^1$. (d) \textit{Inclusion of both the parabolic--hyperbolic and the parabolic--elliptic type of degeneracies}. This requires relaxing the assumptions on the associated functions such as $S'(p),\, \k(s)>0$, assumed for instance in \cite{Ber_El_Al_Mghaz_a_post_Rich_14,baron2017adaptive,cances2014posteriori}  in order to avoid the blow-up due to degeneracy. It poses a challenge particularly since the parabolic--hyperbolic degeneracy, stemming from $\k(s)=0$, causes  a loss of regularity of the solutions. To circumvent this issue, we assume instead that the initial saturation $s_0$ is bounded away from the degenerate value at $0$. With this assumption, a function $S_{\mathrm{m}}: [0,T]\to (0,1]$ is computed using maximum principle such that $S_{\mathrm{m}}(t)\leq s(\bm{x},t) \leq 1$ for all $(\bm{x},t)\in \Om\times [0,T]$. For the parabolic--elliptic degeneracy, a degeneracy estimator is introduced to provide an upper bound on the $L^2(H^1)$ norm of the error. (e) \textit{Rigorous inclusion} of \textit{linearization errors due to inexact solvers}, space and time \textit{adaptive meshes}, and implementation of \textit{adaptive linearization}. (f) It is shown numerically that despite nonlinearities, degeneracies, and heterogeneities, the effectivity index of the estimators lies  between 1 and 3 in most cases, even locally.

 The paper is organized as follows.  \Cref{sec:Richards} serves as a mathematical prologue to the Richards equation. The associated functions, relevant transformations, well-posedness results, and maximum principles are discussed in detail. In \Cref{sec:ErrorResidual}, lower and upper bounds on error by the dual norm of the residual are derived. The upper bound is provided separately for the $H^1(H^{-1})$, $L^2(L^2)$, and the $L^2(H^1)$    errors in a time-smoothened fashion, see \Cref{theo:UpperBound}. In \Cref{sec:FiniteElement}, a finite element  approximation to the Richards problem \eqref{eq:FullSystem} is considered, and some time-interpolations are discussed. These are used in  \Cref{sec:APosteriori} to compute the equilibrated flux and the a posteriori estimators. Reliability and local  space--time efficiency bounds are proven for the estimators. Finally, numerical results are presented in \Cref{sec:Numerical}. The theoretical findings are verified and the corresponding effectivity indices are obtained using a nondegenerate as well as a degenerate case with known exact solutions. To demonstrate the prowess of the estimators,  a realistic degenerate problem is analyzed in a heterogeneous, anisotropic domain, with discontinuous initial condition and mixed boundary condition. To conclude, it is shown in Appendix \ref{App:linear} how to take into account the additional errors from iterative linearization, whereas Appendix \ref{App:proof} collects some technical proofs.

\section{The Richards equation}\label{sec:Richards}
Here, we give a brief introduction to the Richards equation and state some of its properties important for our analysis.
\subsection{Basic notation}
\textbf{Spaces:} Let $\Om\subset \R^d$ be an open polytope with a Lipschitz-continuous boundary.
Let $(\cdot,\cdot)$ and $\|\cdot \|$ represent respectively the $L^2(\Om)$ inner product and norm; $(\cdot,\cdot)_\om$ and $\|\cdot\|_\om$ stand for the $L^2$-inner product and norm with respect to any Lipschitz subdomain $\om\subset \Om$. The Sobolev space $H^1(\Om)$ contains all functions $u\in L^2(\Om)$ such that the weak derivative $\del u\in \bm{L}^2(\Om;\R^d)$, and $H^1_0(\Om)$ is the subspace of $H^1(\Om)$ containing functions vanishing at the boundary $\p\Om$ in the trace sense. The space $H^{-1}(\Om)$ stands for the dual of $H^1_0(\Om)$, and $\langle \cdot,\cdot\rangle$ denotes the corresponding duality pairing.
With final time $T>0$ and  $L^2(0,T;V)$ denoting the $L^2$ Bochner space for a Banach space $V$, we introduce the Hilbert spaces 
\begin{align}
\X:=L^2(0,T;H^1_0(\Om)) \text{ and } \W:= \{u\in L^2(0,T;H^1(\Om)):\; \p_t u \in L^2(0,T;H^{-1}(\Om))\}.\label{eq:BochnerW}
\end{align}
\textbf{Inequalities:} For a Lipschitz subdomain $\om \subseteq \Om$ with diameter $h_\om$, let $u\in H^1(\om)$ be such that either $\int_{\om} u=0$ or the trace of $u$ is zero on a section  of $\p\om$ of nonzero measure. Then the Poincar\'e--Friedrichs inequality states that there exists a constant $C_{\mathrm{P},\om}>0$ such that
\begin{align}\label{eq:Poincare}
\|u\|_\om \leq C_{\mathrm{P},\om} h_\om \,  \|\del u\|_{\om}.
\end{align}
For a convex  $\om$  in the zero mean-value case, $C_{\mathrm{P},\om}$ can be taken as $\pi^{-1}$. 

\textbf{Notation:} Let $[\cdot]_+=\max(\cdot,0)$ and $[\cdot]_-=\min(\cdot,0)$ denote the positive and negative part functions respectively. In our notation, $a\lesssim b$ will refer to the inequality $a\leq C b$, where $C>0$ is a constant that depends solely on the shape-regularity of the spatial meshes in the space dimension $d$, and on the ratio $K_{\mathrm{m}}\slash K_{\mathrm{M}}$ (see \ref{prop:Kabs} below). In particular, it is independent of mesh-size, time-step size, the functions $\k(\cdot),\, S(\cdot)$, $f$, and the polynomial degrees associated with the numerical scheme.

\subsection{Assumptions on the data}
\begin{figure}[h!]
\begin{subfigure}{.48\textwidth}
\includegraphics[scale=.4]{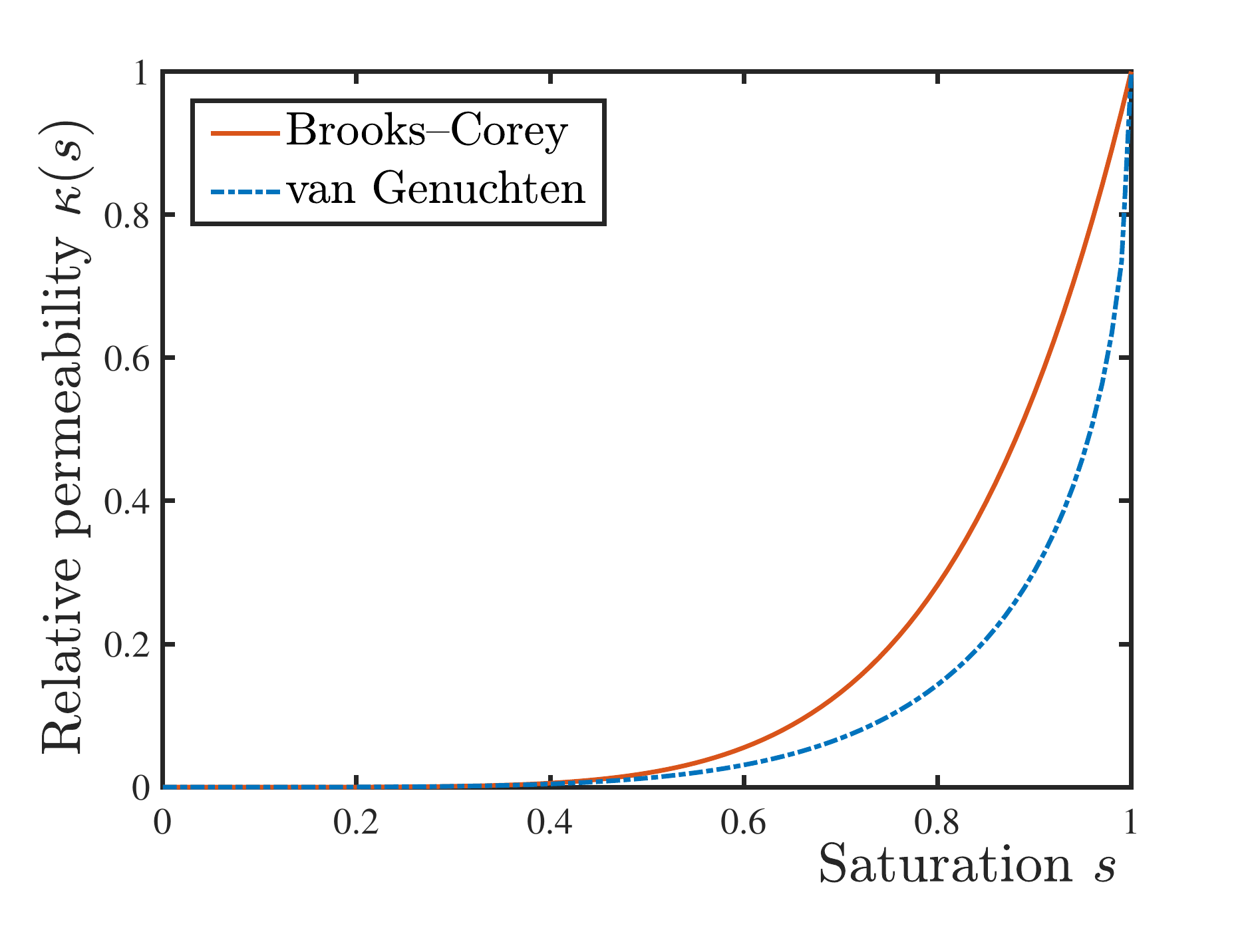}
\end{subfigure}
\begin{subfigure}{.48\textwidth}
\includegraphics[scale=.4]{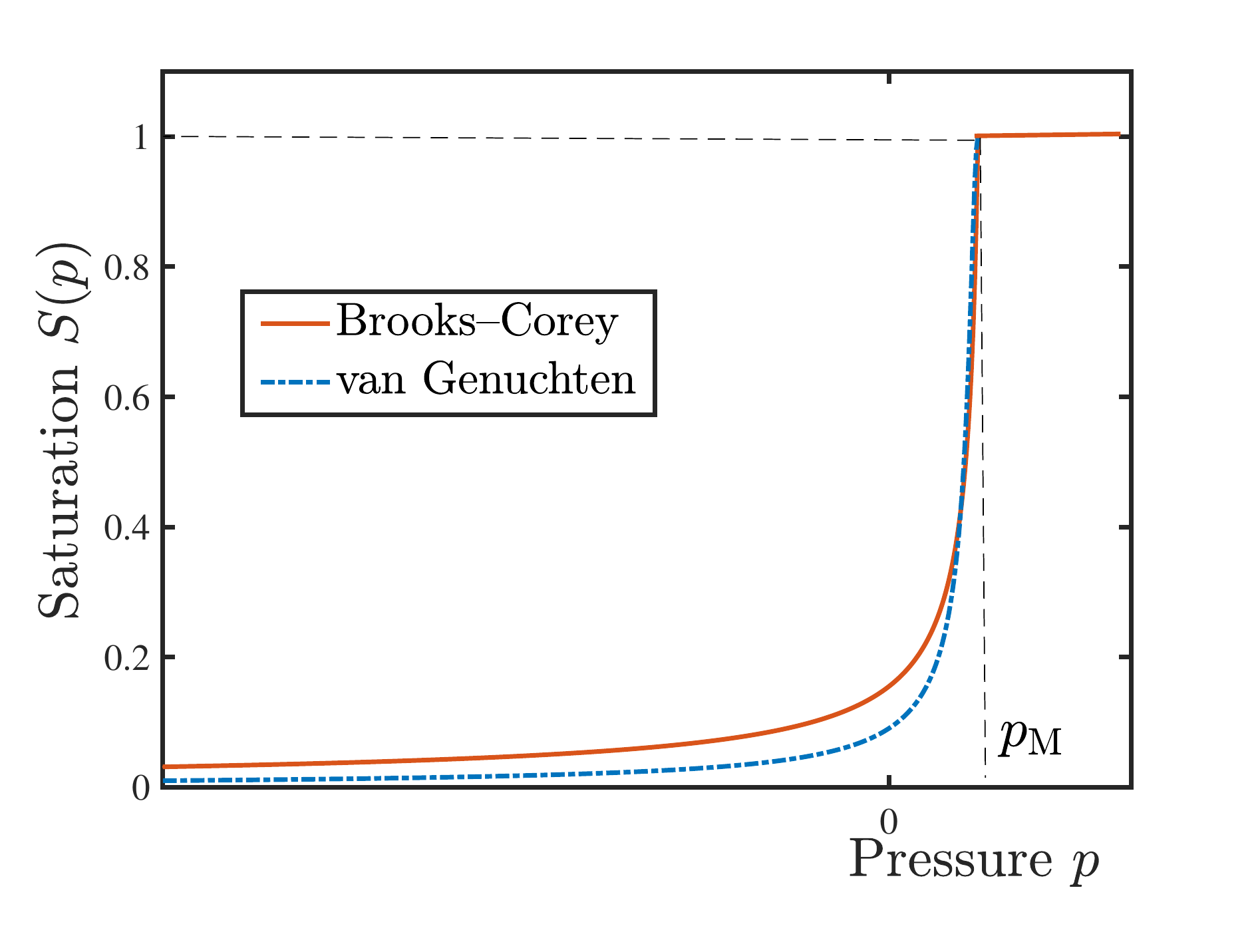}
\end{subfigure}
\caption{The functions $\k(s)$ (left) and $S(p)$ (right) as modeled by the Brooks--Corey \eqref{eq:BrooksCorey} and the van Genuchten \eqref{eq:vanGenuchten} models. The parameters are $\lambda_1=0.75$ and $\lambda_2=2$ taken from \cite{lenhard1989correspondence}, which gives a rather close match between the two models. For the heat equation, for comparison, $\k(s)= 1$ and $S(p)=p$.} \label{fig:KandS}
\end{figure}

We assume the following properties for the data in \eqref{eq:FullSystem}:
\begin{enumerate}[label=(P\theproperties)]
\itemsep .3em 
 \item The relative permeability function $\k$ is of the class $C^1([0,1])$ with $\k(0)\geq 0$, $\k(1)=1$, and $\k(0)< \k(s)< \k(1)$ for all $s\in (0,1)$.\label{prop:k}
 \stepcounter{properties}
  \item The saturation function $S$ is of the class $\mathrm{Lip}(\R)$ with Lipschitz constant $L_S>0$. It is either linear, or there exists a constant $p_\mathrm{M}\in (0,\infty]$ such that $\lim\limits_{p\searrow -\infty} S(p)=0$, and
  
 \begin{itemize}
 \item[(a)] $S|_{(-\infty,p_{\mathrm{M}}]}\in C^2((-\infty,p_{\mathrm{M}}])$, $S'(p) >0 \text{ for all } p< p_\mathrm{M}$, and $\lim\limits_{p\nearrow p_{\mathrm{M}}} S'(p)>0$;
 \item[(b)] $S(p)=1  \text{ and consequently } S'(p)=0 \text{ for all } p> p_\mathrm{M}.$
 \end{itemize}
  
  \label{prop:pc}
\stepcounter{properties}
\item The absolute permeability tensor $\K: \Om\mapsto \R^{d\times d}$ is piecewise constant in $\Om$, bounded, and satisfies the ellipticity condition, i.e., there exists positive constants $K_\mathrm{m},K_\mathrm{M}$ such that for any $\bm{\zeta} \in \R^d$,
$$
K_\mathrm{m} |\bm{\zeta}|^2 \leq 
\bm{\zeta}^\mathrm{T} \K(\bm{x}) \bm{\zeta}\leq K_\mathrm{M} |\bm{\zeta}|^2\quad  \text{ for almost all } \bm{x}\in \Om, 
$$
where $|\bm{\zeta}|$ is the Euclidean norm of $\bm{\zeta}$, i.e.,  $|\bm{\zeta}|=(\sum_{j=1}^d \zeta_j^2)^{\frac{1}{2}}$. Consequently, there exist unique positive-definite tensor-valued functions $\K^{\frac{1}{2}}$, $\K^{-\frac{1}{2}}$, and $\K^{-1}$.\label{prop:Kabs}\stepcounter{properties}\hspace{-.5em}
  \item The source term $f\in C^1([0,1]\times \Om\times \R)$ and there exists a function $f_{\mathrm{m}}\in C^1([0,1])$ such that $f_{\mathrm{m}}(\cdot)\leq \inf_{\bm{x}\in\Om, t\in \R^+} f(\cdot,\bm{x},t)$.
   \label{prop:f} \stepcounter{properties}
  \item The initial condition $s_0\in L^\infty(\Om)$ satisfies
  $${0<\mathrm{ess}\inf_{\bm{x}\in \Om}\{s_0(\bm{x})\}\leq \mathrm{ess}\sup_{\bm{x}\in \Om}\{s_0(\bm{x})\}\leq 1}.$$
   \label{prop:s0}\stepcounter{properties}
\end{enumerate}  
These assumptions are consistent with experiments, see e.g. \cite{helmig1997multiphase}. 
\begin{remark}[Choices for the functions $\k$ and $S$]
Two most commonly used models for the functions $\k(\cdot)$ and $S(\cdot)$ \cite{lenhard1989correspondence} are the Brooks--Corey model,
\begin{equation}
\k(s)=s^{\frac{2+3\lambda_1}{\lambda_1}}, \quad S(p)= (2- p\slash p_\mathrm{M})^{-\lambda_1} \text{ for } p\leq p_\mathrm{M}, \label{eq:BrooksCorey}
\end{equation}
and the van Genuchten model,
\begin{equation}
\k(s)=\sqrt{s}\,(1-(1-s^{1/\lambda_2})^{\lambda_2})^2, \quad S(p)=1\slash (1+ (p_{\mathrm{M}}-p)^{\frac{1}{1-\lambda_2}})^{\lambda_2} \text{ for } p\leq p_\mathrm{M},\label{eq:vanGenuchten}
\end{equation}
where $\lambda_1>0$ and $\lambda_2\in (0,1)$ are parameters. These functions $\k(\cdot)$ and $S(\cdot)$ are plotted in \Cref{fig:KandS} for $\lambda_1=0.75$ and $\lambda_2=2$. Observe that both the models satisfy assumptions \ref{prop:k}--\ref{prop:pc}.
\end{remark}

\subsection{Capillary pressure, diffusivity, total pressure, and auxiliary functions}
Here, we introduce some auxiliary functions that will be useful later.
\subsubsection{Capillary pressure function}
\begin{figure}[h!]
\centering
\includegraphics[scale=.4]{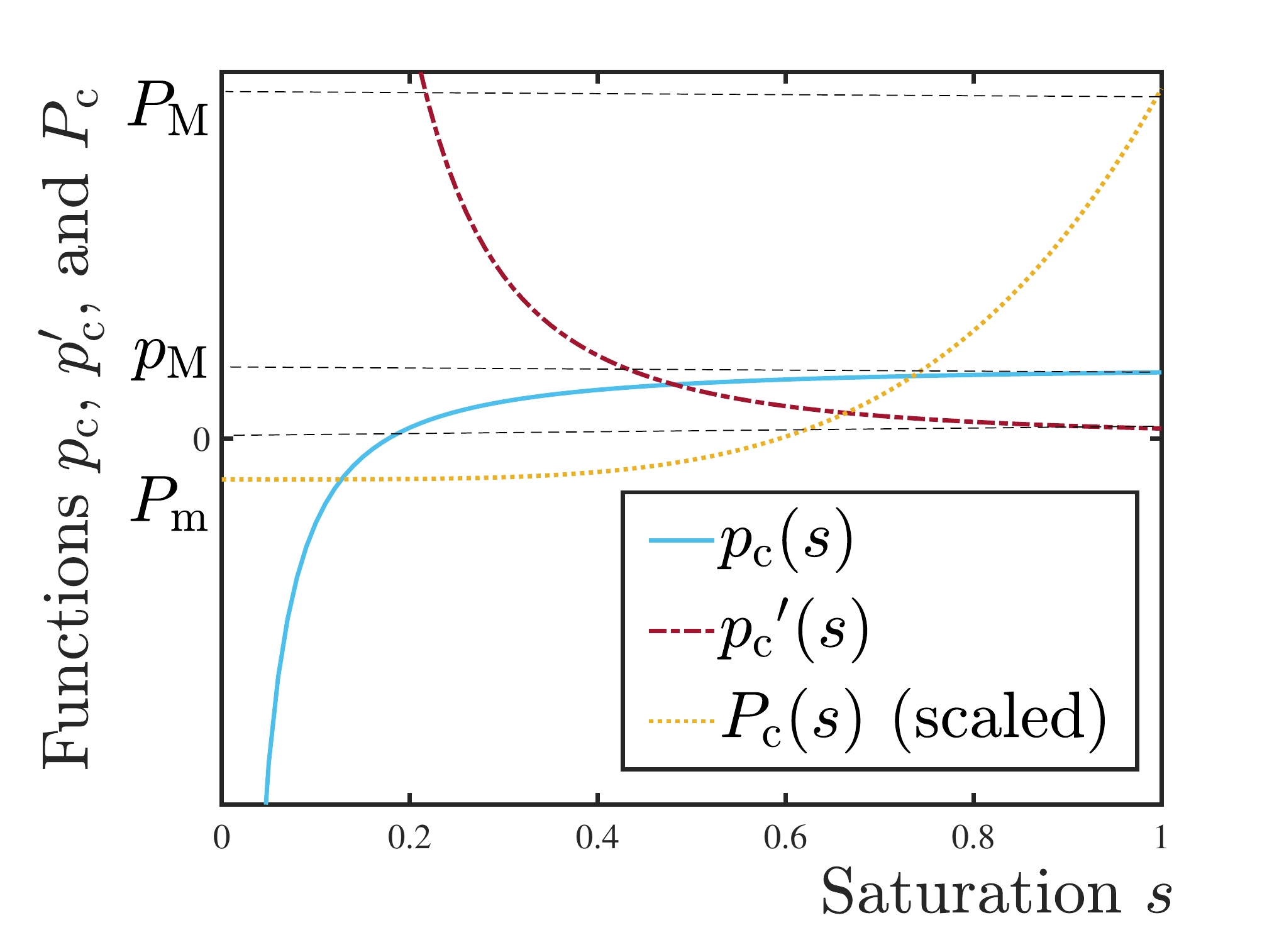}
\caption{The ${p_\mathrm{c}}$, ${p_\mathrm{c}}'$, and ${P_\mathrm{c}}$ functions for the Brooks--Corey model with $\lambda_1=0.75$.}\label{fig:pc}
\end{figure}

Since $S(\cdot)$ is a strictly increasing function in the interval $(-\infty,p_{\mathrm{M}}]$, its inverse 
\begin{subequations}\label{eq:capPressure}
\begin{equation}
p_{\mathrm{c}}(s):=S^{-1}(s)
\end{equation}
 is well-defined for $0<s\leq 1$. This is commonly known as the capillary pressure function. It is strictly increasing and $\lim_{s\searrow 0} p_{\mathrm{c}}(s)=-\infty$, see \Cref{fig:pc}. 
Using $p_{\mathrm{c}}(\cdot)$, the relation \eqref{eq:CapPressure} is alternatively stated as
\begin{equation}
p\begin{cases}
=p_{\mathrm{c}}(s) &\text{ if } 0<s< 1,\\
\in [p_{\mathrm{M}},\infty] &\text{ if } s=1.
\end{cases}
\end{equation}
\end{subequations}

\subsubsection{Diffusivity and total pressure functions} 
\begin{figure}[h!]
\begin{subfigure}{.48\textwidth}
\includegraphics[scale=.4]{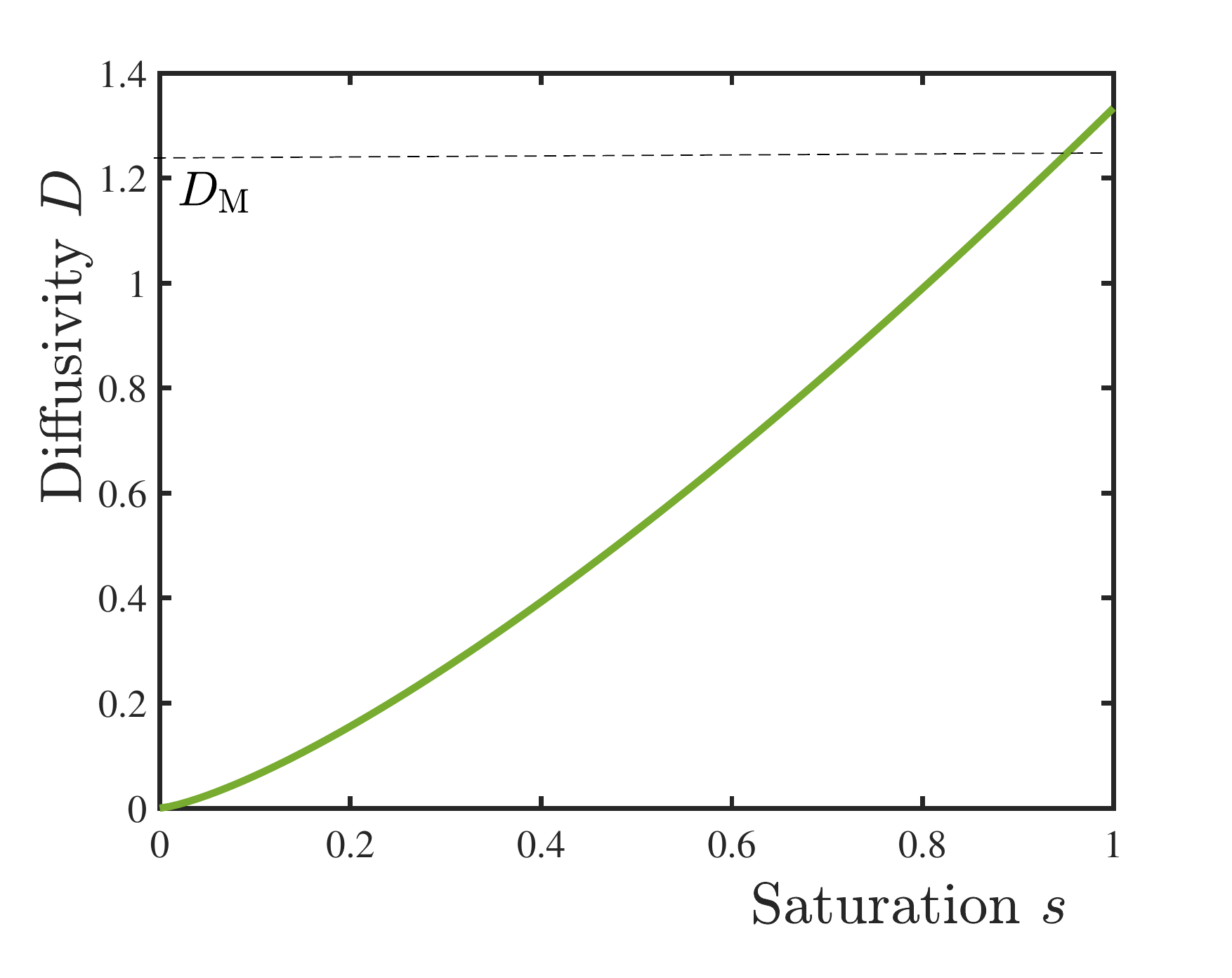}
\end{subfigure}
\begin{subfigure}{.48\textwidth}
\includegraphics[scale=.4]{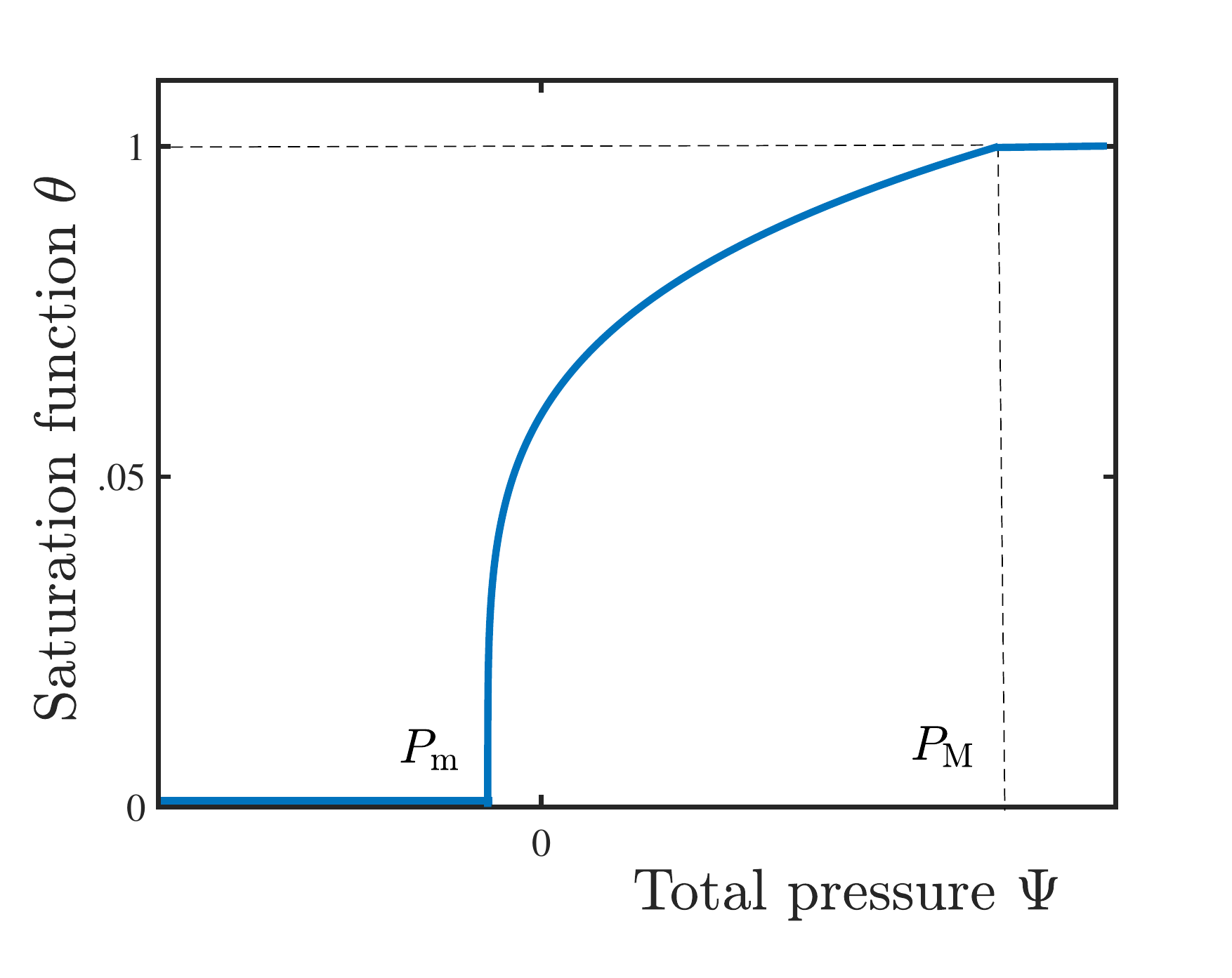}
\end{subfigure}
\caption{The functions $D(\cdot)$ and $\Sf(\cdot)$ for the Brooks--Corey model with $\lambda_1=0.75$.} \label{fig:DandTheta}
\end{figure}

We further introduce the diffusivity function $D:(0,1]\to \R^+$ as
\begin{equation}
D(s):=\k(s)\,{p_{\mathrm{c}}}'(s),\label{eq:DefDiffusivity}
\end{equation}
 and the total pressure function $P_{\mathrm{c}}:(0,1]\to \R$ (see \Cref{fig:pc}) as
\begin{align}
P_{\mathrm{c}}(s):=\int_{S(0)}^{s} D(\vr) \,\mathrm{d}\vr.\label{eq:totalPcS}
\end{align}
 The properties of $D$ and $P_c$ that follow from \ref{prop:k}--\ref{prop:pc} are 
\begin{equation}\label{eq:propD}
D\in C^1((0,1]); \; 0<D(s)<\infty \text{ for all } 0<s\leq 1; \text{ and  } \lim_{s\searrow 0} D(s)\geq0;
\end{equation}
whereas, $P_{\mathrm{c}}\in C^1((0,1])$  is strictly increasing since 
\begin{subequations}\label{eq:propPc}
\begin{align}
 {P_\mathrm{c}}'(s)=D(s),
\end{align}
and there exists fixed  $P_\mathrm{m},\,P_{\mathrm{M}}\in [-\infty,\infty)$ depending only upon $\k(\cdot)$ and $p_\mathrm{c}(\cdot)$ such that
\begin{align}\label{eq:DefPmPM}
 P_\mathrm{m}=\lim_{s \searrow 0} P_{\mathrm{c}}(s), \text{ and } P_\mathrm{M}=P_{\mathrm{c}}(1).
 \end{align}
\end{subequations}
Accordingly, an increasing and continuous function $\Sf: \R\to [0,1]$ is defined by 
\begin{equation}\label{eq:DefSfunc}
 \Sf(\Psi):=
\begin{cases}
0 &\text{ if } \Psi\leq P_\mathrm{m},\\
(P_{\mathrm{c}})^{-1}(\Psi) &\text{ if } P_\mathrm{m}<\Psi<P_\mathrm{M},\\
1 &\text{ if } \Psi\geq P_\mathrm{M}.
\end{cases}
 \end{equation} 
The plots of $D(\cdot)$ and $\Sf(\cdot)$ are shown in \Cref{fig:DandTheta}.

\begin{remark}[Properties of the function $\Sf$]
Observe from \eqref{eq:propPc}--\eqref{eq:DefSfunc} that,  
\begin{align}\label{eq:profSf}
\Sf'(\Psi)= \frac{1}{P'_{\mathrm{c}}(\Sf(\Psi))}=\dfrac{1}{D(\Sf(\Psi))} \text{ for all } \Psi\in (P_\mathrm{m},P_\mathrm{M}]. 
\end{align}
Consequently,  $\Sf|_{(P_\mathrm{m},P_\mathrm{M}]}\in C^1((P_\mathrm{m},P_\mathrm{M}])$. Moreover, it holds for all $\Psi > P_{\mathrm{m}}$ that
\begin{align}\label{eq:PsiandPcS}
\Psi=P_{\mathrm{c}}(\Sf(\Psi)) +[\Psi-P_{\mathrm{M}}]_+.
\end{align}
\end{remark}

\subsubsection{The Kirchhoff transform function}
\begin{figure}[h!]
\centering
\includegraphics[scale=.4]{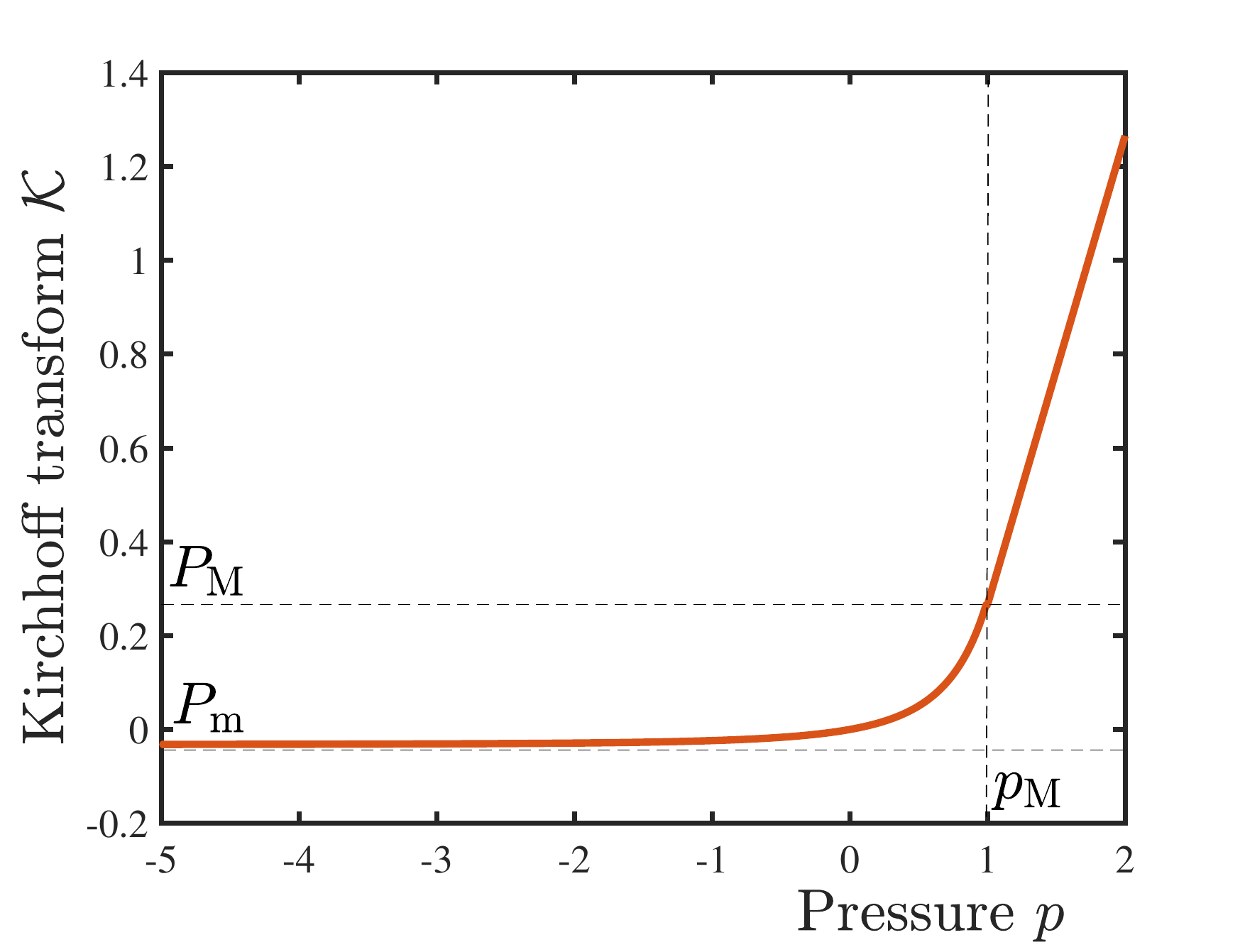}
\caption{The Kirchhoff transform $\calK$ for the Brooks--Corey model with $\lambda_1=0.75$.}\label{fig:Kirchhoff}
\end{figure}

The well-known Kirchhoff transformation \cite{alt1984nonstationary},  $\calK\in C^1(\R)$, is defined by
\begin{equation}\label{eq:TotalPandP}
\mathcal{K}(p):=\begin{cases}
P_\mathrm{c}(S(p))=\int_0^p \k(S(\vr))\,\mathrm{d}\vr &\text{ for } p\leq p_{\mathrm{M}},\\
P_\mathrm{M}+ \k(1)(p-p_\mathrm{M}) &\text{ for } p>p_\mathrm{M}.
\end{cases}
\end{equation}
The plot of $\calK$ is shown in \Cref{fig:Kirchhoff}. Note  from \ref{prop:pc} that $\calK(p)=P_{\mathrm{c}}(S(p))>P_{\mathrm{m}}$. This implies $\Sf\,\circ\,\calK=S$ since $\Sf(\calK(p))=P_{\mathrm{c}}^{-1}(P_{\mathrm{c}}(S(p)))=S(p)$ if $p\leq p_{\mathrm{M}}$, and $\Sf(\calK(p))= \Sf(P_{\mathrm{M}}+ \k(1)(p-p_{\mathrm{M}}))=1=S(p)$ if $p>p_{\mathrm{M}}$ (see \eqref{eq:DefSfunc}). Consequently, 
\begin{align}\label{eq:ConnectionVariables}
\text{taking } \Psi=\calK(p) \text{ there holds } \del \Psi=\k(S(p))\del p, \text{ and }  s=S(p)=\Sf(\Psi).
\end{align}

Explicit expressions of all the functions introduced above can be computed for the Brooks--Corey model. They are stated in \Cref{table:BrooksCorey}.
 
\begin{table}[h!]
\centering
\begin{tabular}{|l|l|l||l|l|l|}
\hline
func. &unit &Brooks--Corey expression &func. &unit &Brooks--Corey expression\\
 \hline
 &&&&&\\[-.8em]
 $\k(s)$ &-- & $s^{\frac{2+3\lambda_1}{\lambda_1}} $&  $S(p)$ &-- &$(2-\tfrac{p}{p_\mathrm{M}})^{-\lambda_1}$\\[.4em]
 $p_\mathrm{c}(s)$ &[Pa] & $ p_\mathrm{M}(2- s^{-\frac{1}{\lambda_1}})$ &$D(s)$ &[Pa] &$\tfrac{p_\mathrm{M}}{\lambda_1} S^{2+ \frac{1}{\lambda_1}}$\\[.4em]
 $P_\mathrm{c}(s)$ &[Pa] & $\tfrac{p_\mathrm{M}}{(1+3\lambda_1)}(s^{3+ \frac{1}{\lambda_1}}- 2^{-(1+3\lambda_1)})$ &$\Sf(\Psi)$ &-- &$[\tfrac{1+3\lambda_1}{p_\mathrm{M}} \Psi + 2^{-(1+3\lambda_1)}]^{\frac{\lambda_1}{1+3\lambda_1}}$
\\ \bottomrule   
\end{tabular}
\caption{The table of the introduced  functions $\k$, $S$, $p_{\mathrm{c}}$, $D$, $P_{\mathrm{c}}$, and $\Sf$ with their physical units and expressions for the Brooks--Corey model. The expressions are valid for $p\leq p_{\mathrm{M}}$, $s\in (0,1]$, and $\Psi\leq P_{\mathrm{M}}$. In addition, $\calK(p)=P_{\mathrm{c}}(S(p))$ for $p\leq p_{\mathrm{M}}$. In the heat equation case, for comparison, $\k(s)=D(s)=1$ and $S,\, p_\mathrm{c},\, P_\mathrm{c},\,  \Sf$, and $\calK$ are all identity functions.}\label{table:BrooksCorey}
\end{table}

\subsection{Weak formulations}\label{sec:Formulations}
We give below two equivalent weak formulations of the problem  \eqref{eq:FullSystem} discussing their strong and weak points. They will both be used to derive the a posteriori error estimates.
\subsubsection{The pressure formulation}
In the pressure formulation of \eqref{eq:FullSystem}, the main unknown is the pressure $p$. It reads: 
solve for $p \in \X$ and $s=S(p)\in H^1(0,T;H^{-1}(\Om))$ such that $s(0)=s_0$ and for all $\f\in \X$,
\begin{align}
\int^T_0 \langle \p_t S(p), \f \rangle +   \int^T_0  (\K \k(S(p)) (\del p + \vg),\del \f)=  \int^T_0 (f(S(p),\bm{x},t),\f).\label{eq:RichardsP}
\end{align}
This formulation has the advantage of generalizing to heterogeneous porous media, where the functions $S$ and $\k$ are defined differently in different subdomains of $\Om$. In particular, since $p$ is a physical quantity that remains continuous across the interfaces of such subdomains,  formulation \eqref{eq:RichardsP} has a conforming nature also in such circumstances.

\subsubsection{The total pressure formulation}
In the total pressure formulation of \eqref{eq:FullSystem}, the main unknown is the total pressure $\calK(p)$ which will henceforth be denoted by $\Psi$. It reads: solve for $\Psi\in \X$ with $s=\Sf(\Psi)\in H^1(0,T;H^{-1}(\Om))$ such that $s(0)=s_0$ and for all $\f\in \X$,
\begin{align}
\int^T_0 \langle \p_t \Sf(\Psi), \f \rangle +   \int^T_0  (\K (\del \Psi + \vg \k(\Sf(\Psi))),\del \f)=  \int^T_0 (f(\Sf(\Psi),\bm{x},t),\f).\label{eq:Richards}
\end{align}
The formulation \eqref{eq:Richards} is derived from \eqref{eq:RichardsP} using the variable transformation \eqref{eq:ConnectionVariables}. 
The total pressure formulation has the advantage of having a linear diffusion term. However, if the definition of $\k$ and $S$ varies inside the domain, for instance, in the case of heterogeneous porous media, then $\Psi$ is not uniformly defined. Moreover, the inverse transform $\Psi\mapsto p$ is often numerically expensive to compute, and $\Psi$ lacks a physical interpretation. We emphasize that, in this study, we have refrained from using $\mathcal{K}^{-1}$.

For $S(p)<1$, a saturation formulation is also valid, where $s=S(p)$ is the primary unknown and $D(s)$ serves as the diffusion coefficient. This formulation, however, breaks down at $s=1$ due to the non-invertibility of $S(p)$ \cite{alt1984nonstationary}.

\subsubsection{Well-posedness}
\begin{proposition}[Existence, uniqueness, and regularity]\label{pros:EUM}
Let \ref{prop:k}--\ref{prop:s0} hold.  Then there exists a unique weak solution $p\in \X$ of \eqref{eq:RichardsP} with $s=S(p)\in \W$  and $s(0)=s_0$. Moreover, there exists a unique weak solution $\Psi\in \X$ of \eqref{eq:Richards} with $\Sf(\Psi)\in \W$ and $\Sf(\Psi(0))=s_0$. Furthermore, the variables $p$, $s$, and $\Psi$  are related through \eqref{eq:ConnectionVariables}. 
\end{proposition}
The existence of a solution of \eqref{eq:RichardsP} for $p\in \X$ with $\p_t S(p)\in L^2(0,T;H^{-1}(\Om))$ has been proven in the seminal papers \cite{alt1983quasilinear,alt1984nonstationary}, whereas uniqueness is proven in \cite{otto1996l1} using $L^1$-contraction.  Since $p\in \X$, and $S(\cdot)$ is Lipschitz continuous, one automatically gets $s\in \W$. From the embedding of $\W$ in $C(0,T;L^2(\Om))$, we have $s\in C(0,T;L^2(\Om))$. The equivalence of the $p$ and the $\Psi$ formulations follows from the uniqueness of the solutions. 

\subsection{Maximum principle}\label{sec:MaxPrin}
In the case of the Richards equation, the saturation $s$ is bounded in $[0,1]$, and $s\searrow 0$ causes parabolic--hyperbolic degeneracy to occur. In this section, we use the maximum principle to obtain computable lower bounds for $s(\bm{x},t)$, bounding it away from $0$. 
For a positive initial saturation, the function $S_\mathrm{m}:\R^+\to (0,1]$ is a lower bound function of $s\in \W$, if 
\begin{equation}
0<S_\mathrm{m}(t)\leq s(\bm{x},t) \text{ for almost all } (\bm{x},t)\in \Om\times [0,T].\label{eq:MaximumPrin}
\end{equation}
To ensure that a lower bound function satisfying  \eqref{eq:MaximumPrin} exists for $S(p)\in \W$ when $p\in \X$ solves \eqref{eq:RichardsP}, additional  restrictions have to be imposed on the source term function $f$.
In particular, note that if $\k(0)=0$ and for some $(\bm{x},t)\in \Om\times [0,T]$ we have $s=0$, then from \eqref{eq:Richards1} $\p_t s=f(0,\bm{x},t)$. Since $s<0$ is unphysical, this forces 
\begin{equation}
f(0,\bm{x},t)\geq 0 \text{ for all } (\bm{x},t)\in \Om\times [0,T].\label{eq:fextra}
\end{equation}
This constraint will be imposed below to obtain computable maximum principle estimates.
 If $f$ is independent of $s$, then  \eqref{eq:fextra} simply implies that $f\geq 0$.
In comparison, in the context of the heat equation, $s$ is not bounded in $[0,1]$, and hence conditions such as \eqref{eq:fextra} are not required.

\subsubsection{A time-dependent lower bound}
Recalling hypothesis \ref{prop:f}, define a function $\bar{S}_\mathrm{m}(t)$ by the integral equation
\begin{align}\label{eq:DefSMSm}
& \bar{S}_\mathrm{m}(t)=\min\left (\mathrm{ess}\inf_{\bm{x}\in \Om}\{s_0(\bm{x})\},S(0)\right ) + \int_0^t f_{\mathrm{m}}(\bar{S}_\mathrm{m}(\vr))\, \mathrm{d}\vr .
\end{align}
Then, we have the following result:
\begin{proposition}[Existence of $\bar{S}_{\mathrm{m}}$ satisfying \eqref{eq:DefSMSm}]
Assume \ref{prop:f}--\ref{prop:s0}.  Additionally, assume that there exists a choice of the function $f_{\mathrm{m}}$ such that an interval $[0,J]$, $J\in (0,1),$ and a constant $C_f\geq 0$ exist for which the inequality
$$
f_{\mathrm{m}}(s)\geq -C_f s \text{ holds } \forall s\in [0,J].
$$
Then, there exists a continuous function $\bar{S}_{\mathrm{m}}:\R^+\to \R^+$ that satisfies \eqref{eq:DefSMSm}. 
\end{proposition}

The existence of $\bar{S}_{\mathrm{m}}$ follows from the Picard--Lidel\"of theorem by the differentiability of the function $f_{\mathrm{m}}$ assumed in \ref{prop:f}. The bound $\bar{S}_{\mathrm{m}}(t)>0$ follows from the inequality $\tfrac{d}{dt} \bar{S}_{\mathrm{m}}(t)\geq -C_f \bar{S}_{\mathrm{m}}(t)$ and $\bar{S}_{\mathrm{m}}(0)>0$.  The constraint $f\geq f_{\mathrm{m}}\geq -C_f s$ embodies and generalises \eqref{eq:fextra}. In practice, $\bar{S}_{\mathrm{m}}(t)$ can be computed to arbitrary precision using numerical approaches such as the Runge--Kutta method.

\begin{proposition}[Time-dependent lower bound of $s$]\label{pros:MaximumPrinciple}
Let \ref{prop:f}--\ref{prop:s0} hold and $p\in \X$ with $s=S(p)\in \W$ and $s(0)=s_0$ be a solution of \eqref{eq:RichardsP}. Moreover, let $\K$ be  constant in $\Om$. Then  $S_\mathrm{m}=\min(\bar{S}_{\mathrm{m}},S(0))$, with $\bar{S}_{\mathrm{m}}$ defined in \eqref{eq:DefSMSm}, is a lower bound function of $s$ satisfying \eqref{eq:MaximumPrin}.
\end{proposition}

Since proving the maximum principle result is not the main focus of this paper, we postpone the proof to Appendix \ref{App:proof}, along with other proofs of this section.

\subsubsection{A space-dependent lower bound}
\Cref{pros:MaximumPrinciple} gives a computable lower bound of $s(\bm{x},t)$ for a given $t\in [0,T]$, provided the absolute permeability $\K$ is constant. The following result also gives a lower bound of $s$, relaxing the assumption of $\K$ being constant.
\begin{proposition}[Existence of a bounded function]\label{pros:IntroVS}
Let \ref{prop:k}--\ref{prop:f} hold. For a constant $J\leq 0$, let $\vs\in H^1(\Om)$ with $\vs=J$ on $\p\Om$ in the trace sense, solve 
\begin{align}
(\K \k(S(\vs))[\del \vs +\vg],\del \f)=\left (\inf_{t\in \R^+} [f(S(\vs),\bm{x},t)]_-,\f \right ),\quad \forall \f\in H^1_0(\Om).\label{eq:VSweak}
\end{align}
Assume that there exists a constant $p_{\mathrm{l}}\leq 0$ such that $f(S(p),\bm{x},t)\geq 0$ for all $p<p_{\mathrm{l}}$. Then 
\begin{align}
\min(p_\mathrm{l},J) + \min_{\bm{x}\in \Om}\{\vg\cdot\bm{x}\}\leq \vs(\bm{x})+  \vg\cdot\bm{x}\leq J +  \max_{\bm{x}\in \Om}\{\vg\cdot\bm{x}\}  \text{ for almost all } \bm{x}\in \Om.
\end{align}
\end{proposition}
\noindent
The existence of $\vs$ follows from \cite{alt1983quasilinear} and the existence of $p_{\mathrm{l}}<0$ is compatible with \eqref{eq:fextra}. The counterpart of \Cref{pros:MaximumPrinciple} for this case is:
\begin{proposition}[Space-dependent lower bound of $s$]\label{pros:MaximumPrincipleB}
Let \ref{prop:f}--\ref{prop:s0} hold and $p\in \X$ with $s=S(p)\in \W$ and $s(0)=s_0$ be a solution of \eqref{eq:RichardsP}. For the constant
$$
J=\mathrm{ess}\inf_{\bm{x}\in \Om}\left ([p_{\mathrm{c}}(s_0(\bm{x}))]_- - \max_{\bm{x}\in \Om}\{\vg\cdot\bm{x}\}+ \vg\cdot\bm{x}\right )\leq 0,
$$ 
let $\vs\in H^1(\Om)$ be obtained from \Cref{pros:IntroVS}. Then $S_\mathrm{m}(t)= \mathrm{ess}\inf\limits_{\bm{x}\in \Om}(S(\vs(\bm{x})))$ for $t>0$  is a lower bound function of $s$ satisfying \eqref{eq:MaximumPrin}.
\end{proposition}

\section{Relations between the error and the residual}
In this section, the dual norm of the residual will be used to bound from above and from below an error metric that we will use  in place of the $\W$-norm in the present nonlinear and degenerate setting. 
\label{sec:ErrorResidual}

\subsection{Residual}
For
\begin{align}
\Psi_{h\t}\in \X \text{ with } s_{h\t}:=\Sf(\Psi_{h\t})\in \W,
\end{align} 
the residual $\mathcal{R}(\Psi_{h\t})\in L^2(0,T;H^{-1}(\Om))$  with respect to the weak formulation \eqref{eq:Richards} is defined as
\begin{equation}
\int_0^T \langle\mathcal{R}(\Psi_{h\t}),\f\rangle:= \int_0^T \left [(f(s_{h\t},\bm{x},t),\f)- \langle \p_t s_{h\t}, \f\rangle -(\K  (\del \Psi_{h\t} + \vg\k(s_{h\t})),\del \f) \right ]\label{eq:DefRhtau}
\end{equation}
for all $\f\in \X$. If $\Psi\in \X$ with $s=\Sf(\Psi)\in \W$ denotes the solution to \eqref{eq:Richards} then $\mathcal{R}(\Psi)=0$.

\subsection{Norms}
On a Lipschitz subdomain $\om\subseteq \Om$, we introduce 
equivalent (semi)norms on $H^1(\om)$, $H^1_0(\om)$ and $H^{-1}(\om)$: 
\begin{subequations}\label{eq:DualNormDef}
\begin{align}
\label{eq:NormHminK}
 &\|\vr\|_{\Kh(\om)}:=\|\K^{\frac{1}{2}}\del \vr\|_{\om} \text{ for } \vr\in H^{1}(\om),\\
 &\|\vr\|_{\Kn(\om)}:=\sup_{\f\in H^1_0(\om)} \{\langle \vr,\f\rangle_{H^{-1}(\om),H^1_0(\om)} \slash \|\f\|_{\Kh(\om)}\} \text{ for } \vr\in H^{-1}(\om).
\end{align}
\end{subequations}
From the properties of $\K$ stated in \ref{prop:Kabs}, it is immediate that
\begin{align}\label{eq:CurrencyKnorms}
K_\mathrm{m}^{\frac{1}{2}} \|\del \vr\|_{\om}\leq \|\vr\|_{\Kh(\om)}\leq K_\mathrm{M}^{\frac{1}{2}} \|\del \vr\|_{\om}\; \text{ and }\; K_\mathrm{M}^{-\frac{1}{2}} \|\vr\|_{H^{-1}(\om)}\leq \|\vr\|_{\Kn(\om)}\leq K_\mathrm{m}^{-\frac{1}{2}} \|\vr\|_{H^{-1}(\om)}.
\end{align}

Let $\a:[0,T]\to [0,\infty)$ denote a bounded non-negative function. For a subdomain $\om\subseteq \Om$, and an interval $I\subseteq [0,T]$, we introduce the distance measure  $\dst^{\a}_{\om,I}$ on the set $\{\psi\in L^2(0,T;H^1(\om)): \Sf(\psi)\in H^1(0,T;H^{-1}(\om))\}$  as
\begin{align}\label{eq:DefLBErrorMeasure}
\dst^{\a}_{\om,I}(\Psi_1,\Psi_2):=&\|\p_t(\Sf(\Psi_1)-\Sf(\Psi_2))\|_{L^2(I;\Kn(\om))}\nonumber\\
&  + \|\a(\Sf(\Psi_1)-\Sf(\Psi_2))\|_{L^2(\om\times I)} + \|\Psi_1-\Psi_2\|_{L^2(I,\Kh(\om))}.
\end{align}
The distance measure combines the $L^2(I;\Kh(\om))$-norm of $\Psi_1-\Psi_2$ with the $H^1(I;\Kn(\om))\cap L^2(\om\times I)$ norms of $\Sf(\Psi_1)-\Sf(\Psi_2)$. Note that for $\a=0$, the middle term disappears.

We also introduce the class of time-integration functionals $\J_\a:L^2([0,T])\to [0,\infty)$ as: for $\vr\in L^2([0,T])$,
 \begin{align}\label{eq:TiemIntegrator}
 \J_\a(\vr):= \left [\exp\left (-\smallint_0^T \a \right )\int_0^T \left (\vr^2(t) + \a(t) \exp\left (\smallint_t^T \a \right ) \int_0^t \vr^2 \right)\dd t \right ]^{\frac{1}{2}}.
 \end{align}
 The operator $\J_\a$ defines a norm and satisfies the triangle inequality. It is actually equivalent to the $L^2([0,T])$-norm, since the inequality $0\leq \int_0^t \vr^2 \dd t \leq \|\vr \|^2_{L^2([0,T])}$ for $t\in [0,T]$, and $\int_0^T \a  \exp(\int_t^T \a)\dd t= \exp(\int_0^T \a)-1$  directly gives 
\begin{align}\label{eq:L2normEquivalent}
  \exp\left (-\tfrac{1}{2}\smallint_0^T \a \right ) \|\vr \|_{L^2([0,T])} \leq \J_\a(\vr)\leq \|\vr \|_{L^2([0,T])}.
 \end{align}
Consequently, it is equal to the $L^2([0,T])$ norm if $\a= 0$.

\subsection{Lower bound on the error by the residual}
Extending Theorem 2.1 of \cite{Ern_Sme_Voh_heat_HO_Y_17} to the present degenerate nonlinear setting, we have
\begin{theorem}[Lower bound on error by the dual norm of the residual]\label{theo:TotalPressureEfficiency}
Let \ref{prop:k}--\ref{prop:s0} hold and let $\Psi\in \X$ with $s=\Sf(\Psi)\in \W$ denote the unique solution of \eqref{eq:Richards}. Let $\om\subseteq \Om$ be a Lipschitz subdomain of $\Om$ and let $I\subseteq [0,T]$ be a time interval. Let the norms $\|\cdot\|_{\Kh}$, $\|\cdot\|_{\Kn}$ and the error measure $\dst^{\a}_{\om,I}(\cdot,\cdot)$ be defined as in \eqref{eq:DualNormDef}--\eqref{eq:DefLBErrorMeasure} for $\a(t)=C_{\mathrm{P},\om}\,h_{\om}\, K_{\mathrm{m}}^{-\frac{1}{2}}\,\max\limits_{[0,1]\times\om\times \{t\}}|\p_s f|  + |\vg| K_{\mathrm{M}}^{\frac{1}{2}} \|\k'\|_{L^\infty([0,1])}$. Then, for any $\Psi_{h\t}\in \X$ with $s_{h\t}= \Sf(\Psi_{h\t})\in \W$,  one has
\begin{align}
\|\mathcal{R}(\Psi_{h\t})\|_{L^2(I;\Kn(\om))}\leq \dst^{\a}_{\om,I}(\Psi,\Psi_{h\t}).\label{eq:EfiiciencyMain1a}
\end{align}
\end{theorem}

\begin{remark}[The linear case]
Observe that in the linear case, $\k= 1$ and $\p_s f= 0$, yielding $\a= 0$.
\end{remark}
\begin{proof} 
From \eqref{eq:DefRhtau}, one has for any $\f\in L^2(I;H^1_0(\om))$, extended to $\Om\setminus \om$ and $[0,T]\setminus I$ by 0, that
\begin{align}
\int_I \langle\mathcal{R}(\Psi_{h\t}),\f\rangle= \int_I &\left [  \langle \p_t (s- s_{h\t}), \f\rangle +(\K\del (\Psi -\Psi_{h\t}),\del \f) \right. \nonumber\\
&\left.+ (f(s_{h\t},\bm{x},t)-f(s,\bm{x},t),\f)+  (\K\vg(\k(s)-\k(s_{h\t})),\del \f) \right ].\label{eq:DefRrewritten}
\end{align}
Then, from the triangle inequality and definitions of the norms $\|\cdot\|_{\Kh}$, $\|\cdot\|_{\Kn}$ we get
\begin{align*}
&\|\mathcal{R}(\Psi_{h\t})\|_{L^2(I;\Kn(\om))}  \leq\| \p_t (s- s_{h\t})\|_{L^2(I;\Kn(\om))} + \|\Psi -\Psi_{h\t}\|_{L^2(I;\Kh(\om))}\\
& \qquad +\sup_{\substack{\scriptscriptstyle{\f\in L^2(I;H^1_0(\om))},\\ \scriptscriptstyle{\|\f\|_{L^2(I;\Kh(\om))}=1}}}\int_I \left [(f(s_{h\t},\bm{x},t)-f(s,\bm{x},t),\f)_\om+  (\K\vg(\k(s)-\k(s_{h\t})),\del \f)_\om\right ]. 
\end{align*}
The result then follows from the definition of $\dst^{\a}_{\om,I}$ and the computation of the last two terms using 
\begin{align*}
&|(f(s_{h\t},\bm{x},t)-f(s,\bm{x},t),\f)_\om|\overset{\eqref{eq:Poincare},\, \eqref{eq:CurrencyKnorms}}\leq  C_{\mathrm{P},\om}\,h_{\om}\,K_{\mathrm{m}}^{-\frac{1}{2}}\,\max\limits_{[0,1]\times\om\times \{t\}}|\p_s f| \|s_{h\t}-s\|_{\om}\|\f\|_{\Kh(\om)},\\
& |(\K\vg(\k(s)-\k(s_{h\t})),\del \f)_{\om}|\overset{\ref{prop:k},\,\ref{prop:Kabs}}\leq |\vg| K_{\mathrm{M}}^{\frac{1}{2}} \|\k'\|_{L^\infty([0,1])}\|s_{h\t}-s\|_{\om}\|\f\|_{\Kh(\om)}.
\end{align*}
\end{proof}
\vspace{-2em}

\subsection{Upper bound on the error by the residual}
For the lower bound function $S_m(t)$ satisfying  \eqref{eq:MaximumPrin}, the diffusivity function $D$ of \eqref{eq:DefDiffusivity}, the saturation function $\Sf$ of \eqref{eq:DefSfunc}, and the source term $f$ of \ref{prop:f}, let
\begin{subequations}
\begin{align}
&D_{\mathrm{m}}(t):=\min\{D(\vr): \vr\in [S_\mathrm{m}(t),1]\},\quad D_\mathrm{M}(t) :=\max\{|D'(\vr)|: \vr\in [S_\mathrm{m}(t),1]\},\\
&\Sf_{\p,\mathrm{M}}(t):=\max\{\Sf'(P_\mathrm{c}(\vr)): \vr\in [S_\mathrm{m}(t),1]\},\\
&f_{\p,M}(t):=\max\{|\p_s f(\vr,\bm{x},t)|: \vr\in [0,1], \, (\bm{x},t)\in \Om\times [0,T]\}.
\end{align}\label{eq:DmFpm}
\end{subequations}
Recalling \eqref{eq:propD}, we have $D_\mathrm{m}(t)>0$ and $D_\mathrm{M}(t)<\infty$. Similarly $\Sf_{\p,\mathrm{M}}(t)<\infty$, $f_{\p,M}(t)<\infty$.
Then, inspired by \cite{Di_Pi_Voh_Yous_a_post_Stef_15} we propose 
\begin{theorem}[Upper bound on error by the dual norm of the residual]\label{theo:UpperBound}
Let \ref{prop:k}--\ref{prop:s0} hold and $\Psi\in \X$ denote the unique solution of \eqref{eq:Richards} with $s=\Sf(\Psi)\in \W$.  Let $\Psi_{h\t}\in \X$ with $s_{h\t}=\Sf(\Psi_{h\t})\in \W$ be arbitrary.
 Assume that a lower bound function $S_{\mathrm{m}}(t)$, satisfying \eqref{eq:MaximumPrin}, exists for $s$ and $s_{h\t}$.
 Recall the definitions of $D_\mathrm{m}$, $D_{\mathrm{M}}$, $f_{\p,M}$, and $\Sf_{\p,\mathrm{M}}$ from \eqref{eq:DmFpm}. 
 Let the residual $\mathcal{R}$,  norms $\|\cdot\|_{H^{\pm 1}_\K}$, and the time-integrator $\J_\a$ be defined in \eqref{eq:DefRhtau}, \eqref{eq:DualNormDef}, and \eqref{eq:TiemIntegrator}  respectively. Then, for any $\lambda:[0,T]\to \R^+$, the following estimates hold:
\begin{subequations}\label{eq:UpperBoundEstimate}
\textbf{Estimate in the } ${L^2(\Om\times [0,T])}$ \textbf{ and } ${L^\infty(0,T;\Kn(\Om))}$ \textbf{ norms:}
\begin{align}
& { e^{-\smallint_0^T (\lambda + \mathfrak{C}_1)} }\|(s-s_{h\t})(T)\|^2_{\Kn(\Om)}+  \J_{\lambda + \mathfrak{C}_1} \left (\Sf_{\p,\mathrm{M}}^{-\frac{1}{2}}\|s-s_{h\t}\| \right )^2\nonumber\\
 \leq &\|s_0-s_{h\t}(0)\|^2_{\Kn(\Om)} +  \J_{\lambda + \mathfrak{C}_1}(\lambda^{-\frac{1}{2}}\|\mathcal{R}(\Psi_{h\t})\|_{\Kn(\Om)})^2.\label{eq:UpperBoundEstimate1st}
\end{align}
\hspace{1.3em}\textbf{Estimate in the } $L^2(0,T;\Kh(\Om))$  \textbf{and } ${L^\infty(0,T;L^2(\Om))}$ \textbf{ norms:}  For \vspace{-.5em}
\begin{align*}
C^\infty_{h\t}(t):=\|\K^{\frac{1}{2}}\del s_{h\t}(t)\|^2_{L^\infty(\Om)},\, \text{ assume that }\smallint_0^T C^\infty_{h\t}(t)\,\mathrm{d}t <\infty.
\end{align*}  
On a Lipschitz subdomain $\Om^{\mathrm{deg}}(t)\supseteq \{s(\bm{x},t)=1\}\cup \{s_{h\t}(\bm{x},t)=1\}$ of  $\Om$ (possibly disconnected), let $D(s)/2 \leq D(s_{h\t})\leq 2 D(s)$ hold, and define the parabolic--elliptic degeneracy estimator $\eta^{\mathrm{deg}}\in L^2([0,T])$ as
\begin{align*}
\eta^{\mathrm{deg}}(t):=& \sqrt{\tfrac{2}{D(1)} }\bigg [\|[\Psi_{h\t}(t)-P_{\mathrm{M}}]_+\|^2_{\Kh(\Om)} \nonumber\\
&  + \left (\|[f(1,\bm{x},t)]_+\|_{\Kn(\Om^{\mathrm{deg}}(t))} +\left  \|\left (\K^{\frac{1}{2}}-\tfrac{\K^{-\frac{1}{2}}}{|\Om^{\mathrm{deg}}(t)|}\smallint_{\Om^{\mathrm{deg}}(t)} \K \right )\bm{g}\right \|_{\Om^{\mathrm{deg}}(t)}\right )^2  \bigg]^{\frac{1}{2}}.
\end{align*}
Then it holds that,
\begin{align}\label{eq:UpperBoundEstimate2nd}
&e^{-\smallint_0^T \mathfrak{C}_2 }\|(s-s_{h\t})(T)\|^2 + \tfrac{1}{2}\J_{\mathfrak{C}_2}\left (\left\| D(s)^{-\frac{1}{2}} \K^{\frac{1}{2}}\del (\Psi-\Psi_{h\t})\right \| \right )^2 \nonumber\\
  \leq  &\|s_0-s_{h\t}(0)\|^2 + \J_{\mathfrak{C}_2}\left ( \eta^{\mathrm{deg}} \right )^2 + 4  \, \J_{\mathfrak{C}_2}\left ( D_{\mathrm{m}}^{-\frac{1}{2}} \|\mathcal{R}(\Psi_{h\t})\|_{\Kn(\Om)} \right )^2.
\end{align}
\hspace{1.3em}\textbf{Estimate in the } $H^1(0,T;\Kn(\Om))$ \textbf{norm:}
\begin{align}\label{eq:UpperBoundEstimate3rd}
&\J_\lambda(\|\p_t (s-s_{h\t})\|_{\Kn(\Om)})^2\nonumber\\
 \leq\; & 3\left [\J_\lambda (\|\Psi-\Psi_{h\t}\|_{\Kh(\Om)} )^2 + \mathfrak{C}_3(T)\, \J_\lambda\left (\|s-s_{h\t}\|\right )^2
 + \J_\lambda (\|\mathcal{R}(\Psi_{h\t})\|_{\Kn(\Om)})^2 \right ].
\end{align}
\end{subequations}
Recalling the Poincar\'e constant $C_{\mathrm{P},\Om}$ from \eqref{eq:Poincare}, the functions $\mathfrak{C}_{1,2,3}:(0,T)\to [0,\infty)$ are 
\begin{subequations}\label{eq:Cestimators}
\begin{align}
&\mathfrak{C}_1(t):= 2 \Sf_{\p,\mathrm{M}} (t)\left[ K_\mathrm{M}|\bm{g}|^2  \|\k'\|_{L^\infty([0,1])}^2 + \tfrac{C^2_{\mathrm{P},\Om} h_\Om^2}{K_\mathrm{m} }  f_{\p,\mathrm{M}}^2(t)\right ],\\[.2em]
&\mathfrak{C}_2(t):= \tfrac{1}{D_{\mathrm{m}}(t)} \left[ D^2_{\mathrm{M}}(t)\, C^\infty_{h\t}(t) +  4 K_{\mathrm{M}} |\bm{g}|^2 \|\k'\|^2_{L^\infty([0,1])}\right ]+  2 f_{\p,\mathrm{M}}(t),\\[.2em]
&\mathfrak{C}_3(t):=(C_{\mathrm{P},\Om} h_\Om\, K^{-\frac{1}{2}}_\mathrm{m}\, \|f_{\p,\mathrm{M}}\|_{L^\infty([0,t])} + K^{\frac{1}{2}}_\mathrm{M} |\bm{g}| \|\k'\|_{L^\infty([0,1])})^2.
\end{align}
\end{subequations}
\end{theorem}

The function $\lambda>0$ in \eqref{eq:UpperBoundEstimate1st} is introduced to optimize the effectivity of the estimates. The reason as well as a possible value of $\lambda$ will be explained in detail in \Cref{rem:lambda}.

\begin{remark}[Degeneracy at $s=1$]
Observe that the estimate \eqref{eq:UpperBoundEstimate2nd} contains the degenereacy estimator $\etaDeg$, despite the estimates \eqref{eq:UpperBoundEstimate1st}, \eqref{eq:UpperBoundEstimate3rd} not including it. This stems from the fact that proving a contraction in $L^2(0,T;H^1(\Om))$ is generally not possible for degenerate problems. However, proving contraction in the  $L^\infty(0,T;L^1(\Om))$ and the $L^2(0,T;H^{-1}(\Om))$ norms are possible  \cite{otto1996l1,kubo2005nonlinear}. The last two components in the definition of $\etaDeg$ represent the two reasons why the parabolic--elliptic degeneracy might occur despite the initial condition $s_0$ being in $(0,1]$, i.e. the positivity of $f$ and the non-uniformity of $\K$.  Additionally, assuming that ${s=1}$ only occurs  in a superset covering $\{s_{h\t}=1\}$, the estimator $\etaDeg$ is fully computable, see \Cref{sec:Numerical}. 
\end{remark}

\begin{remark}[Reduction in the linear case]
Observe that, in the linear heat equation case, $\k(s)$ and $p_{\mathrm{c}}'(s)$ are equal to 1, giving a constant $D(s)=1$ and $D_{\mathrm{M}}(t)=0$. Similarly $\p_s f=0$. Thus, one has $\mathfrak{C}_{1,2,3}=0$. Hence, for the linear case, taking $\lambda=0$ in \eqref{eq:UpperBoundEstimate3rd} exponential terms in \eqref{eq:UpperBoundEstimate2nd}--\eqref{eq:UpperBoundEstimate3rd} vanish,  and they reduce to the estimates provided in \cite{Ern_Sme_Voh_heat_HO_Y_17,Ver_a_post_heat_03}.
\end{remark}

\begin{remark}[Bounds on $\|\p_t(s-s_{h\t})\|_{\Kn(\Om)}$ and $\dst^{\a}_{\Om,[0,T]}(\Psi,\Psi_{h\t})$]
Choosing $\lambda$ in \eqref{eq:UpperBoundEstimate1st} such that $\lambda + \mathfrak{C}_1=\mathfrak{C}_2$ and $\lambda=\mathfrak{C}_2$ in \eqref{eq:UpperBoundEstimate3rd}, we have a complete bound for $\|\p_t(s-s_{h\t})\|_{\Kn(\Om)}$ using the other components of \eqref{eq:UpperBoundEstimate}. Combining \eqref{eq:UpperBoundEstimate}, one obtains an estimate for all components of $\dst^{\a}_{\Om,[0,T]}(\Psi,\Psi_{h\t})$ defined in \eqref{eq:DefLBErrorMeasure}. Hence, \Cref{theo:TotalPressureEfficiency,theo:UpperBound} provide both  lower and upper bounds of $\dst^{\a}_{\Om,[0,T]}(\Psi,\Psi_{h\t})$ in that using \eqref{eq:L2normEquivalent}, one has
\begin{align*}
\int^T_0\|\mathcal{R}(\Psi_{h\t})\|^2_{\Kn(\Om)}&\overset{\eqref{eq:EfiiciencyMain1a}}\leq \dst^{\a}_{\Om,[0,T]}(\Psi,\Psi_{h\t})^2\\
&\overset{\eqref{eq:UpperBoundEstimate}}\lesssim \exp\left (\smallint_0^T \mathfrak{C}_2\right )\left [ \|s_0-s_{h\t}(0)\|^2 +\int_0^T \left ([\eta^{\mathrm{deg}}]^2 + \|\mathcal{R}(\Psi_{h\t})\|^2_{\Kn(\Om)}\right )^2 \right ].
\end{align*}
However, this upper bound is rather rough since it hides its dependence on $D_\mathrm{m}$, $\mathfrak{C}_{1\slash 3}$ and $C_{\mathrm{P},\Om}\, h_\Om$. Note that $\exp\left (\smallint_0^T \mathfrak{C}_2\right )$ may take very large values and might explode as $T\to \infty$, which is the usual consequence of using Gronwall Lemma. This is avoided in our analysis.
\end{remark}

\begin{proof}[\textbf{Proof of \Cref{theo:UpperBound}}]
In the proof, we shorten $\mathcal{R}(\Psi_{h\t})$ to simply $\mathcal{R}$. From \eqref{eq:DefRhtau}, we have for all $\f\in L^2(0,T;H^1_0(\Om))$,
\begin{align}
&\int_0^T \left [\langle \p_t (s-s_{h\t}),\f \rangle + (\K  \del (\Psi-\Psi_{h\t}), \del \f)\right ]\nonumber\\
&\quad = \int_0^T \left [\langle \mathcal{R}, \f\rangle  + (f(s,\bm{x},t)-f(s_{h\t},\bm{x},t)),\f) + (\K \vg(\k(s_{h\t})-\k(s)),\del \f)\right ].\label{eq:RobustnessEqAll}
\end{align}

\textbf{Step 1 (Estimate \eqref{eq:UpperBoundEstimate1st}):}
Let the Green function $G^0_{h\t}\in C(0,T; H^1_0(\Om))$ satisfy for all $t\in [0,T]$ and $\f\in H^1_0(\Om)$,
\begin{equation}\label{eq:GreenFuncUpperB}
(\K \del G^0_{h\t}(t), \del \f)=\langle(s-s_{h\t})(t),\f\rangle.
\end{equation}
The problem is well-defined as $(s-s_{h\t})(t)\in L^2(\Om)$. Moreover, 
\begin{align}\label{eq:HKminG0ht}
\| G^0_{h\t}(t)\|_{\Kh(\Om)}&=\sup_{\|\f\|_{\Kh(\Om)}=1}(\K\del G^0_{h\t}(t),\del \f)\nonumber\\
&=\sup_{\|\f\|_{\Kh(\Om)}=1} \langle (s-s_{h\t})(t), \f\rangle=\|(s-s_{h\t})(t)\|_{\Kn(\Om)}.
\end{align}
Since $\p_t (s-s_{h\t})\in L^2(0,T;H^{-1}(\Om))$, equation \eqref{eq:GreenFuncUpperB} can be differentiated in time, implying that $\p_t G^0_{h\t} \in \X$ exists satisfying
\begin{equation}\label{eq:ptG0ht}
\int_0^T (\K \del \p_t G^0_{h\t}, \del \f)=\int_0^T  \langle \p_t (s-s_{h\t}),\f\rangle \text{ for all } \f\in \X.
\end{equation}
We now insert the test function $\f= G^0_{h\t}$ in \eqref{eq:RobustnessEqAll}. Using \eqref{eq:ptG0ht}, we see
\begin{align}
&\int_0^T \langle \p_t (s-s_{h\t}),G^0_{h\t}\rangle=\int_0^T (\K \del \p_t G^0_{h\t},\del G^0_{h\t})= \tfrac{1}{2}\int_{\Om}\left [ |\K^{\frac{1}{2}}\del G^0_{h\t} (T)|^2 -|\K^{\frac{1}{2}}\del G^0_{h\t} (0)|^2\right ]\nonumber\\
&\overset{\eqref{eq:HKminG0ht}}{=} \tfrac{1}{2}\|G^0_{h\t} (T)\|^2_{\Kh(\Om)} - \tfrac{1}{2}\|s_0-s_{h\t}(0)\|^2_{\Kn(\Om)}\label{eq:GreensMainIneq1st}.
\end{align}
Using the identity \eqref{eq:PsiandPcS} and noting that $([\Psi-P_{\mathrm{M}}]_+-[\Psi_{h\t}-P_{\mathrm{M}}]_+, s -s_{h\t})\geq 0$ which follows from the monotonicity of $[\cdot]_+$, one further has from \eqref{eq:GreenFuncUpperB} that
\begin{align}
& \int_0^T (\K \del (\Psi- \Psi_{h\t}),\del G^0_{h\t})=\int_0^T (\Psi- \Psi_{h\t},s-s_{h\t})\overset{\eqref{eq:PsiandPcS}}\geq \int_0^T (P_{\mathrm{c}}(s)- P_{\mathrm{c}}(s_{h\t}),s-s_{h\t}) \nonumber\\
&=\int_0^T\int_\Om {P_{\mathrm{c}}}'(s-s_{h\t})^2 \overset{\eqref{eq:profSf},\, \eqref{eq:DmFpm}}\geq \int_0^T \frac{1}{\Sf_{\p,M}(t)}\|s-s_{h\t}\|^2.
\end{align}
Recalling the Poincar\'e inequality \eqref{eq:Poincare} and the definitions  \eqref{eq:DualNormDef} of $\|\cdot\|_{\Kh}$, $\|\cdot\|_{\Kn}$ norms, we have
\begin{align}
& \int_0^T \langle \mathcal{R},  G^0_{h\t}\rangle\leq \int_0^T \|\mathcal{R}\|_{\Kn(\Om)}\;\|  G^0_{h\t} \|_{\Kh(\Om)} \leq   \int_0^T \left [\tfrac{1}{2\lambda}\|\mathcal{R}\|^2_{\Kn(\Om)} +  \tfrac{\lambda}{2} \|  G^0_{h\t} \|^2_{\Kh(\Om)} \right ],\label{eq:YoungOnResLambda}
\end{align}
as well as
\begin{align}
& \int_0^T (f(s,\bm{x},t)-f(s_{h\t},\bm{x},t),G^0_{h\t})\overset{\eqref{eq:DmFpm}}\leq \int_0^T   f_{\p,\mathrm{M}}(t) \|s-s_{h\t}\|\|G^0_{h\t}\|\nonumber\\
&\quad \leq \tfrac{1}{4}   \int_0^T \frac{1}{\Sf_{\p,\mathrm{M}}(t)}\|s-s_{h\t}\|^2 +  \int_0^T  \Sf_{\p,\mathrm{M}}(t)\, f_{\p,\mathrm{M}}(t)^2 \, \|G^0_{h\t}\|^2 \nonumber\\
&\quad \overset{\eqref{eq:Poincare},\, \eqref{eq:CurrencyKnorms}}\leq \tfrac{1}{4}   \int_0^T \frac{1}{\Sf_{\p,\mathrm{M}}(t)}\|s-s_{h\t}\|^2 +\tfrac{C^2_{\mathrm{P},\Om} h_\Om^2}{K_\mathrm{m}}\int_0^T \Sf_{\p,\mathrm{M}}(t) \, f_{\p,\mathrm{M}}(t)^2 \, \|  G^0_{h\t} \|^2_{\Kh(\Om)},
\end{align}
and
\begin{align}
& \int_0^T (\K \vg(\k(s_{h\t})-\k(s)),\del G^0_{h\t} ) \leq \tfrac{1}{4 K_\mathrm{M}|\vg|^2 \|\k'\|_{L^\infty([0,1])}^2}  \int_0^T \tfrac{1}{\Sf_{\p,\mathrm{M}}(t)}  \int_{\Om} \vg^{\mathrm{T}}\K \vg(\k(s)-\k(s_{h\t}))^2 \nonumber\\
&\qquad \qquad \phantom{abcdefghijklmnopqrs}+ {\textstyle K_\mathrm{M}|\vg|^2 \|\k'\|_{L^\infty([0,1])}^2}  \int_0^T  \Sf_{\p,\mathrm{M}}(t) \int_{\Om} |\K^{\frac{1}{2}} \del G^0_{h\t}|^2\nonumber\\
&\quad \overset{\ref{prop:Kabs}}\leq \tfrac{1}{4}   \int_0^T \frac{1}{\Sf_{\p,\mathrm{M}}(t)}\|s-s_{h\t}\|^2  + {\small K_\mathrm{M}|\vg|^2 \|\k'\|_{L^\infty([0,1])}^2} \int_0^T  \Sf_{\p,\mathrm{M}}(t)\, \|  G^0_{h\t} \|^2_{\Kh(\Om)}.\label{eq:GreensMainIneqEnd}
\end{align}
Combining \eqref{eq:GreensMainIneq1st}--\eqref{eq:GreensMainIneqEnd} with \eqref{eq:RobustnessEqAll}, one has
\begin{align}
& \|  G^0_{h\t}(T) \|^2_{\Kh(\Om)} +  \int_0^T  \frac{1}{\Sf_{\p,\mathrm{M}}(t)} \|(s-s_{h\t})(t)\|^2 \nonumber\\
& \leq \|s_0-s_{h\t}(0)\|^2_{\Kn(\Om)} +  \int_0^T\tfrac{1}{\lambda}\|\mathcal{R}\|^2_{\Kn(\Om)} + \int_0^T (\lambda+ \mathfrak{C}_1(t))\, \|  G^0_{h\t}(t) \|^2_{\Kh(\Om)}.
\end{align}
Applying the Gronwall Lemma
\begin{align}
u(t)\leq \a(t) + \int_0^t \b(\vr)u(\vr)\mathrm{d}\vr \implies u(t)\leq \a(t) + \int_0^t \b(\vr)\, \a(\vr) \exp\left( \smallint_\vr^t  \b(r) \, \mathrm{d}r\right) \mathrm{d}\vr\label{eq:GronwallLemma}
\end{align}
 with $u(t)= \|  G^0_{h\t} (t)\|^2_{\Kh(\Om)}$, $\a(t)=\|s_0-s_{h\t}(0)\|^2_{\Kn(\Om)} + \int_0^t \lambda^{-1} \|\mathcal{R}\|^2_{\Kn(\Om)} -  \int_0^t \frac{1}{\Sf_{\p,\mathrm{M}}} \|(s-s_{h\t})\|^2$, $\b(t)=\lambda+\mathfrak{C}_1(t)$, and re-normalizing both sides by dividing with $\exp(\int_0^T(\lambda+ \mathfrak{C}_1))$, we have \eqref{eq:UpperBoundEstimate1st}. Observe that the total coefficient of $\|s_0-s_{h\t}(0)\|^2_{\Kn(\Om)}$, after cancellation of terms and subsequent division, becomes unity.
 
\textbf{Step 2 (Estimate \eqref{eq:UpperBoundEstimate2nd}:} We choose the test function $\f=s-s_{h\t}\in \X$ in \eqref{eq:RobustnessEqAll}. Termwise, this gives 
\begin{align}\label{eq:RobustMainIneq1}
&\int_0^T \langle \p_t (s-s_{h\t}),s-s_{h\t} \rangle = \tfrac{1}{2}\|s(T)-s_{h\t}(T)\|^2- \tfrac{1}{2}\|s_0-s_{h\t}(0)\|^2,\\
 &\int_0^T \langle \mathcal{R}, s-s_{h\t}\rangle \overset{\eqref{eq:DualNormDef}}\leq  \int_0^T \|\mathcal{R}\|_{\Kn(\Om)}\|s-s_{h\t}\|_{\Kh(\Om)}\nonumber\\
&  \leq  \int_0^T \tfrac{2}{D_\mathrm{m}(t)} \|\mathcal{R}\|^2_{\Kn(\Om)} +  \int_0^T \tfrac{D_\mathrm{m}(t)}{8}  \|s-s_{h\t}\|_{\Kh(\Om)}^2\nonumber\\
& \overset{\eqref{eq:DmFpm}}\leq 2 \int_0^T \tfrac{1}{D_\mathrm{m}(t)} \|\mathcal{R}\|^2_{\Kn(\Om)} +  \tfrac{1}{8} \int_0^T\int_{\Om} D(s)|\K^{\frac{1}{2}}\del (s-s_{h\t})|^2,
\end{align}\begin{align}%%\\%%\\
& \int_0^T (f(s,\bm{x},t)-f(s_{h\t},\bm{x},t)),s-s_{h\t})\overset{\eqref{eq:DmFpm}}\leq \int_0^T f_{\p,\mathrm{M}}(t)  \|s-s_{h\t}\|^2,\\
& \int_0^T (\K \vg(\k(s_{h\t})-\k(s)),\del (s-s_{h\t})) \leq  \int_0^T  \left [\tfrac{2\vg^{\mathrm{T}}\K \vg}{D_{\mathrm{m}}(t)}\|\k(s)-\k(s_{h\t})\|^2 + \tfrac{D_{\mathrm{m}}(t)}{8} \|s-s_{h\t}\|_{\Kh(\Om)}^2 \right ]\nonumber\\
&\quad \overset{\ref{prop:Kabs}}\leq 2 K_\mathrm{M}|\bm{g}|^2 \|\k'\|_{L^\infty([0,1])}^2 \int_0^T \tfrac{1}{D_{\mathrm{m}}(t)}\|s-s_{h\t}\|^2 + \tfrac{1}{8}\int_0^T \int_{\Om} D(s) |  \K^{\frac{1}{2}} \del (s-s_{h\t})|^2.
\end{align}

 To estimate $\|\Psi-\Psi_{h\t}\|_{L^2(0,T;\Kh(\Om))}$, we need to also consider the parabolic-elliptic degeneracy. Consider the domains $\Om^1(t):=\{\bm{x}\in \Om :s(\bm{x},t),\, s_{h\t}(\bm{x},t)<1\}$, $\Om^2(t):=\{\bm{x}\in \Om :s(\bm{x},t)=1, s_{h\t}(\bm{x},t)<1\}$, $\Om^3(t):=\{\bm{x}\in \Om :s(\bm{x},t)<1, s_{h\t}(\bm{x},t)=1\}$, and $\Om^4(t):=\{\bm{x}\in \Om :s(\bm{x},t)=s_{h\t}(\bm{x},t)=1\}$ where the equalities and inequalities are satisfied in an almost everywhere sense inside the domains. We divide  accordingly the remaining term of \eqref{eq:RobustnessEqAll}
 $$
\int_0^T (\K\del(\Psi-\Psi_{h\t}),\del(s-s_{h\t}))=2 T_1 + T_2 + T_3 + T_4,
 $$ 
 where the terms $T_{1,2,3,4}$ are explained below.
 
$\bullet\quad$  
Observing that $\Sf(\Psi),\,\Sf(\Psi_{h\t})< 1$ a.e. in $\Om^1(t)$,   the first term $T_1$ is divided into two parts 
\begin{subequations}
\begin{align}
T_1&:=\tfrac{1}{2}\int_0^T(\K\del (\Psi-\Psi_{h\t}),  \del(s-s_{h\t}))_{\Om^{1}}=\tfrac{1}{2}\int_0^T(\K\del (\Psi-\Psi_{h\t}),  \del(\Sf(\Psi)-\Sf(\Psi_{h\t})))_{\Om^{1}}\nonumber\\
&=\tfrac{1}{2}\int_0^T(\K\del (\Psi-\Psi_{h\t}),  (\Sf'(\Psi)\del \Psi-\Sf'(\Psi_{h\t})\del \Psi_{h\t}))_{\Om^{1}}\nonumber\\
&=\tfrac{1}{2}\int_0^T(\K\del (\Psi-\Psi_{h\t}),  \Sf'(\Psi)\del (\Psi-\Psi_{h\t})+ (\Sf'(\Psi)-\Sf'(\Psi_{h\t}))\del \Psi_{h\t})_{\Om^{1}}\nonumber\\
&\overset{\eqref{eq:profSf}}= \tfrac{1}{2}\int_0^T\int_{\Om^{1}} \tfrac{|\K^{\frac{1}{2}}\del(\Psi-\Psi_{h\t})|^2}{D(\Sf(\Psi))} + \tfrac{1}{2}\int_0^T(\K\del (\Psi-\Psi_{h\t}), (\Sf'(\Psi)-\Sf'(\Psi_{h\t}))\del \Psi_{h\t})_{\Om^{1}}.
\end{align}
The second term on the right  is estimated as 
\begin{align}
&\tfrac{1}{2}\int_0^T (\K\del (\Psi-\Psi_{h\t}), (\Sf'(\Psi)-\Sf'(\Psi_{h\t}))\del \Psi_{h\t})_{\Om^{1}}\nonumber\\
&\overset{\eqref{eq:profSf}}= -\tfrac{1}{2}\int_0^T (\K\del (\Psi-\Psi_{h\t}), \left (\tfrac{D(s)-D(s_{h\t})}{D(s)D(s_{h\t})}\right )\del \Psi_{h\t})_{\Om^{1}}\nonumber\\
&= -\tfrac{1}{2}\int_0^T \left (\frac{1}{D(s)} \K\del (\Psi-\Psi_{h\t}), (D(s)-D(s_{h\t}))\del s_{h\t}\right )_{\Om^{1}}\nonumber\\
&\geq -\tfrac{1}{2}\int_0^T \|\K^{\frac{1}{2}} \del s_{h\t}\|_{L^\infty(\Om_1)} \int_{\Om^{1}} \frac{1}{D(\Sf(\Psi))} |D(s)-D(s_{h\t})||\K^{\frac{1}{2}} \del (\Psi-\Psi_{h\t}))|\nonumber\\
&\geq -\tfrac{1}{4}\int_0^T \|\K^{\frac{1}{2}} \del s_{h\t}\|^2_{L^\infty(\Om)} \int_{\Om^{1}} \frac{|D(s)-D(s_{h\t})|^2}{D(\Sf(\Psi))}  - \tfrac{1}{4}\int_0^T \int_{\Om^{1}} \tfrac{|\K^{\frac{1}{2}} \del (\Psi-\Psi_{h\t})|^2}{D(\Sf(\Psi))}\nonumber\\
&\geq  -\int_0^T   \tfrac{C^\infty_{h\t}(t) D^2_{M}(t)}{4 D_{\mathrm{m}}(t)}\, \int_{\Om^{1}} |s-s_{h\t}|^2 - \tfrac{1}{4}\int_0^T \int_{\Om^{1}}\tfrac{|\K^{\frac{1}{2}} \del (\Psi-\Psi_{h\t})|^2}{D(\Sf(\Psi))}.\label{eq:MiddleOfBigIneq}
\end{align}
\end{subequations}
Hence, we have
\begin{align}
T_1\geq \int_0^T \int_{\Om^{1}}\tfrac{|\K^{\frac{1}{2}} \del (\Psi-\Psi_{h\t})|^2}{4 D(\Sf(\Psi))}-\int_0^T   \tfrac{C^\infty_{h\t}(t) D^2_{M}(t)}{4 D_{\mathrm{m}}(t)}\, \|s-s_{h\t}\|^2.
\end{align}

$\bullet\quad$ We estimate $T_1$ once again. 
Recall that $s,\,s_{h\t}< 1$ a.e. in $\Om^1(t)$ implying $\Psi_{h\t}=P_\mathrm{c}(s_{h\t})$ and $\Psi=P_\mathrm{c}(s)$ in $\Om^{1}(t)$. Note from \eqref{eq:propPc} that ${P_\mathrm{c}}'=D$. Hence, we have
\begin{align}
T_1&=\tfrac{1}{2} \int_0^T(\K\del (\Psi-\Psi_{h\t}), \del (s-s_{h\t}))_{\Om^{1}}=\tfrac{1}{2}\int_0^T(\K \del (P_\mathrm{c}(s)-P_{\mathrm{c}}(s_{h\t})), \del (s-s_{h\t}))_{\Om^{1}}\nonumber\\
&=\tfrac{1}{2}\int_0^T (\K (D(s)\del s-D(s_{h\t})\del s_{h\t}), \del(s-s_{h\t}))_{\Om^{1}}\nonumber\\
&=\tfrac{1}{2}\int_0^T\int_{\Om^{1}} D(s)|\K^{\frac{1}{2}}\del (s- s_{h\t})|^2 + \tfrac{1}{2} \int_0^T (\K(D(s)-D(s_{h\t}))\del s_{h\t},\del (s-s_{h\t}))_{\Om^{1}}.\label{eq:2ndMainIneqLeadingTerm}
\end{align}
Similar to \eqref{eq:MiddleOfBigIneq}, the second term is estimated as
\begin{align}
& \tfrac{1}{2}\int_0^T ((D(s)-D(s_{h\t}))\K \del s_{h\t},\del (s-s_{h\t}))_{\Om^{1}}\nonumber\\
%&\quad \geq -\tfrac{1}{2}\int_0^T \|\K^{\frac{1}{2}} \del s_{h\t}\|_{L^\infty(\Om)} \int_\Om |D(s)-D(s_{h\t})|\,|\K^{\frac{1}{2}} \del (s-s_{h\t}))|\nonumber\\
%&\quad \geq -\tfrac{1}{2}\int_0^T \tfrac{\|\K^{\frac{1}{2}} \del s_{h\t}\|^2_{L^\infty(\Om)}}{D_{\mathrm{m}}(t)}\int_{\Om}(D(s)-D(s_{h\t}))^2 - \tfrac{1}{8}\int_0^T  D_{\mathrm{m}}(t) \int_{\Om}|\K^{\frac{1}{2}}\del (s-s_{h\t})|^2\nonumber\\
&\quad \geq -\int_0^T   \tfrac{C^\infty_{h\t}(t) D^2_{\mathrm{M}}(t)}{4 D_{\mathrm{m}}(t)} \|s-s_{h\t}\|^2 - \tfrac{1}{4}\int_0^T \int_{\Om^{1}}  D(s) |\K^{\frac{1}{2}} \del (s-s_{h\t})|^2.
\end{align}
Hence, we have so far that
\begin{equation}
2T_1 \geq \tfrac{1}{4}\int_0^T \int_{\Om^{1}}\left [ \tfrac{|\K^{\frac{1}{2}} \del (\Psi-\Psi_{h\t})|^2}{D(\Sf(\Psi))} + D(s) |\K^{\frac{1}{2}}\del (s- s_{h\t})|^2 \right ]-\int_0^T   \tfrac{C^\infty_{h\t}(t) D^2_{M}(t)}{2 D_{\mathrm{m}}(t)}\, \|s-s_{h\t}\|^2.
\end{equation}

$\bullet\quad$ Observe that $s=1$ in $\Om^2(t)$ and $\Sf'(\Psi_{h\t})=1\slash D(\Sf(\Psi_{h\t}))$ using \eqref{eq:profSf}. Also, $\int_{\Om^2} |\K^{\frac{1}{2}}\del\Psi|^2 \leq \int_\Om  |\K^{\frac{1}{2}}\del[\Psi-P_{\mathrm{M}}]_+|^2$. Moreover,  
 $\Om^2(t)\subseteq \Om^{\mathrm{deg}}(t)$ implying  that $D(s)/2 \leq D(s_{h\t})\leq 2 D(s)$ in $\Om^2(t)$ from the assumptions of \Cref{theo:UpperBound}. Using these, we have
\begin{align}
 T_2&:=\int_0^T(\K\del (\Psi-\Psi_{h\t}),  \del (s-s_{h\t}))_{\Om^2}=\int_0^T(\Sf'(\Psi_{h\t}) \K\del (\Psi_{h\t}-\Psi),  \del \Psi_{h\t})_{\Om^2} \nonumber\\
&= \tfrac{1}{2}\int_0^T\int_{\Om^2} \left [ \Sf'(\Psi_{h\t}) |\K^{\frac{1}{2}}\del \Psi_{h\t}|^2 + \Sf'(\Psi_{h\t}) |\K^{\frac{1}{2}}\del (\Psi_{h\t}-\Psi)|^2 - \Sf'(\Psi_{h\t}) |\K^{\frac{1}{2}}\del \Psi|^2\right ]\nonumber\\
&= \tfrac{1}{2}\int_0^T\int_{\Om^2} \left [ D(s_{h\t})|\K^{\frac{1}{2}}\del s_{h\t}|^2 + \tfrac{ |\K^{\frac{1}{2}}\del (\Psi-\Psi_{h\t})|^2}{D(s_{h\t})} \right ]- \tfrac{1}{2}\int_0^T\int_{\Om^2} \tfrac{ |\K^{\frac{1}{2}}\del\Psi|^2}{D(s_{h\t})}\nonumber\\
&\geq  \tfrac{1}{4}\int_0^T\int_{\Om^2} D(s)|\K^{\frac{1}{2}}\del (1- s_{h\t})|^2 + \int_0^T\int_{\Om^2} \tfrac{ |\K^{\frac{1}{2}}\del (\Psi-\Psi_{h\t})|^2}{4 D(s)}- \int_0^T\int_{\Om} \tfrac{ |\K^{\frac{1}{2}}\del[\Psi-P_{\mathrm{M}}]_+|^2}{D(1)}.
\end{align}
\vspace{1em}
In the above inequality, the identity $(a-b)a=\frac{1}{2}[a^2 + (a-b)^2 -b^2]$ has been used. 
 
$\bullet\quad$ With the same manipulations one has (note that $D(s)\leq 2 D(1)$ in $\Om^3(t)$ from the assumptions of \Cref{theo:UpperBound})
 \begin{align}
 T_3&:=\int_0^T(\K\del (\Psi-\Psi_{h\t}),  \del (s-s_{h\t}))_{\Om^3} \nonumber\\
&\geq  \tfrac{1}{2}\int_0^T\int_{\Om^3} D(s)|\K^{\frac{1}{2}}\del (1- s)|^2 + \int_0^T\int_{\Om^3} \tfrac{ |\K^{\frac{1}{2}}\del (\Psi-\Psi_{h\t})|^2}{2 D(s)}- \int_0^T\int_{\Om} \tfrac{ |\K^{\frac{1}{2}}\del[\Psi_{h\t}-P_{\mathrm{M}}]_+|^2}{D(1)}.
\end{align}
$\bullet\quad$ Finally, in $\Om^4(t)$ one has $s=s_{h\t}=1$, thus giving
\begin{align}
 T_4:=\int_0^T(\K\del (\Psi-\Psi_{h\t}),  \del (s-s_{h\t}))_{\Om^4}=0.
 \end{align} 
  With this, we have 
\begin{align}
&2T_1+ T_2+ T_3+T_4\geq \tfrac{1}{4}\int_0^T \int_{\Om}\left [ \tfrac{|\K^{\frac{1}{2}} \del (\Psi-\Psi_{h\t})|^2}{D(\Sf(\Psi))} + D(s)|\K^{\frac{1}{2}}\del (s- s_{h\t})|^2 \right ]\nonumber\\
 &\qquad -\int_0^T   \tfrac{C^\infty_{h\t}(t) D^2_{M}(t)}{2 D_{\mathrm{m}}(t)}\, \|s-s_{h\t}\|^2 - \tfrac{1}{D(1)} \int_0^T \left (\|[\Psi-P_{\mathrm{M}}]_+\|^2_{\Kh(\Om)} + \|[\Psi_{h\t}-P_{\mathrm{M}}]_+\|^2_{\Kh(\Om)} \right ).
 \end{align} 
 
$\bullet\quad$ To estimate $\|[\Psi-P_{\mathrm{M}}]_+\|_{\Kh(\Om)}$ insert $\f=[\Psi-P_{\mathrm{M}}]_+$ in \eqref{eq:Richards}. Note that $\p_t \Sf(\Psi)=0$ and $f(s,\bm{x},t)=f(1,\bm{x},t)$ if $\Psi> P_{\mathrm{M}}$. Also, $\int_0^T (\bm{c},\del[\Psi-P_{\mathrm{M}}]_+)=\int_0^T\int_{\p\Om} \bm{c}\cdot \hat{\bm{n}}_{\p\Om}[\Psi-P_{\mathrm{M}}]_+=0$ for the constant vector $\bm{c}=\smallint_{\Om^{\mathrm{deg}}(t)} \K\bm{g}$. Moreover, $f [\Psi-P_{\mathrm{M}}]_+\leq [f]_+[\Psi-P_{\mathrm{M}}]_+$. Using these relations leads to
\begin{align*}
&\int_0^T \|[\Psi-P_{\mathrm{M}}]_+\|_{\Kh(\Om)}^2 = \int_0^T (f(1,\bm{x},t),[\Psi-P_{\mathrm{M}}]_+) - \int_0^T (\K\bm{g}  ,\del [\Psi-P_{\mathrm{M}}]_+)\\
&\leq \int_0^T  ([f(1,\bm{x},t)]_+,[\Psi-P_{\mathrm{M}}]_+)-\int_0^T (\K\bm{g}-\tfrac{1}{|\Om^{\mathrm{deg}}|}\smallint_{\Om^{\mathrm{deg}}} \K\bm{g},\del [\Psi-P_{\mathrm{M}}]_+) \nonumber\\
&\leq \int_0^T \left ( \|[f(1,\bm{x},t)]_+\|_{\Kn(\Om^{\mathrm{deg}})}+ \|\K^{-\frac{1}{2}}(\K\bm{g}-\tfrac{1}{|\Om^{\mathrm{deg}}|}\smallint_{\Om^{\mathrm{deg}}} \K\bm{g})\|_{\Om^{\mathrm{deg}}}\right ) \|[\Psi-P_{\mathrm{M}}]_+\|_{\Kh(\Om)}.
\end{align*}
Using Young's inequality on the right hand side and recalling the definition of $\etaDeg$ we estimate
\begin{align}
\tfrac{1}{D(1)} \int_0^T \left (\|[\Psi-P_{\mathrm{M}}]_+\|^2_{\Kh(\Om)} + \|[\Psi_{h\t}-P_{\mathrm{M}}]_+\|^2_{\Kh(\Om)} \right )\leq \tfrac{1}{2} \int_0^T [\etaDeg]^2.
\end{align}

$\bullet\quad$ Combining all the estimates above, one obtains
\begin{align}
&\|s(T)-s_{h\t}(T)\|^2 + \tfrac{1}{2}\int_0^T \int_{\Om} \tfrac{|\K^{\frac{1}{2}}\del (\Psi-\Psi_{h\t})|^2}{D(\Sf(\Psi))} \nonumber\\
&\leq \|s_0-s_{h\t}(0)\|^2 + 4 \int_0^T \tfrac{1}{D_\mathrm{m}(t)} \|\mathcal{R}\|^2_{\Kn(\Om)} + \int_0^T [\etaDeg]^2+  \int_0^T \mathfrak{C}_2(t) \|s-s_{h\t}\|^2.
\end{align}
Since $\mathfrak{C}_2(t)>0$, one has \eqref{eq:UpperBoundEstimate2nd} from applying the Gronwall Lemma  \eqref{eq:GronwallLemma}, where $u(t)= \|s(t)-s_{h\t}(t)\|^2$, $\b(t)= \mathfrak{C}_2(t)$ and 
\begin{align*}
\a(t)=\|s_0-s_{h\t}(0)\|^2 +  4 \int_0^t\tfrac{1}{D_\mathrm{m}(t)}  \|\mathcal{R}\|^2_{\Kn(\Om)} + \int_0^t [\etaDeg]^2 -  \tfrac{1}{2}\int_0^t  \int_{\Om}\tfrac{|\K^{\frac{1}{2}} \del (\Psi-\Psi_{h\t})|^2}{D(\Sf(\Psi))}.
\end{align*}

\textbf{Step 3 (Estimate \eqref{eq:UpperBoundEstimate3rd}):}
Using the definition of $H^{\pm 1}_\K$-norms in \eqref{eq:RobustnessEqAll}, we have
\begin{align*}
\int_0^t \|\p_t (s-s_{h\t})\|^2_{\Kn(\Om)}\leq 3\int_0^t \left [ \|\Psi-\Psi_{h\t}\|^2_{\Kh(\Om)} +  \|\res\|^2_{\Kn(\Om)} + \mathfrak{C}_3(T)\|s-s_{h\t}\|^2 \right ],
\end{align*}
for any $t\in (0,T]$. Multiplying the above inequality with $\lambda(t)\exp(\int_t^T \lambda)$, integrating on $[0,T]$, and adding the above inequality for $t=T$, we get from the first term 
\begin{align*}
&\int_0^T\left [\|\p_t (s-s_{h\t})\|^2_{\Kn(\Om)} +\lambda(t)\exp\left (\smallint_t^T \lambda\right ) \int_0^t \|\p_t (s-s_{h\t})\|^2_{\Kn(\Om)}\right ]\\
\overset{\eqref{eq:TiemIntegrator}}=  &\exp\left (\smallint_0^T \lambda \right )\J_\lambda\left (\|\p_t (s-s_{h\t})\|_{\Kn(\Om)}\right )^2,
\end{align*}
and similar for the other terms. The estimate \eqref{eq:UpperBoundEstimate3rd} follows then by cancelling the $\exp\left (\int_0^T \lambda \right )$ multipliers.
\end{proof}

\begin{remark}[Upper bound on $\|D(s)^{\frac{1}{2}}\K^{\frac{1}{2}}\del(s-s_{h\t})\|$] From the step 2 of the proof of \Cref{theo:UpperBound}, it is 
evident that the error component $\|D(s)^{\frac{1}{2}}\K^{\frac{1}{2}}\del(s-s_{h\t})\|$ can be estimated as well through slight changes in coefficients of the right hand side of \eqref{eq:UpperBoundEstimate2nd}. However, to have symmetry between the lower and the upper bounds of \Cref{theo:TotalPressureEfficiency,theo:UpperBound}, this has not been pursued.
\end{remark}

\section{Finite element discretization}
\label{sec:FiniteElement}
We describe in this section the discretization of the Richards problem \eqref{eq:FullSystem} by the finite element method.

\subsection{Time steps}
For the time-interval $(0,T)$, we introduce $N+1$ discrete times $\bm{t}_N=(t_n)_{n=0}^N$ where $t_0=0<t_1<\cdots< t_n<\cdots<t_N=T$. Let $I_n:=(t_{n-1},t_n] $ denote the time intervals and $\t_n:= t_{n}-t_{n-1}$ the lengths of the time steps for $n\in \{1,\dots, N\}$. Note that, we allow nonuniform time stepping. Further, for a vector space $V$, $\mathcal{Q}_{1}\left(I_n;V\right)$ denotes  the space of $V$-valued affine functions over the time-step interval $I_n$.

\subsection{Space meshes}\label{sec:fem_mesh}
For the time sequence $\bm{t}_N$, let $\{\calT_n\}_{n=1}^N$ denote the sequence of matching and uniformly shape regular simplical meshes for the domain $\Om$. The meshes are allowed to undergo refinement or coarsening between time steps. Henceforth, discontinuities of $\K$ are only allowed to happen along internal edges of the mesh.
For each element~$K\in\calT_n$, let $h_K:= \mathrm{diam} \{K\}$ denote the diameter of $K$ and let $\mathfrak{p}_n\geq 1$ denote the spatial polynomial degree associated with $\calT_n$. Our results are generalizable to  polynomial degrees depending locally on $K\in\calT_n$. However, to keep the notation simple, we only consider $\mathfrak{p}_n$ changing between time steps here. For full $hp$-adaptive  algorithm, we refer to \cite{Ern_Sme_Voh_heat_HO_Y_17}.

\subsection{Approximation spaces}\label{sec:fem_approximation}
On a time step $n\in \{1,\dots,N\}$, we define the $H^1_0(\Om)$-conforming $hp$-finite element space $V_{n,h}$ as:
\begin{equation}\label{eq:conforming_space_def}
V_{n,h}:=\left\{u_h \in H^1_0(\Om),\; \eval{u_h}{K} \in \calP_{\mathfrak{p}_n}(K)\quad\forall\,K\in\calT_n\right\},
\end{equation}
where $\calP_{\mathfrak{p}_n}(K)$ denotes the polynomial space of degree $\mathfrak{p}_n\in \N$ on $K$.
Further, let $\Pi_{n,h}:L^2(\Om)\to V_{n,h}$ and $\Lam_{n,h}:L^2(\Om)\to \calP_{\mathfrak{p}_n}(\calT_n)$ represent the $L^2$-orthogonal projection operator with respect to the spaces $V_{n,h}$ and $\calP_{\mathfrak{p}_n}(\calT_n)$, i.e.,
\begin{subequations}\label{eq:DefPiL}
\begin{align}
&\Pi_{n,h}u\in V_{n,h} \text{ for } u\in L^2(\Om) \text{ is such that } (\Pi_{n,h} u, \f_h)=(u,\f_h), \text{ for all } \f_h\in V_{n,h};\\
&\Lam_{n,h}u\in \calP_{\mathfrak{p}_n}(\calT_n) \text{ for } u\in L^2(\Om) \text{ is such that } (\Lam_{n,h} u, \f_h)=(u,\f_h), \text{ for all } \f_h\in \calP_{\mathfrak{p}_n}(\calT_n).
\end{align}
\end{subequations}

\subsection{Finite element discretization}\label{sec:num_scheme}
In \Cref{sec:ErrorResidual}, formulation \eqref{eq:Richards} is used to derive the estimates. However, since \eqref{eq:RichardsP} is the  most general and commonly used formulation, we propose the finite element scheme for \eqref{eq:RichardsP}. We will still be able to apply the analysis of \Cref{sec:ErrorResidual}.
For time discretization, we consider the backward Euler scheme. The problem for each $n\in\{1,\dots,N\}$ and a given $S_{n-1,h}\in L^2(\Om)$ is to find $p_{n,h}\in V_{n,h}$ which satisfies for all $\f_{h} \in\Vn$,
\begin{align}\label{eq:TimeDiscrete}
& \tfrac{1}{\t_n}(S(p_{n,h})-S_{n-1,h}, \f_{h}) +  (\K\k(S(p_{n,h}))[\nabla p_{n,h} +\vg], \nabla \f_{h}) = (f(S(p_{n,h})),\bm{x},t_n),\f_{h}).
\end{align}
For $n=1$, we set $S_{n-1,h}:=\Pi_{1,h} s_0$, whereas, for $n>1$, $S_{n-1,h}:=S(p_{n-1,h})$. The existence of $p_{n,h}$ solving \eqref{eq:TimeDiscrete} is discussed in \cite{douglas1970galerkin} for the nondegenerate case ($\k(0)>0$). The degenerate case is covered  in \cite{oulhaj2018numerical} for the control volume finite element method. In practice, since the problem \eqref{eq:TimeDiscrete} is nonlinear, the  exact $p_{n,h}$ is generally not known, and linear iterations have to be used to approximate $p_{n,h}$. This is discussed at length in Appendix \ref{App:linear}.

From the sequence $\{p_{n,h}\}_{n=1}^N$, we define the space--time discrete total pressure and saturation for all $n\in \{1,\dots,N\}$ as
\begin{align}\label{eq:PsiDiscrete}
\Psi_{n,h}:=\calK(p_{n,h})\in H^1_0(\Om)  \text{ and } S_{n,h}:=\Sf(\Psi_{n,h})\overset{\eqref{eq:ConnectionVariables}}=S(p_{n,h})\in H^1(\Om).
\end{align}
The choice $\Psi_{0,h}=P_{\mathrm{c}}(S_{0,h})=P_{\mathrm{c}}(\Pi_{1,h} s_0)$ is used for extending the definition of $\Psi_{n,h}$ to $n=0$.

\subsection{Time-continuous solutions}\label{sec:time_continuous}
There are multiple ways to define a time-continuous total pressure $\Psi_{h\t}\in \X$ and saturation $s_{h\t}\in \W$, satisfying the requirements of \Cref{theo:UpperBound,theo:TotalPressureEfficiency}, starting from $\{S_{n,h}\}_{n=1}^N$ and $\{\Psi_{n,h}\}_{n=1}^N$ introduced in \eqref{eq:PsiDiscrete}. Here, for $t\in I_n$, we choose
\begin{subequations}\label{eq:timeContSol}
\begin{align}\label{eq:timeContSolPsi}
\Psi_{h\t}(t):=&P_{\mathrm{c}}\left (\tfrac{t-t_{n-1}}{\t_n} S_{n,h} + \tfrac{t_{n}-t}{\t_n} S_{n-1,h} \right )+\left [\tfrac{t-t_{n-1}}{\t_n}  \Psi_{n,h} +\tfrac{t_{n}-t}{\t_n} \Psi_{n-1,h}-P_{\mathrm{M}} \right ]_+,\\
s_{h\t}(t):=& \Sf(\Psi_{h\t}(t)). 
\end{align}
\end{subequations}
Observe that, $\Psi_{h\t}$ and $s_{h\t}$ defined this way satisfy
\begin{subequations}\label{eq:propPsiShtau}
\begin{align}
&\Psi_{h\t}\in C(0,T;H^1_0(\Om))\subset \X, &\text{ and } s_{h\t}\in W^{1,\infty}(0,T;H^1(\Om))\subset \W;\\
&\Psi_{h\t}(t_n)=\Psi_{n,h}, &\text{ and } s_{h\t}(t_n)=S_{n,h},\quad \forall \, n\in \{1,\dots,N\}.\label{eq:propPsiShtauB}
\end{align}
\end{subequations}
The relation \eqref{eq:propPsiShtauB} even holds when $\Psi_{n,h}>P_{\mathrm{M}}$ since $P_{\mathrm{c}}(S_{n,h})=P_{\mathrm{M}}$ in this case and the other contribution from \eqref{eq:timeContSolPsi} adds $[\Psi_{n,h} -P_{\mathrm{M}}]_+$.
Another advantage of this interpolation is that, using \eqref{eq:DefSfunc}, if both $\Psi_{n,h},\, \Psi_{n-1,h}\leq P_{\mathrm{M}}$ (nondegenerate case), or $\Psi_{n,h},\, \Psi_{n-1,h}\geq P_{\mathrm{M}}$ (degenerate case), i.e.,
\begin{align}\label{eq:NonDegDS}
\text{ if either } \Psi_{h\t}\leq P_{\mathrm{M}} \text{ or } \Psi_{h\t}\geq P_{\mathrm{M}}  \text{ in } I_n, \text{ then } \p_t s_{h\t}= \tfrac{1}{\t_n}(S_{n,h}-S_{n-1,h}).
\end{align}

\section{A posteriori error estimates}
\label{sec:APosteriori}
We apply here the developments of \Cref{sec:ErrorResidual} to perform a posteriori error analysis of the finite element discretization of \Cref{sec:FiniteElement}.

\subsection{Equilibrated flux}\label{sec:Equi_Fl}
The objective of this section is to design an equilibrated flux $\sth\in \Hdiv$ that satisfies the mass balance property
\begin{align}\label{eq:sigma_th_equilibration}
\int_K \left [\tfrac{1}{\t_n}(S_{n,h}-S_{n-1,h}) + \nabla{\cdot} \sth  - f(S_{n,h},\bm{x},t_n)\right ] = 0  \text{ for all } K\in \calT_n.
\end{align}

\subsubsection{Local mixed finite element spaces}\label{sec:MFEMspaces}
For the construction of $\sth$, we introduce some standard mixed finite element spaces.
For each $n\in \{1,\dots,N\}$, let $\calVh$ denote the set of vertices of the mesh $\calT_n$, where we distinguish the set of interior vertices $\calVhint$ and the set of  boundary vertices $\calVhext$. For $K\in\calT_n$, $\calV_K\subset \calVh$ denotes the set of vertices of $K$. For each $\ta \in \calVh$, let $\psia$ denote the hat function associated with $\ta$ and $\oma$ the interior of the support of $\psia$, with the associated diameter $h_{\oma}$. 
Furthermore, let $\Ta$ denote the restriction of the mesh $\CR$ to~$\oma$.

For a polynomial degree $\mathfrak{p}\geq 0$, the local spaces $\calP_{\mathfrak{p}}(\Ta)$ and $\RTN_{\mathfrak{p}}(\Ta)$ are defined by
\begin{align*}
\calP_{\mathfrak{p}}(\Ta) &:= \{ u_h \in L^2(\oma),\quad u_h|_{\elCR} \in \calP_{\mathfrak{p}}(\elCR)\quad\forall\,\elCR\in\Ta\},\\
 \RTN_{\mathfrak{p}}(\Ta) &:= \{ \bm{v}_h \in \bm{L}^2(\oma;\R^d),\quad \bm{v}_h|_{K} \in \RTN_{\mathfrak{p}}(K)\quad\forall\elCR\in\Ta\},
\end{align*}
where  $\RTN_{\mathfrak{p}}(\elCR) :=  \calP_{\mathfrak{p}}(\elCR;\R^d) + \calP_{\mathfrak{p}}(\elCR)\bm{x}$ denotes the Raviart--Thomas--N\'ed\'elec space of order $\mathfrak{p}$ on $\elCR$. We use a similar notation on the whole mesh $\calT_n$, and introduce the local mixed finite element spaces $\Va$ and $\Qa$ as
\begin{align}
\Va & :=
\begin{cases}
    &\left\{\bm{v}_h \in \RTN_{\mathfrak{p}_n + 1}(\Ta),\; \bm{v}_h\in \Hdivoma ,\; \bm{v}_h\cdot \bm{n} =0\text{ on }\p\oma \right\}\; \text{if }\ta\in\calVhint,\\
    &\left\{\bm{v}_h \in \RTN_{\mathfrak{p}_n + 1}(\Ta),\; \bm{v}_h\in \Hdivoma ,\; \bm{v}_h\cdot \bm{n} =0\text{ on }\p\oma\setminus\p\Om \right\}\;\text{if }\ta\in\calVhext,
\end{cases}
 \nonumber\\
\Qa & :=
\begin{cases}
     \left\{ u_h\in \calP_{\mathfrak{p}_n +1 }(\Ta) ,\quad (u_h,1)_\oma = 0\right\} & \hspace{4.13cm}\text{if }\ta\in\calVhint,\\
     \; \calP_{\mathfrak{p}_n+1}(\Ta)  & \hspace{4.13cm}\text{if }\ta\in\calVhext.
\end{cases}\label{eq:spacetime_mixed_space_def}
\end{align}
The projector $\PiH: \bm{L}^2(\Om;\R^d)\to \RTN_{\mathfrak{p}_n}(\calT_n)$ is then defined as:
\begin{align}\label{eq:DefProjection} \text{ for }  \bm{u}\in \bm{L}^2(\Om;\R^d),\quad (\K \PiH \bm{u},\bm{v}_h)=(\K\,\bm{u},\bm{v}_h),\;  \text{ for all } \bm{v}_h\in \RTN_{\mathfrak{p}_n}(\calT_n).
\end{align}
Note that it is computed elementwise.

\subsubsection{Flux reconstruction}\label{sec:flux_reconstruction_def}
For each $n\in \{1,\dots,N\}$, we unify the numerical source-like terms and flux-like terms of \eqref{eq:TimeDiscrete} in 
$\Gf_{n,h}\in L^2(\Om)$ and $\bm{F}_{n,h}\in \bm{L}^2(\Om;\R^d)$,
\begin{align}\label{eq:SourceFluxTerms}
\Gf_{n,h} &:=f(S_{n,h},\bm{x},t_n) - \tfrac{1}{\t_n}(S_{n,h}-S_{n-1,h}), \quad \bm{F}_{n,h}:= \nabla \Psi_{n,h} +\vg \k(S_{n,h}).
\end{align}
Observe that the terms defined above are constant in time in $I_n$.
Recalling the projection operators $\PiL$, $\Lam_{n,h}$, and $\PiH$  from \eqref{eq:DefPiL} and \eqref{eq:DefProjection},  the scalar function $\gtautha \in \calP_{\mathfrak{p}_n +1 }(\Ta)$ and the vector field $\tautha \in \RTN_{\mathfrak{p}_n + 1}(\Ta)$ are defined as
\begin{align}\label{eq:tau_g_def}
\gtautha:= (\psia\,\Lam_{n,h}\Snt -  \nabla  \psia\cdot \,\K\, \PiH\,\Fnt)|_{\oma}, \quad \tautha:= - (\psia \K\,\PiH\, \Fnt)|_{\oma}.
\end{align}
Since $\psia \in \Vn$, using $\f_{h}=\psia$ in
\eqref{eq:TimeDiscrete} we get directly for all $\ta\in \calVhint$ that $(\gtautha,1)_{\oma} = 0.$

\begin{definition}[Equilibrated flux $\sth$]\label{def:flux_construction_1}
For a given time-step  $n\in \{1,\dots,N\}$ and for each vertex $\ta\in\calVh$, let the mixed finite element spaces $\Va$ and $\Qa$ be defined by \eqref{eq:spacetime_mixed_space_def}.
For the time discrete solutions introduced in  \Cref{sec:num_scheme}, let $\Snt$ and $\Fnt$ be defined  in \eqref{eq:SourceFluxTerms}. 
Let $\gtautha$ and $\tautha$ be defined by \eqref{eq:tau_g_def}.
Furthermore, let $\stha \in \Va$ be defined by
\begin{equation}\label{eq:stha_minimization_def}
\stha := \argmin_{\substack{ \bm{v}_h \in \Va, \\ \nabla{\cdot} \bm{v}_h = \gtautha}} \|\K^{-\frac{1}{2}}(\bm{v}_h - \tautha)\|_{\oma}.
\end{equation}
Then, after extending $\stha$ by zero from $\oma$ to $\Om$ for each $\ta \in \calVh$, we define the equilibriated flux as
\begin{equation}\label{eq:flux_reconstruction_1}
\sth:= \sum_{\ta \in \calVh} \stha.
\end{equation}
\end{definition}
The well-posedness of $\sth$ follows from Theorem 4.2 of \cite{Ern_Sme_Voh_heat_HO_Y_17}, see also references therein, and it satisfies  \eqref{eq:sigma_th_equilibration} since 
\begin{subequations}\label{eq:GnhSigh}
\begin{align}\label{eq:GnhSigh1}
&\del\cdot\sth\overset{\eqref{eq:flux_reconstruction_1}}= \sum_{\ta \in \calVh} \del\cdot\stha\overset{\eqref{eq:tau_g_def}}= \sum_{\ta \in \calVh} [(\psia\,\Lam_{n,h}\Snt -  \nabla  \psia\cdot \,\K\, \PiH\,\Fnt)]= \Lam_{n,h}\Snt,\\
&\text{ and consequently } (\del\cdot\sth-\Snt,\f_h)_K\overset{\eqref{eq:DefPiL}}=0,\quad \forall K\in \calT_n \text{ and } \, \f_h\in \calP_{\mathfrak{p}_n}(\calT_n).\label{eq:GnhSigh2}
\end{align}
\end{subequations}
Here, the partition of unity property, $\sum_{\ta\in\calV_K}\psia=1$ is used. Practically, $\stha$ are computed by solving the following mixed finite element problems  \cite{Ern_Sme_Voh_heat_HO_Y_17} locally in $\oma$: find $\stha\in \Va$ and $r^{\ta}_{n,h}\in \Qa$ such that
\begin{align*}
& (\K^{-1}\stha, \bm{v}_h)_{\oma}- (\del \cdot \bm{v}_h, r^{\ta}_{n,h})=(\K^{-1} \tautha,\bm{v}_h)_{\oma},\quad \forall  \bm{v}_h\in \Va,\\
& (\del \cdot \stha, u_h)_{\oma}= (\gtautha,u_h)_{\oma} ,\quad \forall  u_h\in \Qa.
 \end{align*}

\subsection{A posterori error estimators}\label{sec:main_aposteriori}
Recalling the definition of time-continuous solutions $(\Psi_{h\t},s_{h\t})$ from \Cref{sec:time_continuous}, we introduce the following a posteriori error estimators: Take $n\in \{1,\dots,N\}$, an open polytope $\om\subseteq \Om$, and $t\in I_n$. Then,
\begin{subequations}\label{eq:estimators}
\begin{align}
\eta^{\mathrm{F}}_{n,h,\om}(t) &:=  \|\K^{-\frac{1}{2}}\sth + \K^{\frac{1}{2}}(\del\Psi_{h\t} + \vg\,\k(s_{h\t}))(t)\|_{\om} \label{eq:etaEq_def}
\end{align}
measures the lack of $\Hdiv$-conformity of the numerical flux $\K(\del\Psi_{h\t} + \vg\,\k(s_{h\t}))$. The quadrature error estimator arising from $\Gf_{n,h}$  not being polynomial (see  \eqref{eq:SourceFluxTerms}) is 
\begin{align}
\eta^{\mathrm{qd},\Gf}_{n,h,\om} & := \frac{h_{\om} }{\sqrt{K_{\mathrm{m}}}\pi} \| \Gf_{n,h} - \Lam_{n,h} \Snt\|_{\om}.\label{eq:etaOscS_S} 
\end{align}
The time-quadrature error of $\p_t s_{h\t}$ is measured by  the estimator
\begin{align}
\eta^{\mathrm{qd},t}_{n,h,\om}(t):=\|\p_t s_{h\t}- \tfrac{1}{\t_n}(S_{n,h}-S_{n-1,h}) \|_{\Kn(\om)}.\label{eq:etaOscS_t} 
\end{align}
Observe that it estimates quadrature since $\int_{I_n} \p_t s_{h\t}\overset{\eqref{eq:propPsiShtauB}}= S_{n,h}-S_{n-1,h}$, and it vanishes in both purely degenerate and nondegenerate regimes  due to \eqref{eq:NonDegDS}.
 The temporal oscillation in data $f$ is measured  by
\begin{align}
\etaOsc(t) &:=  \|f(s_{h\t}(t_n),\bm{x},t_n)-f(s_{h\t}(t),\bm{x},t)\|_{\Kn(\om)}.
\end{align}
The errors in the approximation of the initial condition $s_0$ are accounted by
\begin{align}
\etaOscInitL & :=  \|s_0-\Pi_{1,h} s_0\|,\quad \etaOscInitH:=  \|s_0-\Pi_{1,h} s_0\|_{\Kn(\Om)}.\label{eq:etaIni_def}
\end{align}
The projectors $\Pi_{n,h}$ and $\Lam_{n,h}$ were defined in \eqref{eq:DefPiL},  and the norm $\|\cdot\|_{\Kn(\cdot)}$ was introduced in \eqref{eq:NormHminK}.  With the above definitions, the total estimator is computed as
\begin{align}
\eta_{\mathcal{R}}(t):= \left [ \sum_{K\in \calT_n} [\etaEq(t)  + \etaOscS]^2 \right ]^{\frac{1}{2}}+\etaOscTom(t)+ \eta^{\mathrm{osc}}_{n,\Om}(t).
\end{align}
\end{subequations}

\begin{remark}[Inverse of the Kirchhoff transform] The inverse of the Kirchhoff transform $\calK^{-1}$ (see \eqref{eq:TotalPandP}) does not need to be evaluated for computing the estimators. 
\end{remark}

\subsection{Global reliability}
Complementing \Cref{theo:UpperBound}, our a posteriori error estimate on the error in the finite element discretization \eqref{eq:TimeDiscrete} of the Richards equation \eqref{eq:FullSystem} is
\begin{theorem}[Global reliability]
Recall the definitions and assumptions stated in \Cref{theo:UpperBound}. Let $\{\Psi_{n,h}\}_{n=1}^{N}\subset H^1_0(\Om)$ and $\{S_{n,h}\}_{n=1}^{N}\subset H^1(\Om)$ be defined using the finite element discretization \eqref{eq:TimeDiscrete}--\eqref{eq:PsiDiscrete} and let $\Psi_{h\t}\in C(0,T;H^1_0(\Om))\subset \X$ with $s_{h\t}=\Sf(\Psi_{h\t})\in W^{1,\infty}(0,T;H^1(\Om))\subset \W$  be their time-continuous interpolates as defined  in \eqref{eq:timeContSol}. Let the a posteriori error estimators be defined in~\eqref{eq:estimators}.  Then, for any time $t\in [0,T]$,
\begin{subequations}\label{eq:GlobalReliabilityEstimates}
\begin{align}\label{eq:ResidualBoundedByEtaR}
 \|\res(\Psi_{h\t}(t))\|_{\Kn(\Om)}\leq \eta_{\mathcal{R}}(t).
\end{align}
Consequently, the errors of $s_{h\t}$ and $\Psi_{h\t}$ satisfy: 
\begin{align}\label{eq:EL2etaL2}
\calE_{L^2}^2 &:
={ e^{-\smallint_0^T (\lambda + \mathfrak{C}_1)} }\|(s-s_{h\t})(T)\|^2_{\Kn(\Om)}
+\J_{\lambda + \mathfrak{C}_1}(\Sf_{\p,\mathrm{M}}^{-\frac{1}{2}}\, \|s-s_{h\t}\|)^2 \nonumber\\
&\leq [\etaOscInitH]^2+ \J_{\lambda + \mathfrak{C}_1}(\lambda^{-\frac{1}{2}} \eta_{\mathcal{R}})^2 =:\eta_{L^2}^2,\\
\calE_{H^1}^2&:=
e^{-\smallint_0^T \mathfrak{C}_2 }\|(s-s_{h\t})(T)\|^2 +
 \tfrac{1}{2}\J_{\mathfrak{C}_2}(\|D(s)^{-\frac{1}{2}}\K^{\frac{1}{2}}\del(\Psi-\Psi_{h\t})\|)^2\nonumber\\
&\leq [\etaOscInitL]^2+ \J_{\mathfrak{C}_2}\left ( \eta^{\mathrm{deg}} \right )^2 + 4  \, \J_{\mathfrak{C}_2}\left ( D_{\mathrm{m}}^{-\frac{1}{2}} \eta_{\res} \right )^2 =:\eta_{H^1}^2.\label{eq:EH1etaH1}
\end{align}
\end{subequations}
\label{theo:GlobalReliability}
\end{theorem}\vspace{-2em}

\begin{proof}
From the regularity of $\Psi_{h\t}$ and $\p_t s_{h\t}\in L^\infty(0,T;L^2(\Om))$ one has that $\res(\Psi_{h\t})\in L^\infty(0,T;H^{-1}(\Om))$. Hence, for all $t\in I_n$, $n\in \{1,\dots,N\}$, adding and subtracting $(\sth,\del \f)$,
\begin{align}
&\langle \res(\Psi_{h\t}),\, \f\rangle=(f(s_{h\t},\bm{x},t)-\p_t s_{h\t}-\del\cdot\sth,\f) - (\sth+\K[\del \Psi_{h\t} + \vg\k(s_{h\t})],\del \f) \nonumber\\
&= (\Snt-\del\cdot\sth,\f)  + (f(s_{h\t},\bm{x},t)-f(S_{n,h},\bm{x},t_n),\f) \nonumber\\
&\quad + (\tfrac{1}{\t_n}(S_{n,h}-S_{n-1,h})-\p_t s_{h\t},\f) - \sum_{K\in\calT_n}(\sth+\K[\del \Psi_{h\t} + \vg\k(s_{h\t})],\del \f)_K,\label{eq:UpperBoundOfResidual}
\end{align}
where $\Gf_{n,h}$ is defined in \eqref{eq:SourceFluxTerms}. For the first term on the right, following \eqref{eq:GnhSigh}, we use 
\begin{subequations}\label{eq:EstimatesForReliability}
\begin{align}
&(\Gf_{n,h}-\del\cdot\sth,\f)=\sum_{\elCR \in \CR} (\Gf_{n,h}-\del\cdot\sth,\f)_{\elCR} \overset{\eqref{eq:GnhSigh2}}= \sum_{\elCR\in \CR}(\Gf_{n,h} -\del\cdot\sth,\, \f-\tfrac{1}{|\elCR|}\smallint_{\elCR}\f)_{\elCR}
\nonumber\\
&\overset{\eqref{eq:GnhSigh1}}= \sum_{K\in \calT_n} (\Gf_{n,h}-\Lam_{n,h}\Snt,\,\f-\tfrac{1}{|\elCR|}\smallint_{\elCR}\f)_{\elCR} \overset{\eqref{eq:Poincare},\, \eqref{eq:CurrencyKnorms}}\leq \sum_{K\in \calT_n} [\etaOscS ]\, \|\K^{\frac{1}{2}}\del \f\|_K.
\end{align}
For the rest of the terms in \eqref{eq:UpperBoundOfResidual}, we note 
\begin{align}
&(f(s_{h\t},\bm{x},t)-f(S_{n,h},\bm{x},t_n),\f) \leq   [\eta^{\mathrm{osc}}_{n,\Om}(t)]\, \|\f\|_{\Kh(\Om)},\\
&(\tfrac{1}{\t_n}(S_{n,h}-S_{n-1,h})-\p_t s_{h\t},\f)\leq [\etaOscTom(t)] \, \|\f\|_{\Kh(\Om)},\\
& \sum_{K\in\calT_n}(\sth+\K[\del \Psi_{h\t} + \vg\k(s_{h\t})],\del \f)_K \leq \sum_{K\in\calT_n}[\etaEq(t)]\,\|\K^{\frac{1}{2}}\del \f\|_K.
\end{align}
\end{subequations}
Combining the inequalities of \eqref{eq:EstimatesForReliability} in \eqref{eq:UpperBoundOfResidual} and using the Cauchy--Schwarz inequality, one arrives at \eqref{eq:ResidualBoundedByEtaR}.
Estimates \eqref{eq:EL2etaL2}--\eqref{eq:EH1etaH1} then follow from inserting \eqref{eq:ResidualBoundedByEtaR} in \Cref{theo:UpperBound}. 
\end{proof}

\subsection{Quadrature and temporal discretization estimators}
For providing the efficiency bound, a few more estimators need to be introduced. For $n\in \{1,\dots,N\}$, an open polytope $\om\subseteq \Om$, and  $t\in I_n$, the quadrature estimator for the numerical flux is  defined as
\begin{subequations}\label{eq:estimators2}
\begin{align}
\eta^{\mathrm{qd},\bm{F}}_{n,h,\om} & :=\|\K^{\frac{1}{2}}(\bm{F}_{n,h}-\PiH\Fnt)\|_{\om}. \label{eq:etaOscS_F} 
\end{align}
To measure the temporal discretization error of the numerical solutions $\Psi_{n,h}$ and $S_{n,h}$, we  further introduce for $t\in I_n$ the estimators:
\begin{align}
\eta^{\mathrm{J},H^1}_{n,h,\om}(t) := \|\Psi_{h\t}(t)-\Psi_{n,h}\|_{\Kh(\om)}, \quad \eta^{\mathrm{J},L^2}_{n,h,\om}(t) := \|s_{h\t}(t)-S_{n,h}\|_\om.   \label{eq:etaJ_def}
\end{align}
\end{subequations}

\subsection{Local-in-space and in-time efficiency}
\begin{theorem}[Local and global efficiency]\label{theo:Y_norm_guaranteed_efficiency}
Let $\Psi\in \X$ with $s=\Sf(\Psi)\in \W$ be the weak solution of \eqref{eq:Richards}.   Let $\{\Psi_{n,h}\}_{n=1}^{N}\subset H^1_0(\Om)$ and $\{S_{n,h}\}_{n=1}^{N}\subset H^1(\Om)$ be defined using the finite element discretization \eqref{eq:TimeDiscrete}--\eqref{eq:PsiDiscrete} and let $\Psi_{h\t}\in C(0,T;H^1_0(\Om))\subset \X$ with $s_{h\t}=\Sf(\Psi_{h\t})\in W^{1,\infty}(0,T;H^1(\Om))\subset \W$,  be their time-continuous interpolates as defined  in \eqref{eq:timeContSol}. 
Let $\sth$ denote the equilibrated flux of \Cref{def:flux_construction_1}. Let the a posteriori error estimators be defined in~\eqref{eq:estimators} and \eqref{eq:estimators2}.  Let $\dst^{\a}_{\om,I}$ be defined   in~\eqref{eq:DefLBErrorMeasure} for $\a(t)=\max\limits_{\ta\in \calV_n}\{C_{\mathrm{P},\oma}\,h_{\oma}\}\,K_{\mathrm{m}}^{-\frac{1}{2}}\,\max\limits_{[0,1]\times\Om\times \{t\}}|\p_s f|  + |\vg| K_{\mathrm{M}}^{\frac{1}{2}} \|\k'\|_{L^\infty([0,1])}$ and $t\in I_n$. Then, for each discrete time step $n\in\{1,\dots,N\}$ and mesh element $K\in \calT_n$, the indicators satisfy the following local-in-space and in-time efficiency bound: 
\begin{align}\label{eq:LocalEfficiencyMainEstimate}
&\int_{I_n}( [\etaEq]^2 + [\etaJh]^2 )\nonumber\\
&
\lesssim\; \sum_{\ta\in \calV_K}\left ( \int_{I_n} \left [\sum_{j\in \{\Gf,\bm{F},t\}}[\eta^{\mathrm{qd},j}_{n,h,\oma}]^2 +  [\eta^{\mathrm{osc}}_{n,\oma}]^2  + \a^2\,[\eta^{\mathrm{J},L^2}_{n,h,\oma}]^2 + [\eta^{\mathrm{J},H^1}_{n,h,\oma}]^2\right ] +  \dst^{\a}_{\oma,I_n}(\Psi,\Psi_{h\t})^2 \right  ).
\end{align} 
Furthermore, we have the following global-in-space efficiency bound:
\begin{align}\label{eq:Y_norm_guaranteed_efficiency_lower_global}
&[\eta^n_{\mathrm{LB}}]^2:= \int_{I_n}([\eta^{\mathrm{F}}_{n,h,\Om}]^2 + [\eta^{\mathrm{J},H^1}_{n,h,\Om}]^2)\nonumber\\
& \lesssim  \int_{I_n} \left (\sum_{j\in \{\Gf,\bm{F},t\}}[\eta^{\mathrm{qd},j}_{n,h,\Om}]^2 + [\eta^{\mathrm{osc}}_{n,\Om}]^2+  \a^2[\eta^{\mathrm{J},L^2}_{n,h,\Om}]^2 +  [\eta^{\mathrm{J},H^1}_{n,h,\Om}]^2 \right ) + \dst^{\a}_{\Om,I_n}(\Psi,\Psi_{h\t})^2.
\end{align}
\end{theorem}

\begin{remark}[The linear heat equation case]\label{rem:LinearEquv}
\textbf{Quadrature:} In the absence of non-lineaities, $\etaOscF=\etaOscT=0$. To see this, note that if $S(p)=p$ and $\k$ is linear with respect to $s$, then the numerical solutions $\Psi_{n,h}$ and $S_{n,h}$ are in the same polynomial space as $p_{n,h}$. Thus, the quadrature terms above vanish. Moreover, $\eta^{\mathrm{qd},\Gf}_{n,h,K}=\frac{h_K}{\sqrt{K_{\mathrm{m}}} \pi}\|f(\cdot,t_n)-\Lam_{n,h}f(\cdot,t_n)\|$ becomes a data oscillation term. \textbf{Equivalence with estimates in \cite{Ern_Sme_Voh_heat_HO_Y_17}:} The bounds \eqref{eq:LocalEfficiencyMainEstimate}--\eqref{eq:Y_norm_guaranteed_efficiency_lower_global}
are equivalent to the efficiency bounds presented in \cite[Theorem 5.2]{Ern_Sme_Voh_heat_HO_Y_17} for the linear heat equation since $\int_{I_n} [\etaJh]^2$ is equivalent to $\int_{I_n} \|\del(\mathcal{I} u_{h\t}-u_{h\t})\|_K^2$ defined  in \cite[Section 5]{Ern_Sme_Voh_heat_HO_Y_17}. Additionally, in the linear case $\a=0$.
The grouping of terms in \eqref{eq:LocalEfficiencyMainEstimate}--\eqref{eq:Y_norm_guaranteed_efficiency_lower_global} is particularly useful since the quantity $(\sum_{n=1}^N (\int_{I_n }[\eta^{\mathrm{J},H^1}_{n,h,\Om}]^2 + \dst^{0}_{\Om,I_n}(\Psi,\Psi_{h\t})^2))^{\frac{1}{2}}$ directly relates to the $\|\Psi-\Psi_{h\t}\|_{\mathcal{E}_Y}$ error introduced in \cite[Section 5]{Ern_Sme_Voh_heat_HO_Y_17} which provides an estimate of  $\dst^{0}_{\Om,[0,T]}(\Psi,\Psi_{h\t})$ as proved in \cite[Theorem 5.1]{Ern_Sme_Voh_heat_HO_Y_17}. Thus, the $( \int_{I_n }[\eta^{\mathrm{J},H^1}_{n,h,\Om}]^2 + \dst^{0}_{\Om,I_n}(\Psi,\Psi_{h\t})^2)^{\frac{1}{2}}$ terms can themselves be considered error measures.  \label{rem:VanishingEta}
\end{remark}

\begin{proof} Observe from \eqref{eq:estimators}, \eqref{eq:estimators2}, and the definition of $\Fnt$ in \eqref{eq:SourceFluxTerms} that
\begin{align}\label{eq:ExpandingEtaEq}
[\etaEq]&\leq \|\K^{-\frac{1}{2}}\sth + \K^{\frac{1}{2}} \Fnt\|_K + \|\Psi_{h\t}-\Psi_{n,h}\|_{\Kh(K)} + \| \K^{\frac{1}{2}}\vg (\k(s_{h\t})-\k(S_{n,h}))\|_K \nonumber\\
&\overset{\eqref{eq:estimators}, \eqref{eq:estimators2}}\leq \left (\|\K^{-\frac{1}{2}}\sth + \K^{\frac{1}{2}}\PiH \Fnt\|_K + [\etaOscF]\right) + [\etaJh] + \a\, [\etaJl].
\end{align}
Note that, $[\etaOscF + \etaJh + \a\, \etaJl]\lesssim \sum_{\ta\in\calV_K}[\etaOscFa + \eta^{\mathrm{J},H^1}_{n,h,\oma} + \a\,\eta^{\mathrm{J},L^2}_{n,h,\oma}]$. 
For the first term on the right-hand side of \eqref{eq:ExpandingEtaEq}, one has
\begin{align}
&\|\K^{-\frac{1}{2}}\sth + \K^{\frac{1}{2}}\PiH \Fnt\|_K =\left\| \sum_{\ta\in \calV_K} (\K^{-\frac{1}{2}}\stha + \psia \K^{\frac{1}{2}}\,\PiH \Fnt)\right \|_{K} \nonumber\\
&\leq \sum_{\ta\in\calV_K} \|\K^{-\frac{1}{2}}\stha + \psia\K^{\frac{1}{2}}  \,\PiH\Fnt\|_{K}  \leq \sum_{\ta\in\calV_K} \|\K^{-\frac{1}{2}}\stha + \psia\K^{\frac{1}{2}}  \,\PiH\Fnt\|_{\oma}.\label{eq:WritingKinOma}
\end{align}
By denoting $R^\ta (\f):= (\Gf_{n,h},\f)_{\oma}-(\K\bm{F}_{n,h},\del \f)_{\oma}$, we apply Theorem 1.2 of \cite{Ern_Smears_Voh_H-1_lift_17} (also see Lemma 10 of \cite{Ern_Sme_Voh_heat_HO_Y_17}) to get from \eqref{eq:stha_minimization_def} that
\begin{align}\label{eq:flux_reconstruction_stability}
&\|\K^{-\frac{1}{2}}\stha + \psia \,\K^{\frac{1}{2}} \PiH\Fnt\|_\oma\lesssim \sup_{\scriptscriptstyle{\f \in H^1_0(\oma),\; \|\f\|_{\Kh(\oma)}=1}} [(\Lam_{n,h}\Snt,\f)_{\oma}-(\K\,\PiH\Fnt,\del \f)_{\oma}],\nonumber\\
&=  \sup_{\scriptscriptstyle{\f \in H^1_0(\oma),\; \|\f\|_{\Kh(\oma)}=1}} \left [(\Lam_{n,h} \Snt - \Snt,\f)_{\oma}-(\K(\PiH\Fnt- \Fnt),\del \f)_{\oma} + R^\ta (\f)\right ]\nonumber\\
&\overset{\eqref{eq:DefPiL}}=  \sup_{\substack{\scriptscriptstyle{\f \in H^1_0(\oma),}\\
\scriptscriptstyle{\|\f\|_{\Kh(\oma)}=1}}} \left [ \sum_{\elCR\in \Ta}\left (\Lam_{n,h}\Snt - \Snt,\f - \tfrac{1}{|\elCR|}\smallint_{\elCR}\f\right )_{\elCR}-(\K(\PiH\Fnt- \Fnt),\del \f)_{\oma} + R^\ta (\f) \right ]\nonumber\\
&\lesssim \etaOscSa + \etaOscFa + \sup_{\substack{\scriptscriptstyle{\f \in H^1_0(\oma),}\,
\scriptscriptstyle{\|\f\|_{\Kh(\oma)}=1}}} R^\ta (\f).
\end{align}
Focusing on the final term, and recalling \eqref{eq:SourceFluxTerms},  one obtains for $\|\f\|_{\Kh(\oma)}=1$ and $t\in I_n$ that
\begin{align}
&R^\ta (\f)=(\Gf_{n,h},\f)_\oma-(\K\bm{F}_{n,h},\del \f)_\oma\nonumber\\
&=(f(S_{n,h},\bm{x},t_n)-\tfrac{1}{\t_n}(S_{n,h}-S_{n-1,h}),\f)_\oma-(\K[\nabla  \Psi_{n,h} +\vg \k(S_{n,h})],\del \f)_\oma\nonumber\\
&\overset{\eqref{eq:DefRhtau}}=\langle \res(\Psi_{h\t}),\f\rangle + (f(S_{n,h},\bm{x},t_n)-f(s_{h\t},\bm{x},t),\f)_\oma  + (\p_t s_{h\t} -\tfrac{1}{\t_n}(S_{n,h}-S_{n-1,h}) ,\f)_\oma  \nonumber\\
&\qquad + (\K[\nabla (\Psi_{h\t}- \Psi_{n,h}) +\vg (\k(s_{h\t})-\k(S_{n,h}))],\del \f)_\oma\nonumber\\
&\overset{\eqref{eq:estimators}}\leq \|\res(\Psi_{h\t})\|_{\Kn(\oma)} +  \eta^{\mathrm{osc}}_{n,\oma} + \etaOscTa + \eta^{\mathrm{J},H^1}_{n,h,\oma} + \a\,\eta^{\mathrm{J},L^2}_{n,h,\oma}.\label{eq:ResidualEstimateForLowerBound}
\end{align}
Recall from \Cref{theo:TotalPressureEfficiency} that $\int_{I_n}\|\res(\Psi_{h\t})\|_{\Kn(\oma)}^2 \leq \dst^{\a}_{\oma,I_n}(\Psi,\Psi_{h\t})^2$. Thus
combining \eqref{eq:ExpandingEtaEq}--\eqref{eq:ResidualEstimateForLowerBound}, squaring both sides,  and integrating over $I_n$, we have \eqref{eq:LocalEfficiencyMainEstimate}. To get the global efficiency bound \eqref{eq:Y_norm_guaranteed_efficiency_lower_global}
 we sum \eqref{eq:LocalEfficiencyMainEstimate} over all mesh elements and note that $\sum_{\ta\in\calV^n}\|u\|^2_{\oma}\lesssim \|u\|^2$ and $\sum_{\ta\in\calV_n}\|u\|^2_{H^{\pm 1}_{\K}(\oma)}\lesssim \|u\|^2_{H^{\pm 1}_{\K}(\Om)}$, see \cite[Lemma 3.5]{Coh_DeVore_Noch_cvg_Hmo_12}. 
\end{proof}

\section{Numerical results}\label{sec:Numerical}
We choose the unit square $\Om=(0,1)^2$ as the simulation domain, $T=1$ as the final time, and both uniform and non-uniform triangulations $\calT_n$ with the discretization levels: 
\begin{equation}\label{eq:DiscrtDetail}
(h,\t)=(h_0,\t_0)/\ell \; \text{ where } \ell\in \{1,2,4\},\, h_0=0.2,\, \t_0=0.04.
\end{equation}
The mesh and the time-step size remain fixed between time steps.  Piecewise linear finite elements are used for obtaining  the solutions, i.e., $\mathfrak{p}_n=1$ in \Cref{sec:fem_mesh}. Iterative linearization is discussed in Appendix \ref{app:num_linear}.

We consider the following three test cases:
\begin{itemize}
\item \Cref{sec:NumTestNonDeg}: Nonlinear but nondegenerate problem with known exact solution.
\item \Cref{sec:NumTestDeg}: Nonlinear and degenerate problem  in the total pressure formulation \eqref{eq:Richards} with known exact solution.
\item \Cref{sec:NumTestReal}: Realistic case, nonlinear, degenerate with heterogeneous and anisotropic $\K$, mixed boundary conditions (Neumann + Dirichlet), discontinuous initial condition, non-uniform mesh, and no known exact solution.
\end{itemize}
The code is implemented in FreeFem++ and can be accessed through this \href{https://github.com/koondax/Aposteriori.git}{\textcolor{blue}{link}}.

\begin{remark}[Choice of $\lambda$ in \Cref{theo:UpperBound,theo:GlobalReliability}]\label{rem:lambda}
The choice of $\lambda:[0,T]\to \R^+$ in \eqref{eq:UpperBoundEstimate1st} is important in our simulations since in \eqref{eq:YoungOnResLambda}, the $\|\calR\|_{\Kn(\Om)}$ term is much larger than $\|G^0_{h\t}\|_{\Kh(\Om)}=\|s-s_{h\t}\|_{\Kn(\Om)}$. Hence, choosing  $\lambda=1$ leads to a significant overestimation of the error. Here, we have used $\lambda=200$ in \Cref{sec:NumTestNonDeg} and $\lambda=100$ in  \Cref{sec:NumTestDeg}. These yield close to minimum values of the effectivity indices defined in \eqref{eq:EffInd_UB}. The optimal value of  $\l$ can also be roughly estimated by minimizing the  right hand side of the Young's inequality in \eqref{eq:YoungOnResLambda}. This gives  $\lambda\sim \|\calR\|_{\Kn(\Om)}\slash \|s-s_{h\t}\|_{\Kn(\Om)}\sim \eta_{\mathcal{R}}\slash \eta^{\mathrm{ini},H^{-1}}$ (see \eqref{eq:estimators}) which yields $\lambda$ in the same order of magnitude as the $\lambda$ chosen in our simulations.
\end{remark}

\subsection{Nonlinear nondegenerate case with known solution}\label{sec:NumTestNonDeg}
\begin{figure}[h!]
\begin{subfigure}{.35\textwidth}
\includegraphics[width=\textwidth]{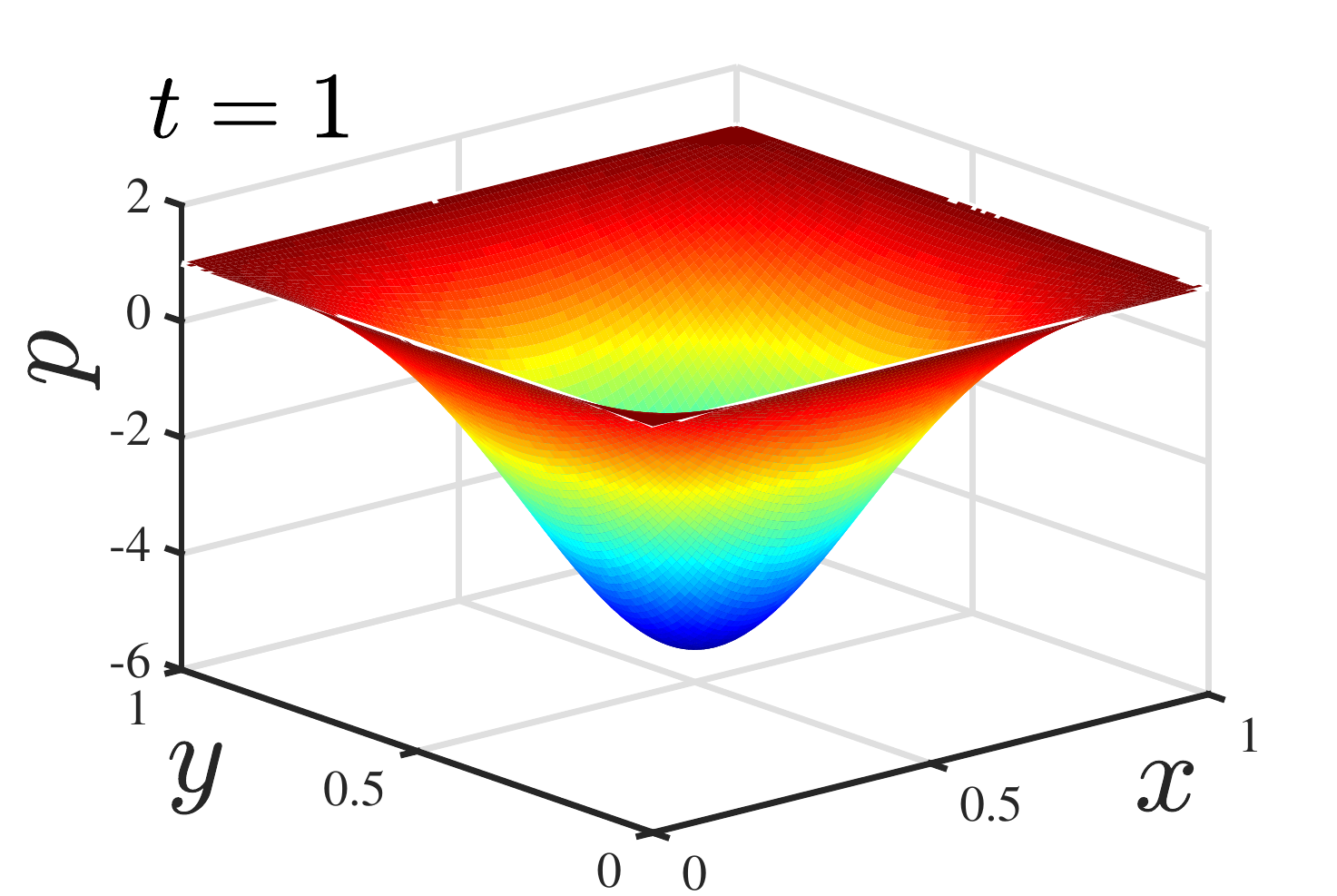}
\end{subfigure}
\begin{subfigure}{.35\textwidth}
\includegraphics[width=\textwidth]{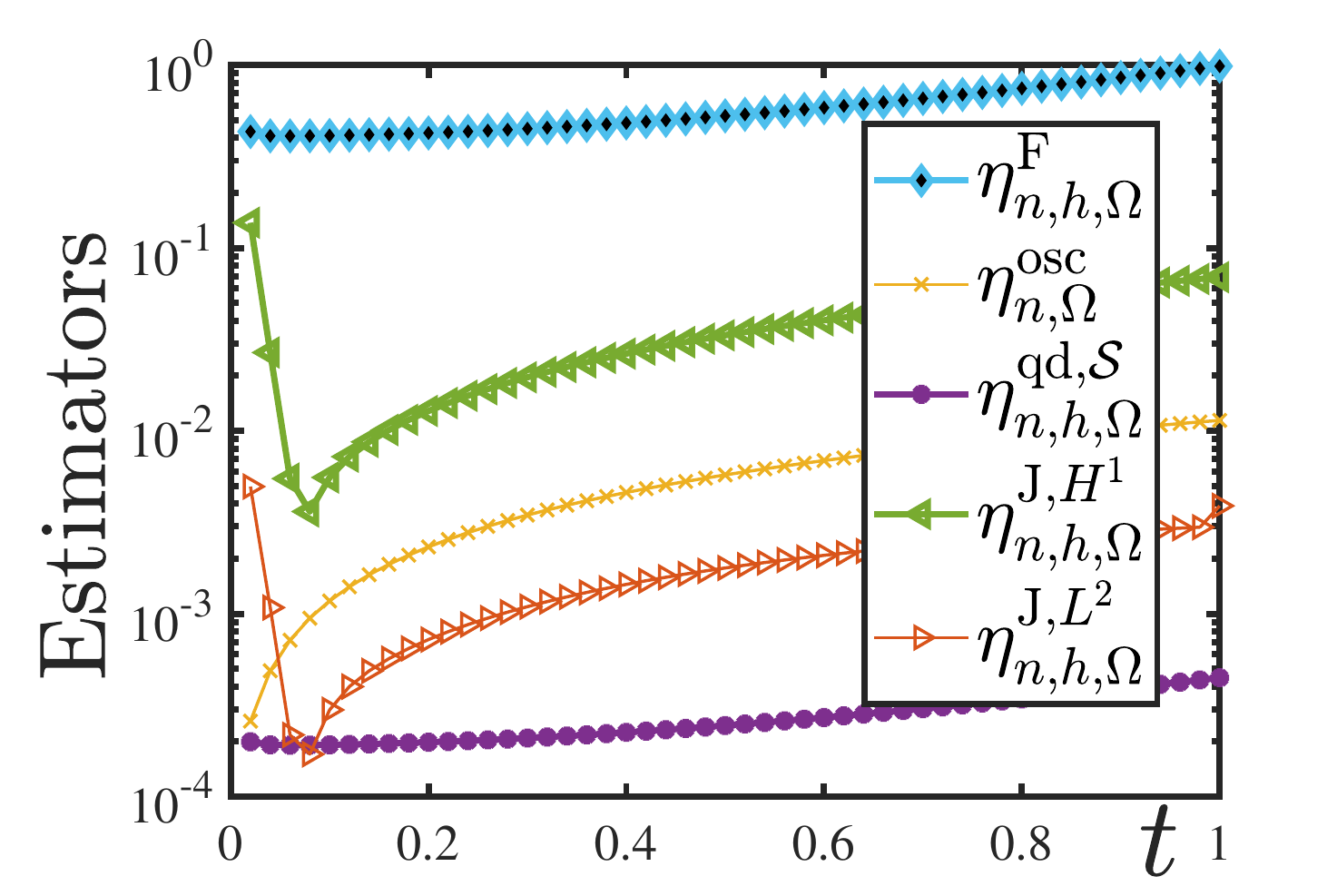}
\end{subfigure}
\begin{subfigure}{.28\textwidth}
\includegraphics[width=\textwidth]{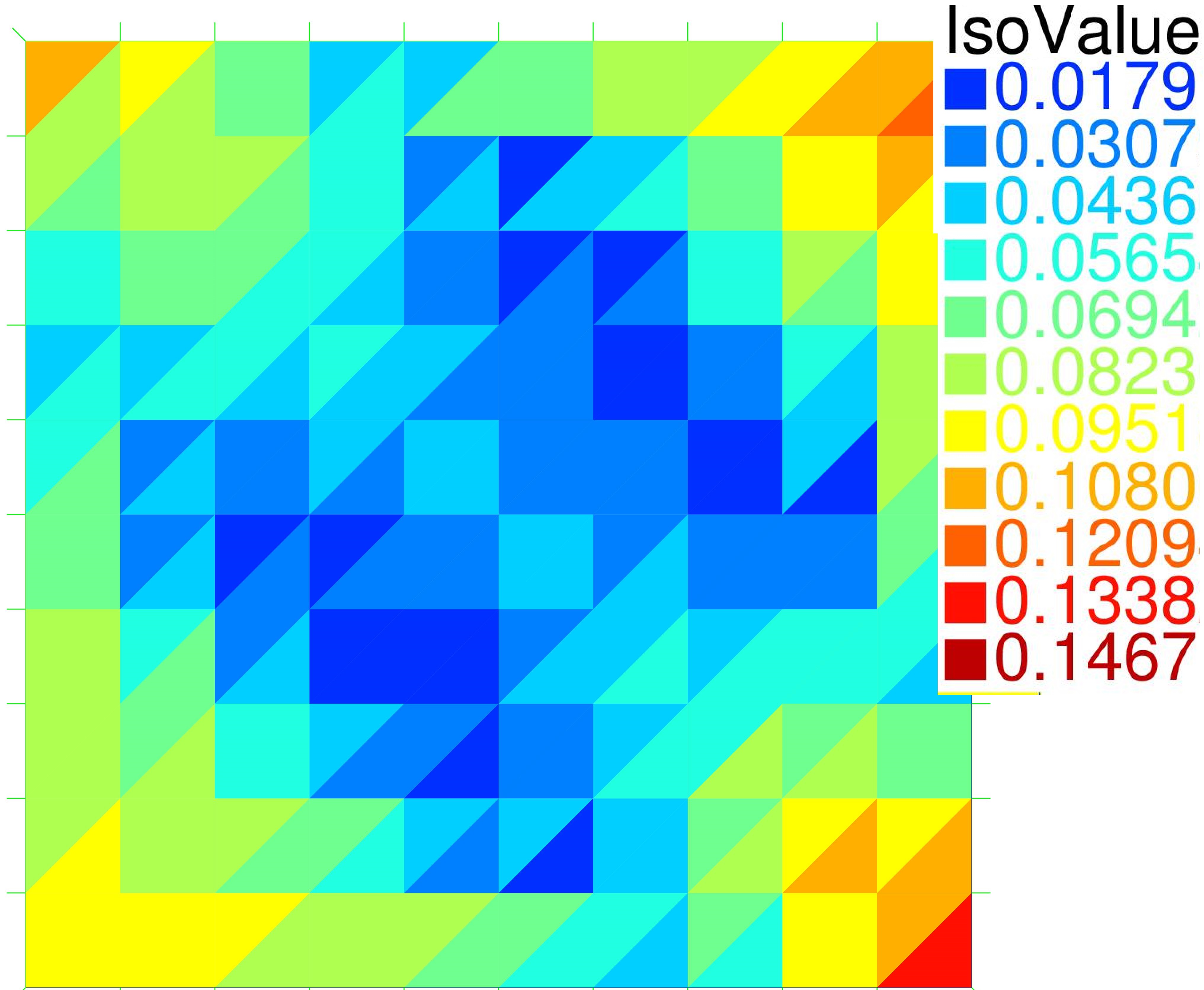}
\end{subfigure}
\caption{[\Cref{sec:NumTestNonDeg}] Exact solution $p_{\mathrm{exact}}$ of \eqref{eq:exactsolNonDeg} at time $t=1$ (left). Evolution of the 5 most significant estimators for $\ell=2$  (center). The elementwise flux estimators $\etaEq$ for $\ell=2$ and $t_n=1$ (right).}\label{fig:Dom&Estimators}
\end{figure}

For this case, $\K=\mathbb{I}$, and with $\hat{\bm{e}}_x$ representing the unit vector along $x$-axis, we specify
\begin{align}\label{eq:NonlinearityTest1}
 \vg=-\hat{\bm{e}}_x,\quad \k(s)=s^3, \text{ and } S(p)=\begin{cases}
\frac{1}{(2-p)^{\frac{1}{3}}} &\text{ if } p<1,\\
1 &\text{ if } p\geq 1.
\end{cases}
\end{align}
These nonlinearities resemble the Brooks--Corey parametrization \eqref{eq:BrooksCorey}. An exact solution  
\begin{align}\label{eq:exactsolNonDeg}
p_{\mathrm{exact}}(x,y,t)=2- \exp(16\,(1+t^2)\, x\,y\, (1-x)\, (1-y))
\end{align}
is fixed, see \Cref{fig:Dom&Estimators} (left).  The source term $f$ is independent of $s$, and is adjusted together with the initial condition $s_0$, and the inhomogeneous Dirichlet boundary condition so that $p_{\mathrm{exact}}$ indeed solves \eqref{eq:RichardsP}.

 Evolution of the different estimators for the case $\ell=2$ is presented in  \Cref{fig:Dom&Estimators} (center), which shows that $\etaEq$ is the dominant estimator followed by $\eta^{\mathrm{J},H^1}_{n,h,\Om}(t)$ for this test. The spatial distribution of $\etaEq$ is shown in \Cref{fig:Dom&Estimators} (right). 
 The time-quadrature and the degeneracy estimators, $\etaOscTom$  and $\eta^{\mathrm{deg}}$, vanish in this case as the problem is nondegenerate. 

\begin{figure}[h!]
\begin{subfigure}{.48\textwidth}
\includegraphics[width=.9\textwidth]{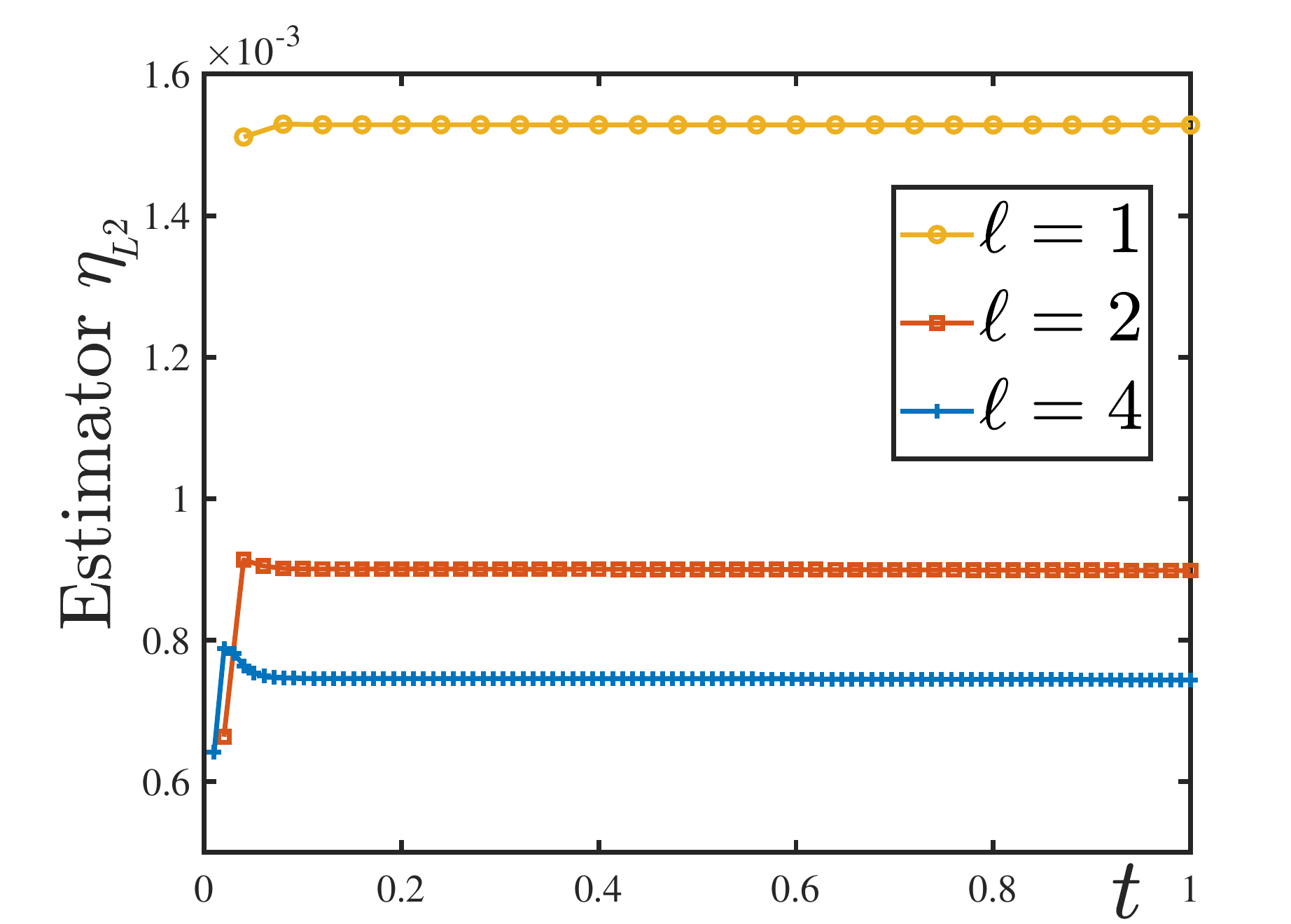}
\end{subfigure}
\begin{subfigure}{.48\textwidth}
\includegraphics[width=.9\textwidth]{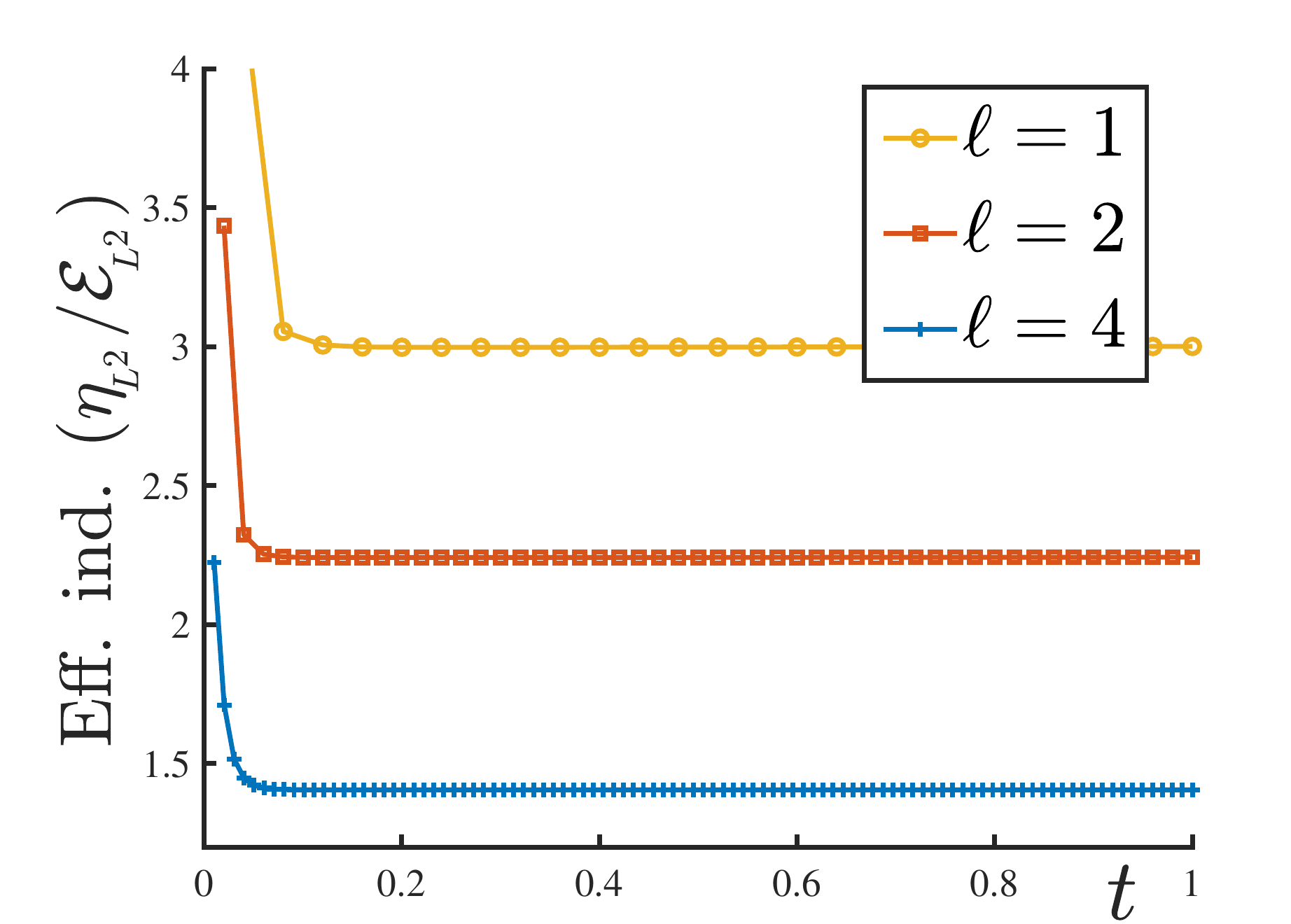}
\end{subfigure}
\caption{[\Cref{sec:NumTestNonDeg}] Estimator $\eta_{L^2}$  of  \eqref{eq:EL2etaL2}   as a function of the final time (left), and the corresponding effectivity index (right). }\label{fig:L2errorTest1}
\end{figure}

\begin{figure}[h!]
\begin{subfigure}{.48\textwidth}
\includegraphics[width=.9\textwidth]{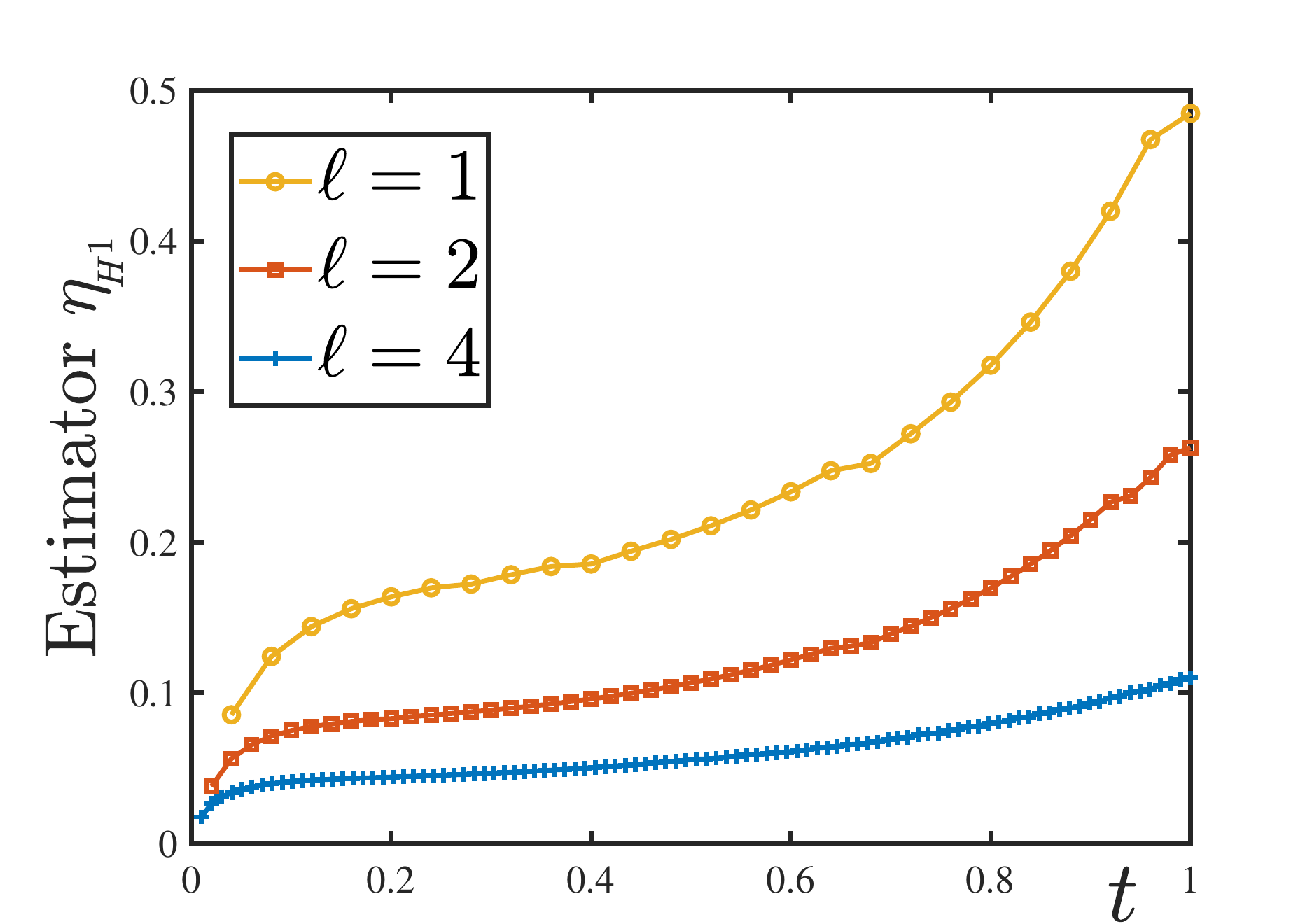}
\end{subfigure}
\begin{subfigure}{.48\textwidth}
\includegraphics[width=.9\textwidth]{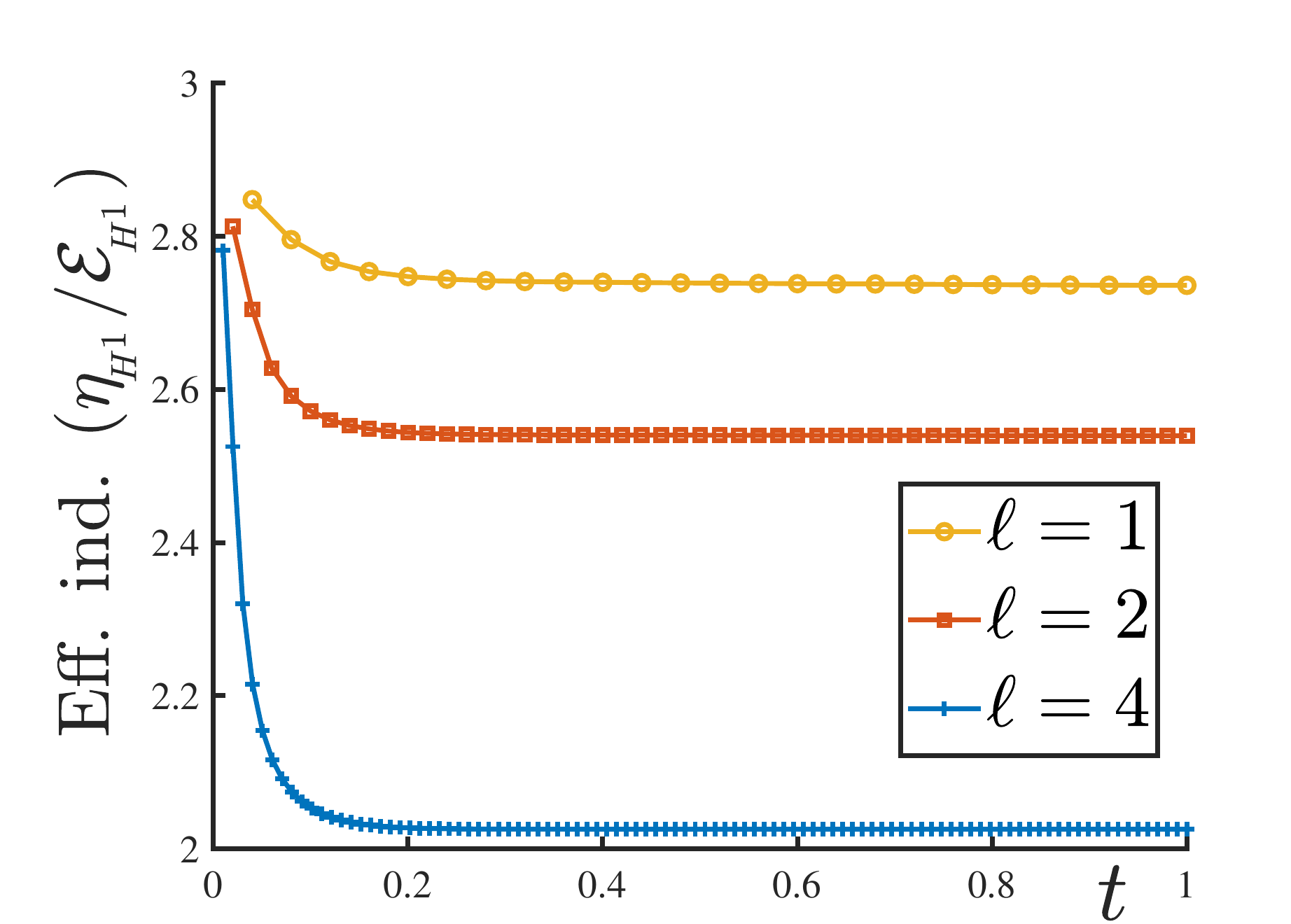}
\end{subfigure}
\caption{[\Cref{sec:NumTestNonDeg}] Estimator $\eta_{H^1}$ of  \eqref{eq:EH1etaH1}   as a function of the final time (left), and the corresponding effectivity index (right).}\label{fig:H1errorTest1}
\end{figure}

Next, we numerically investigate the quality of the upper bound  from \Cref{theo:GlobalReliability}. For this purpose, we introduce the effectivity index defined as
\begin{align}\label{eq:EffInd_UB}
\text{effectivity index}:=\text{upper bound}\slash \text{error}=\begin{cases}
\eta_{L^2}\slash \calE_{L^2},\\
\eta_{H^1}\slash \calE_{H^1}.
\end{cases}
\end{align}
Effectivity index close to 1 is desirable. \Cref{fig:L2errorTest1} (left) shows the evolution of $\eta_{L^2}$ (see \eqref{eq:EL2etaL2}) as a function of time and discretization level $\ell$. The upper bound $\eta_{L^2}$ reaches a constant state after an initial transition period. This is since $\mathfrak{C}_1(t)$ is almost constant for this case, and the error $\|\res\|_{\Kn(\Om)}$ increases exponentially with a rate much smaller than $\lambda + \mathfrak{C}_1$. Hence, a near constant $\eta_{L^2}$ is expected from carrying out the integrals.  The (right) plot shows the  effectivity indices of $\eta_{L^2}$. The effectivity varies between  1.4 and 3.1, and improves with the discretization level $\ell$.  \Cref{fig:H1errorTest1} is the same plot presented for $\eta_{H^1}$. The estimator $\eta_{H^1}$ increases with $t$ as $\mathfrak{C}_2$ increases rapidly with time. The effectivity index again improves as the discretization is refined. 

\begin{figure}[h!]
\begin{subfigure}{.48\textwidth}
\includegraphics[width=.9\textwidth]{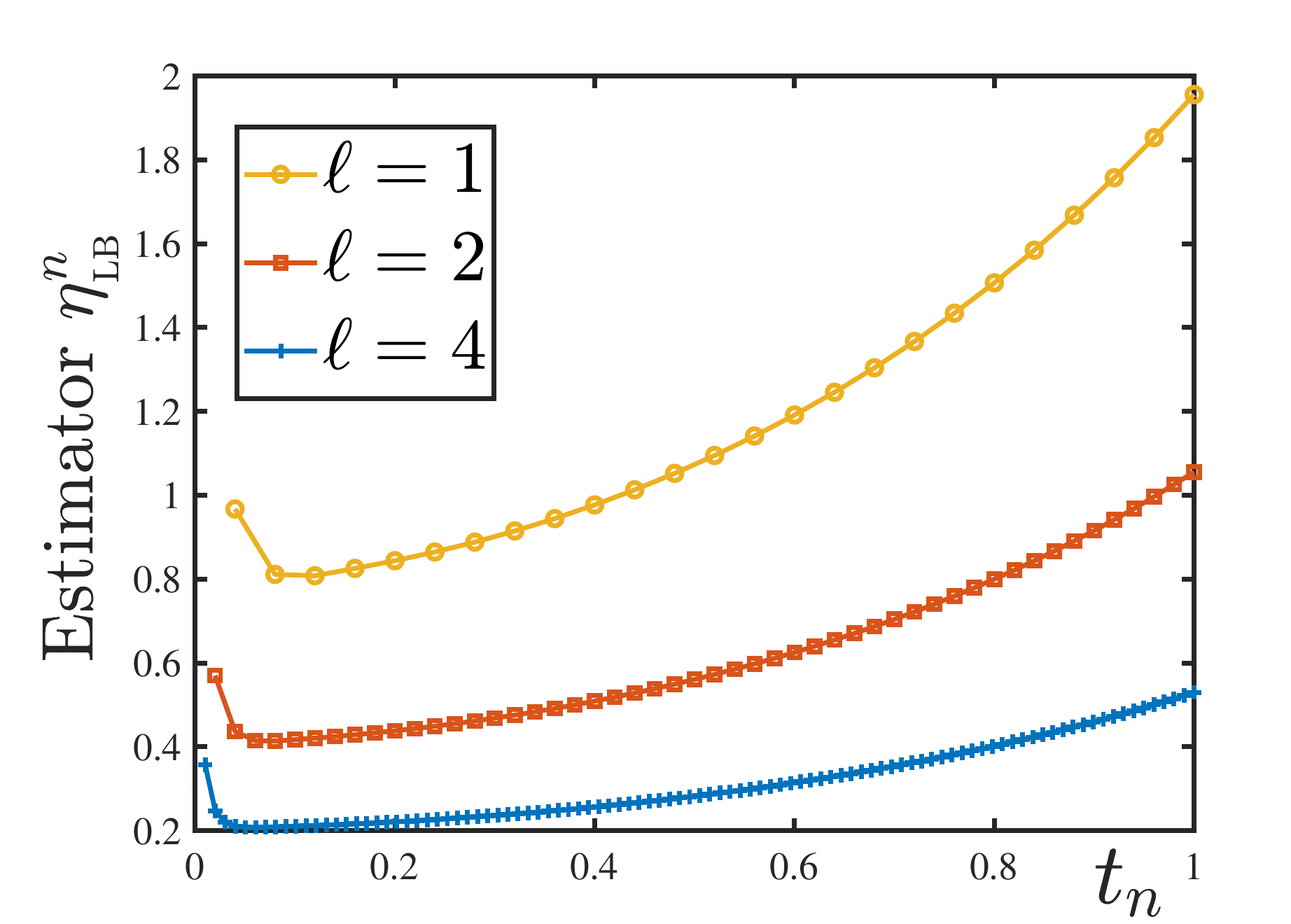}
\end{subfigure}
\begin{subfigure}{.48\textwidth}
\includegraphics[width=.9\textwidth]{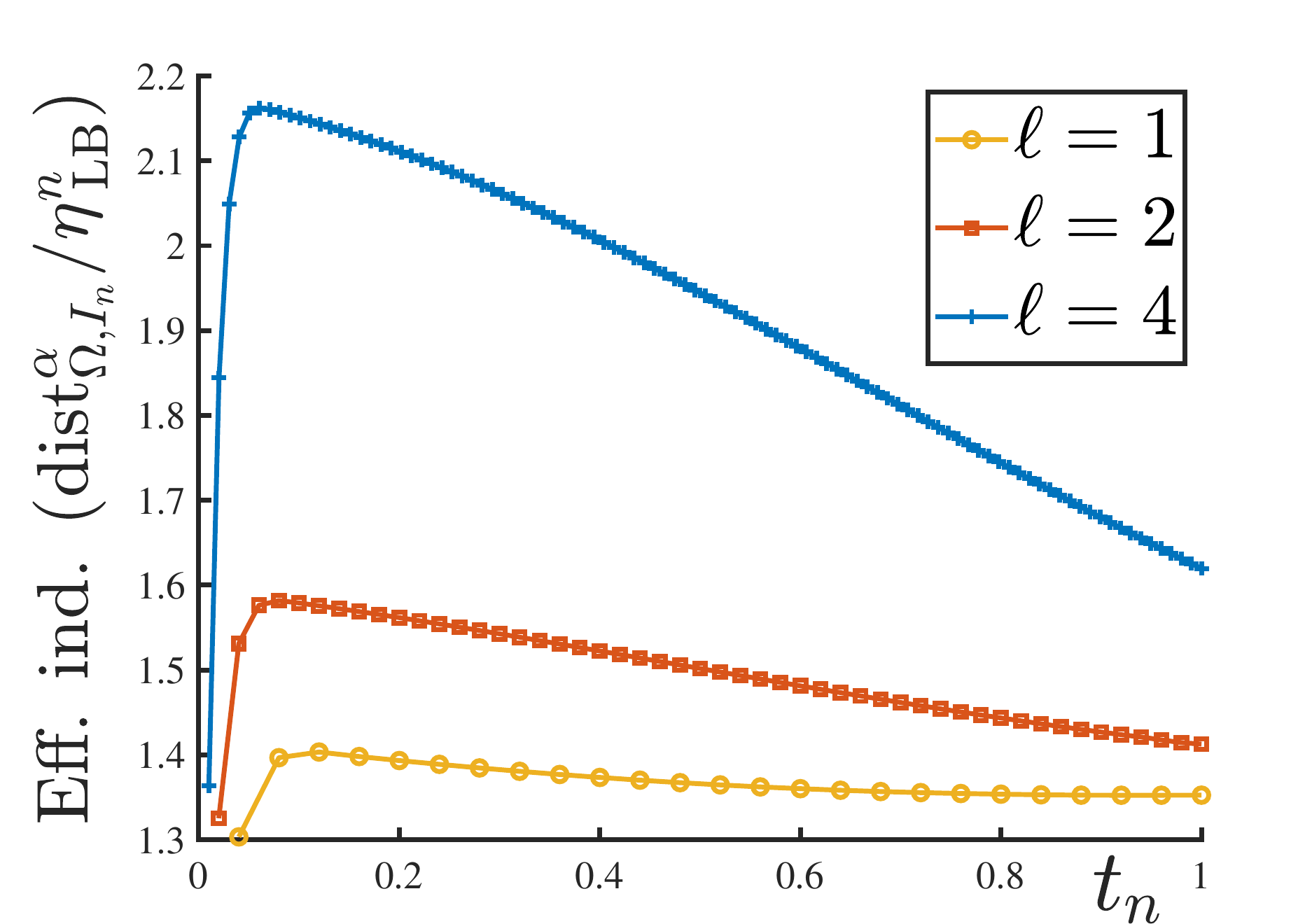}
\end{subfigure}
\caption{[\Cref{sec:NumTestNonDeg}] Estimator $\eta^n_{\mathrm{LB}}$ of  \eqref{eq:Y_norm_guaranteed_efficiency_lower_global} as a function of $t_n$ (left). Its effectivity indices computed using \eqref{eq:EffInd_LB} (right). }\label{fig:LBTest1}
\end{figure}

\begin{figure}[h!]
\centering
\includegraphics[scale=.35]{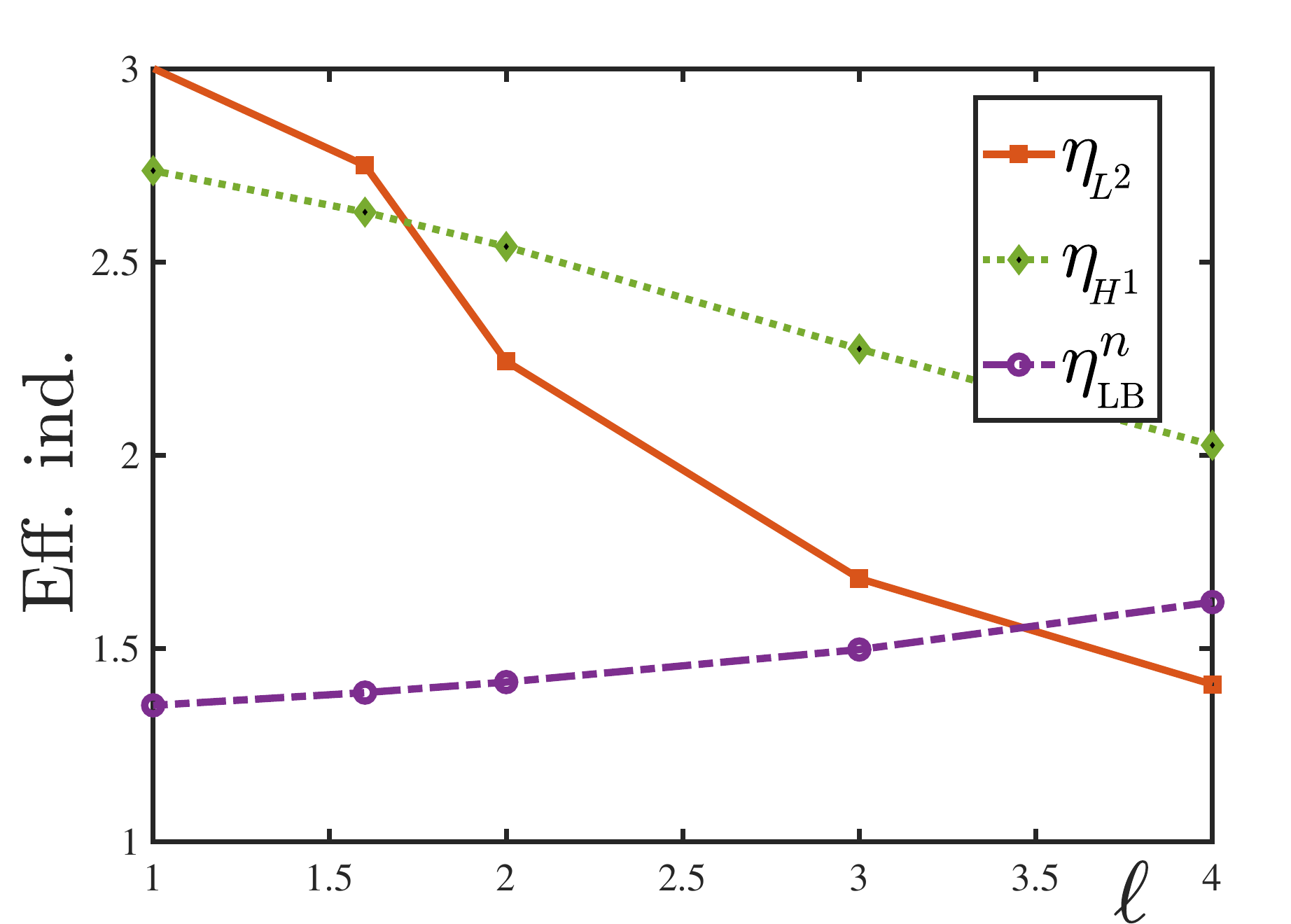}
\caption{[\Cref{sec:NumTestNonDeg}] Effectivity indices  of the estimators $\eta_{L^2}$, $\eta_{H^1}$, and $\eta^n_{\mathrm{LB}}$ at the final time $T=1$ varying with $\ell$. }\label{fig:Discretization}
\end{figure}

We now turn to the lower bound of \Cref{theo:Y_norm_guaranteed_efficiency}. The effectivity index in  this context is computed as 
\begin{align}\label{eq:EffInd_LB}
\text{effectivity index}:= \text{error}\slash \text{lower bound}= \mathrm{dist}^{\a}_{\Om,I_n}(\Psi,\Psi_{h\t}) \slash \eta^n_{\mathrm{LB}},
\end{align}
where $\eta^n_{\mathrm{LB}}$ is given by \eqref{eq:Y_norm_guaranteed_efficiency_lower_global}. The reversed order is to make the effectivity index comparable to the effectivity index of the upper bound.  \Cref{fig:LBTest1} shows the lower bound  $\eta^n_{\mathrm{LB}}$ and its effectivity indices. The effectivity increases with $\ell$ in this case, though only varying between $1.4$ and $2.2$. A  higher effectivity index close to $t=0$ is also observed. This is explained by the fact that the lower bound estimator $\eta^n_{\mathrm{LB}}$  does not incorporate the initial errors, and thus, is more susceptible to inaccuracies close to $t=0$. \Cref{fig:Discretization} shows the variation of the effectivity indices with $\ell$ at the final time $T=1$, for both the reliability and the efficiency estimates.

\begin{figure}[h!]
\begin{subfigure}{.32\textwidth}
\includegraphics[width=\textwidth]{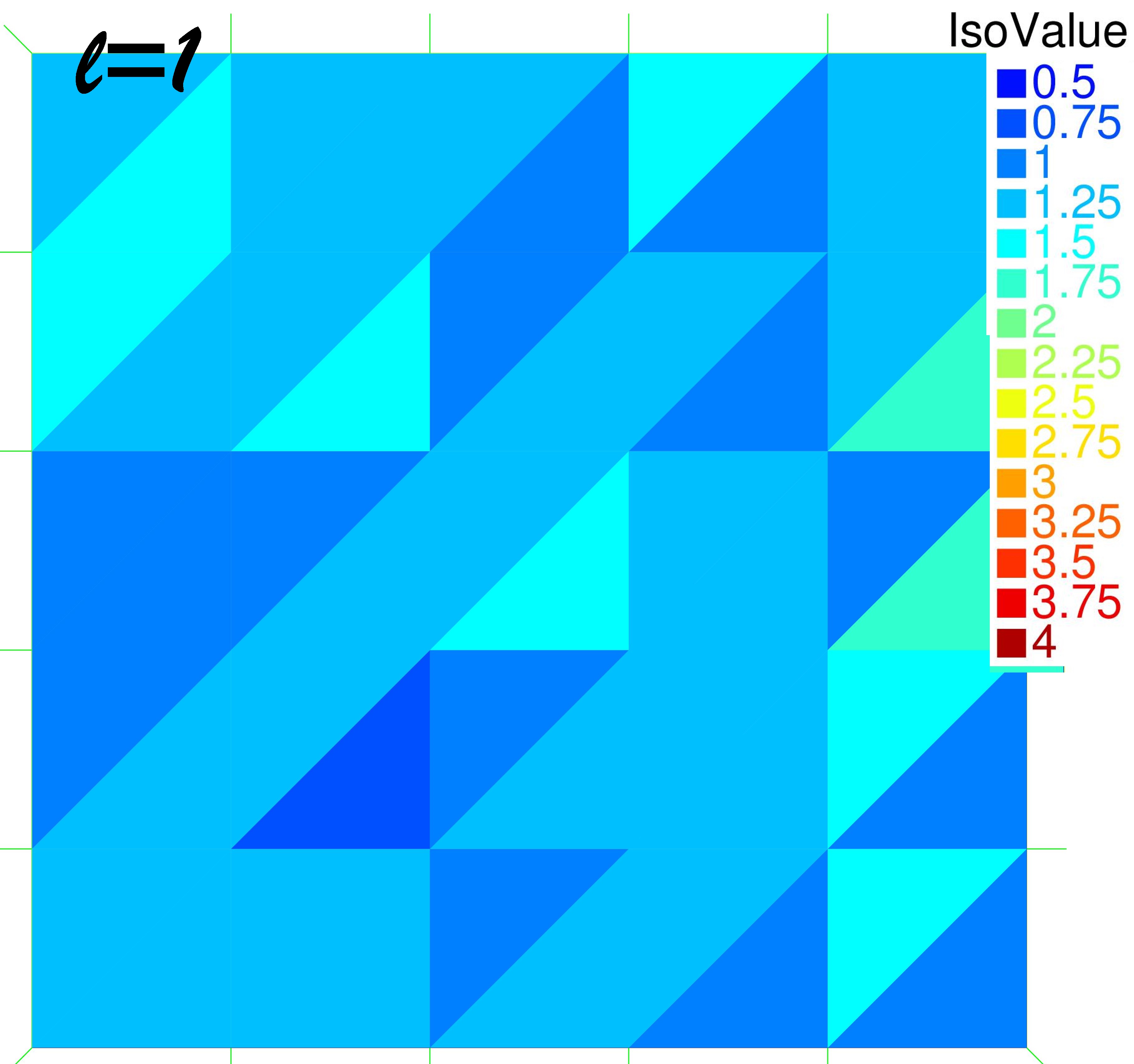}
\end{subfigure}
\begin{subfigure}{.32\textwidth}
\includegraphics[width=\textwidth]{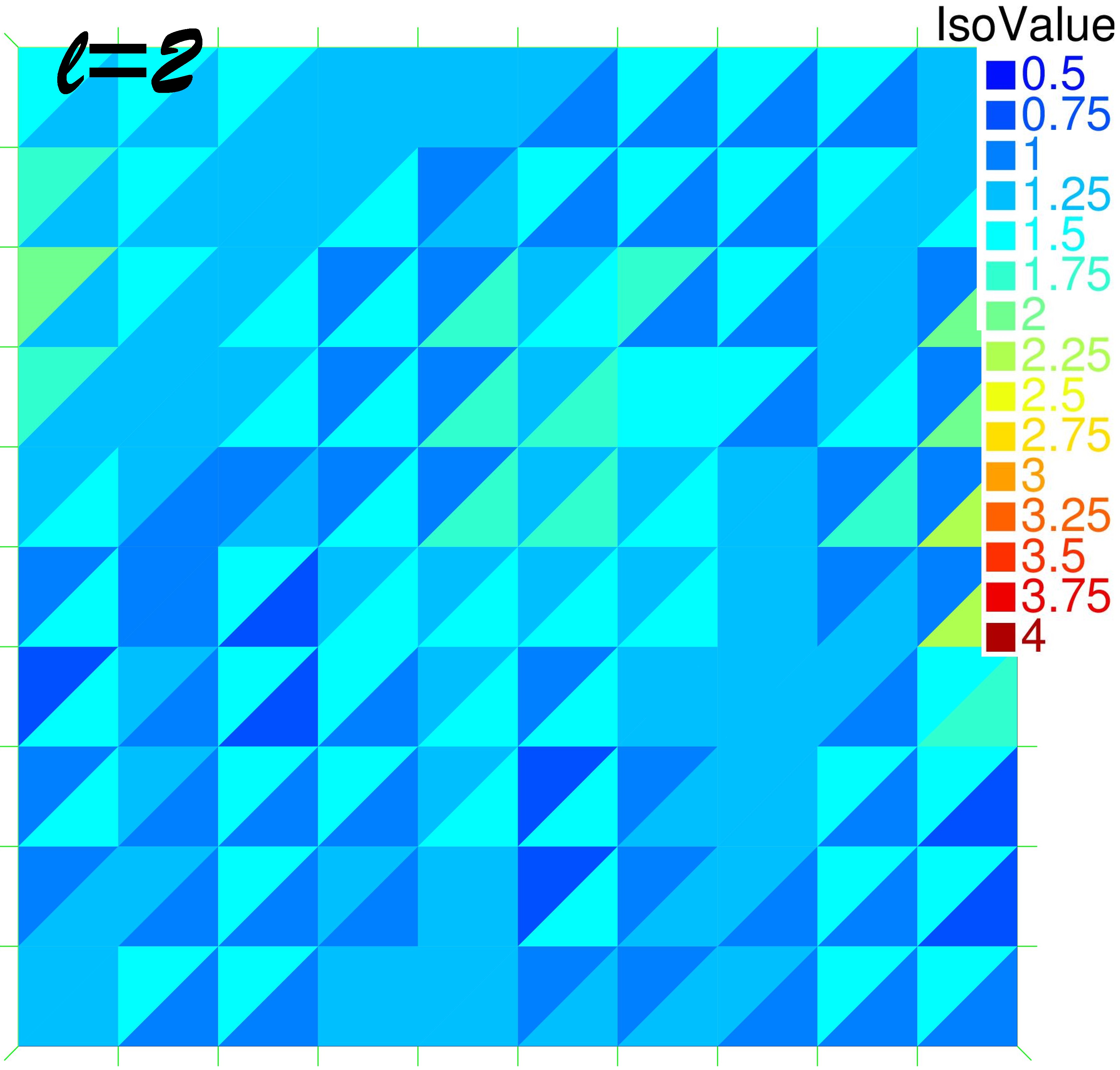}
\end{subfigure}
\begin{subfigure}{.32\textwidth}
\includegraphics[width=\textwidth]{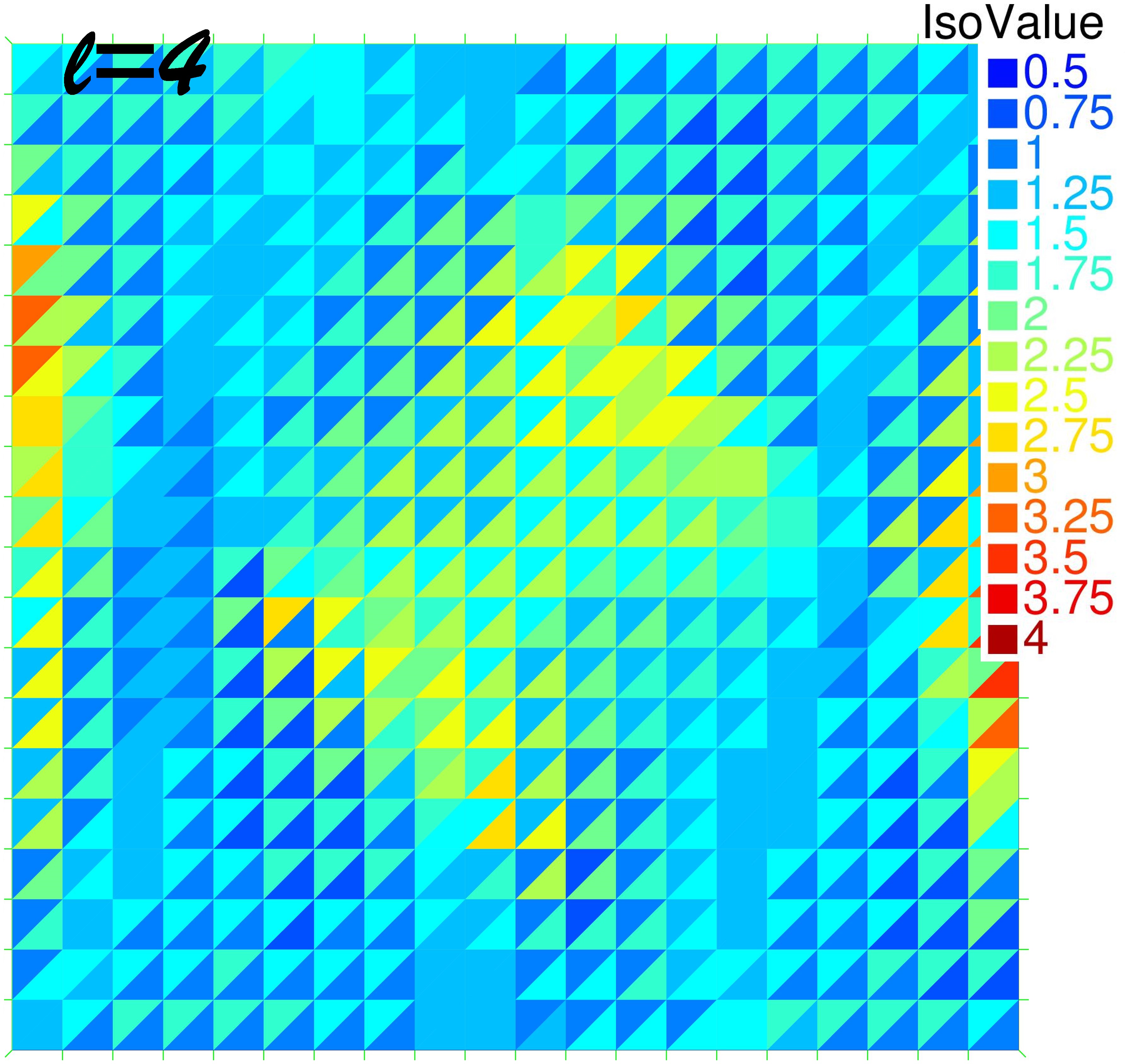}
\end{subfigure}
\caption{[\Cref{sec:NumTestNonDeg}] Local effectivity indices \eqref{eq:EffInd_Loc} at the final time $T=1$ for different values of $\ell$.}\label{fig:LocalEfficiency}
\end{figure}

Inspired by \Cref{theo:Y_norm_guaranteed_efficiency},  the local-in-space and in-time  effectivity indices are computed as
\begin{align}\label{eq:EffInd_Loc}
\text{(eff. ind.)}_{n,K}:=\dst^{\a}_{K,I_n}(\Psi,\Psi_{h\t})\slash \left (\smallint\!_{I_n} ([\etaEq]^2 + [\etaJh]^2)\right )^{\frac{1}{2}},
\end{align}
for all $K\in\calT_n$. From \Cref{fig:Dom&Estimators} (right), it is observed that $\etaEq$ varies with the mesh elements $K$ by a factor of about 10. However, \Cref{fig:LocalEfficiency} shows that the local effectivity indices are in the range $0.6$--$1.8$ for $\ell=1$,  $0.8$--$2.4$ for $\ell=2$, and $0.8$--$3.8$ for $\ell=4$, which we consider excellent. Observe that, $\text{(eff. ind.)}_{n,K}<1$ does not violate \Cref{theo:Y_norm_guaranteed_efficiency} since the error $\dst^{\a}_{\oma,I_n}(\Psi,\Psi_{h\t})$ and the sign `$\lesssim$' (up to a constant) was used there.

\subsection{Nonlinear degenerate case with known solution}\label{sec:NumTestDeg}
\begin{figure}[h!]
\begin{subfigure}{.32\textwidth}
\includegraphics[width=\textwidth]{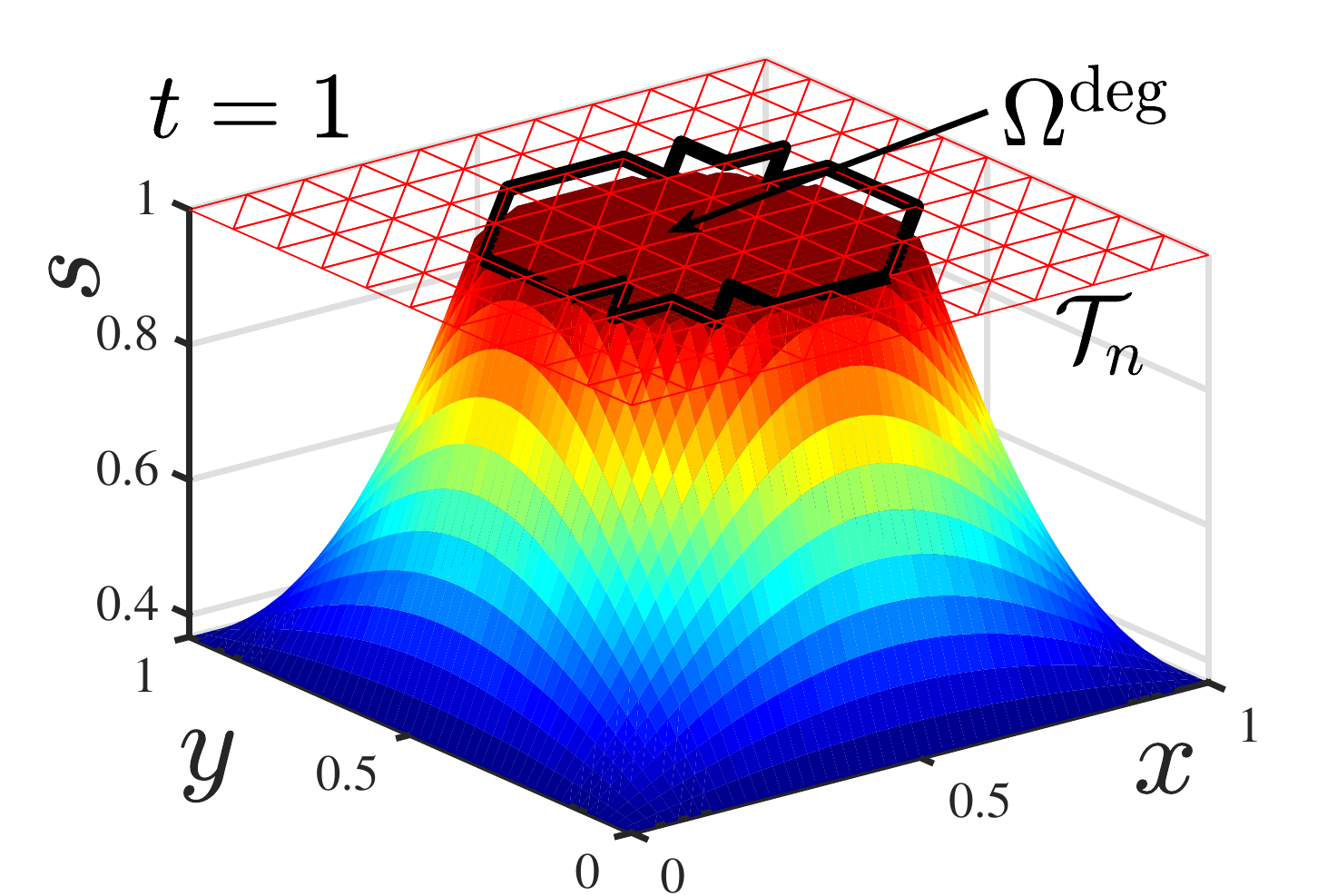}
\end{subfigure}
\begin{subfigure}{.36\textwidth}
\includegraphics[width=\textwidth]{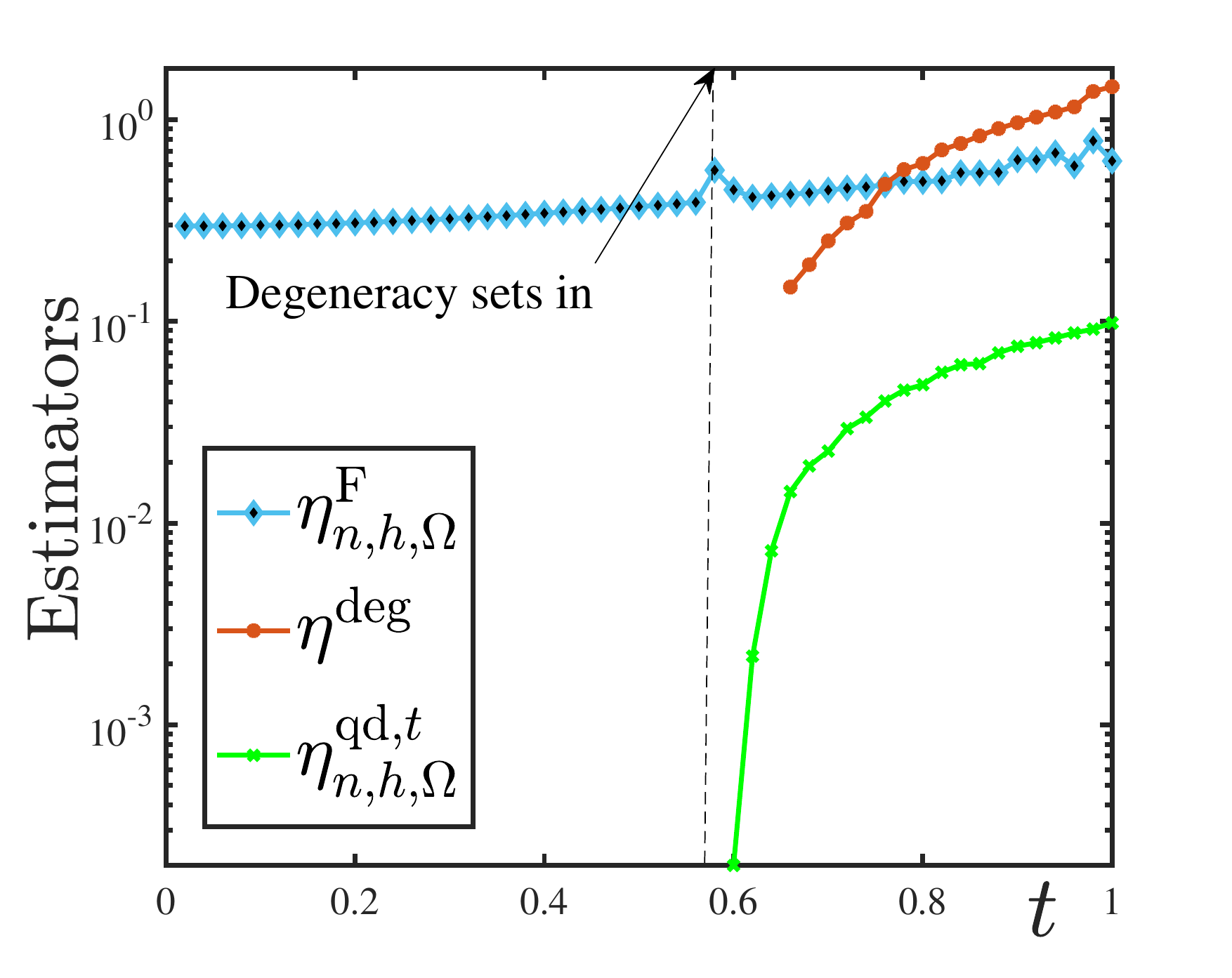}
\end{subfigure}
\begin{subfigure}{.28\textwidth}
\includegraphics[width=\textwidth]{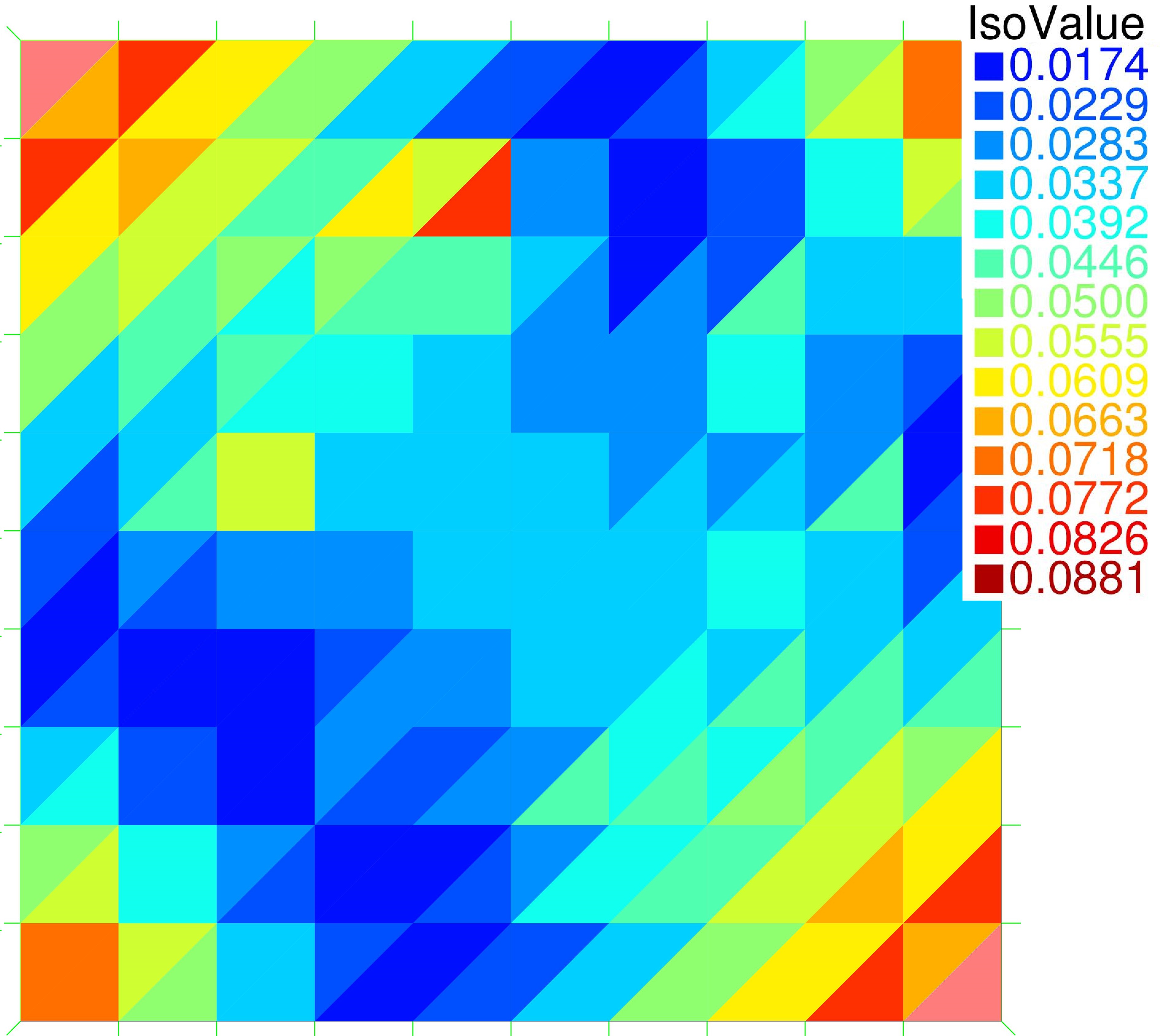}
\end{subfigure}
\caption{[\Cref{sec:NumTestDeg}] Saturation of the exact solution $\Psi_{\mathrm{exact}}$ and the domain $\Om^{\mathrm{deg}}(t)$ at $t=1$ (left). The principal estimators $\eta^{\mathrm{F}}_{n,h,\Om}(t)$,  $\etaDeg(t)$, and $\etaOscTom(t)$ for $\ell=2$ (center), and the elementwise estimators $\etaEq$ at $t_n=1$ (right).}\label{fig:DegDom}
\end{figure}
This test on purpose focuses on a system where degeneracy is the dominant effect. For a change, the total pressure formulation \eqref{eq:Richards} is used here. The nonlinearities are set as
 \begin{align}\label{eq:NonlinearityTest2}
\k(s)=s, \quad \Sf(\Psi)=\begin{cases}
\exp({\Psi-1}) &\text{ if } \Psi<1,\\
1 &\text{ if } \Psi\geq 1,
\end{cases}
\end{align}
and $\K=\mathbb{I}$, with the exact solution used being 
\begin{align}
\Psi_{\mathrm{exact}}(x,y,t)= 12\, (1+t^2)\, x\, y\, (1-x)\,(1-y).
\end{align}
Appropriate source function $f$ (independent of $s$), initial and boundary conditions are again imposed.

The solution is initially nondegenerate and contains a degenerate region after $t=0.58$, see \Cref{fig:DegDom} (left). This is caused by the source term $f$ since $\K$ here is uniform. The domain $\Om^{\mathrm{deg}}(t)$ is approximately computed for $t\in I_n$ as
\begin{align}
\Om^{\mathrm{deg}}(t)=\cup\{K\in\calT_n: K\cap \{\Psi_{h\t}(t)>P_{\mathrm{M}}=1\}\not= \emptyset\},
\end{align}
pointed out in \Cref{fig:DegDom} (left). This replaces here the generally unknown $\Om^{\mathrm{deg}}$ from \Cref{theo:UpperBound}. 
\Cref{fig:DegDom} (center) shows the estimators $\eta^{\mathrm{F}}_{n,h,\Om}$, $\etaDeg$, and $\etaOscTom(t)$ for the case $\ell=2$. The degeneracy estimator $\etaDeg$ defined in \Cref{theo:UpperBound} quickly rises in value as degeneracy sets in. In \Cref{fig:DegDom} (right) we see the distribution of $\etaEq$. The flux estimator $\etaEq$ stays relatively unaffected by the onset of degeneracy.

\begin{figure}[h!]
\begin{subfigure}{.48\textwidth}
\includegraphics[width=.8\textwidth]{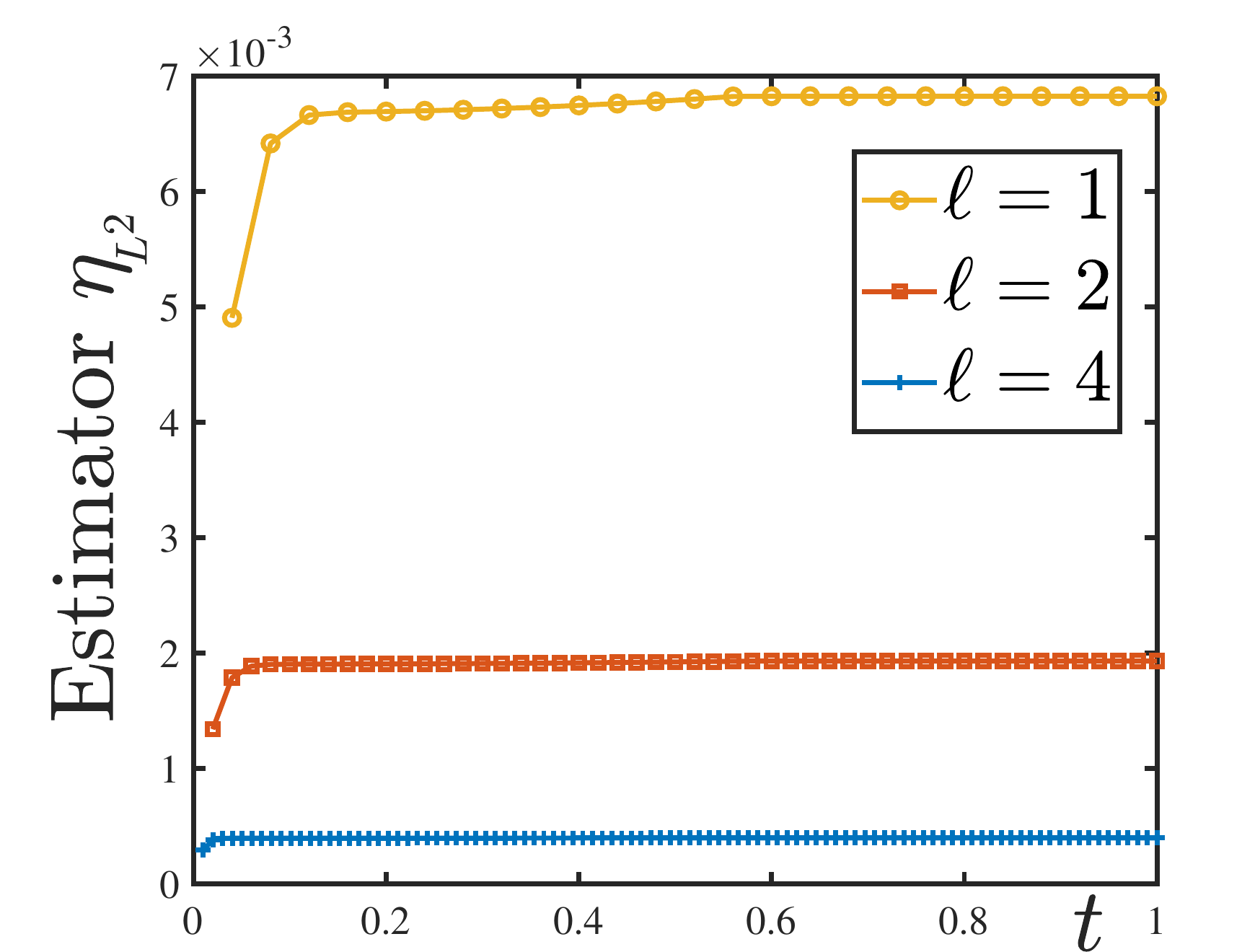}
\end{subfigure}
\begin{subfigure}{.48\textwidth}
\includegraphics[width=.8\textwidth]{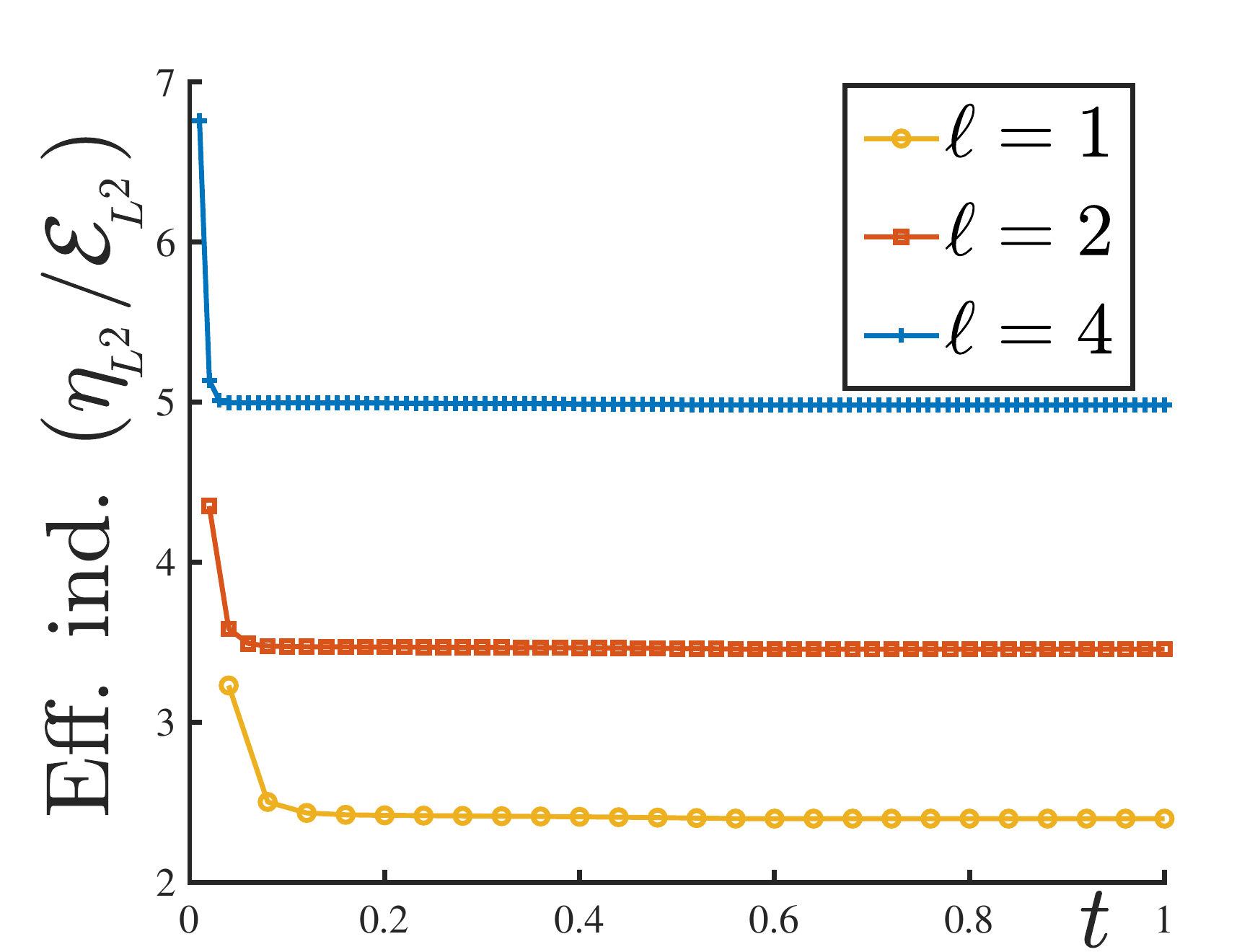}
\end{subfigure}
\caption{[\Cref{sec:NumTestDeg}] Reliability estimator $\eta_{L^2}$ (left) and its effectivity index (right). }\label{fig:DegEffL2}
\end{figure}
\begin{figure}[h!]
\begin{subfigure}{.48\textwidth}
\includegraphics[width=.8\textwidth]{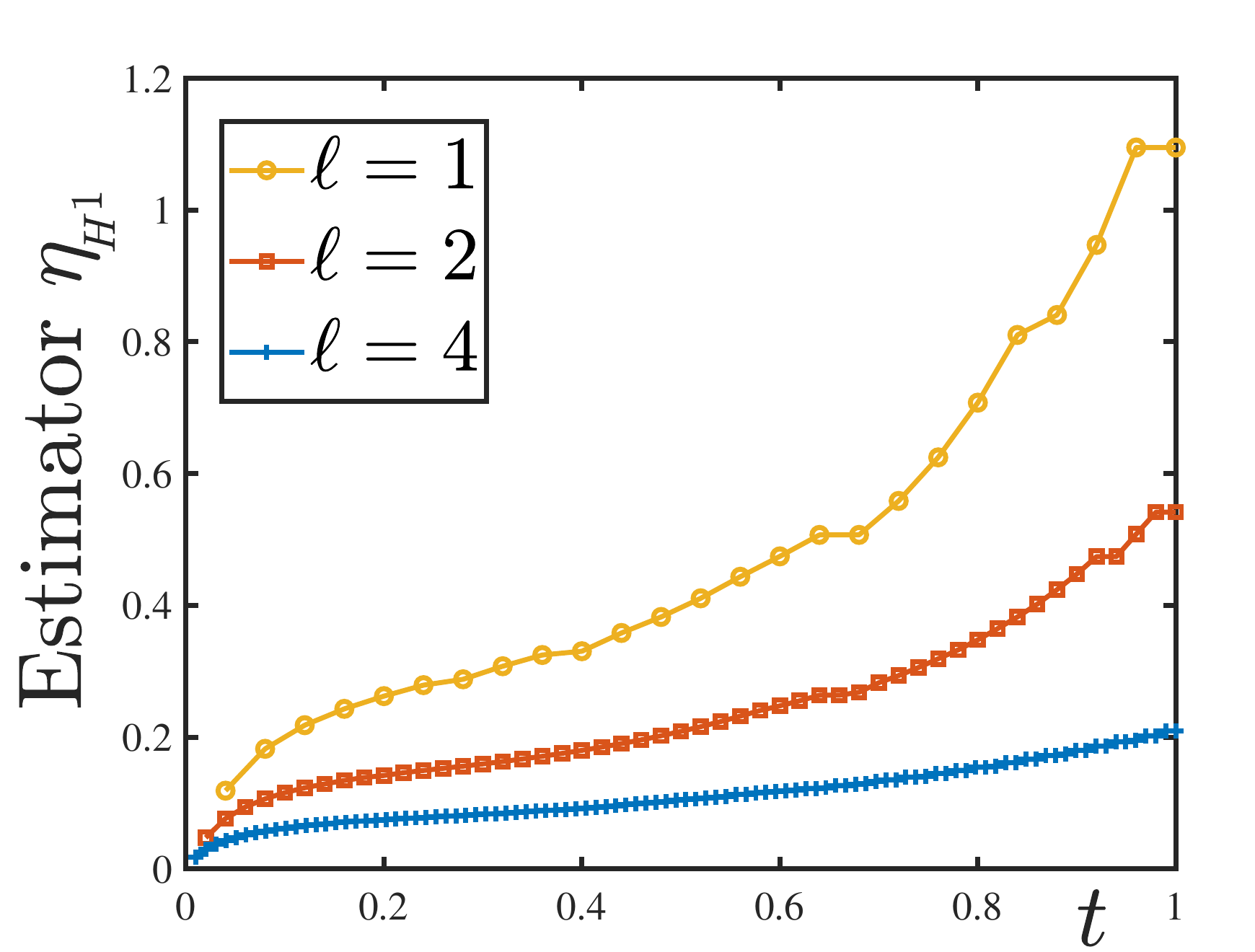}
\end{subfigure}
\begin{subfigure}{.48\textwidth}
\includegraphics[width=.8\textwidth]{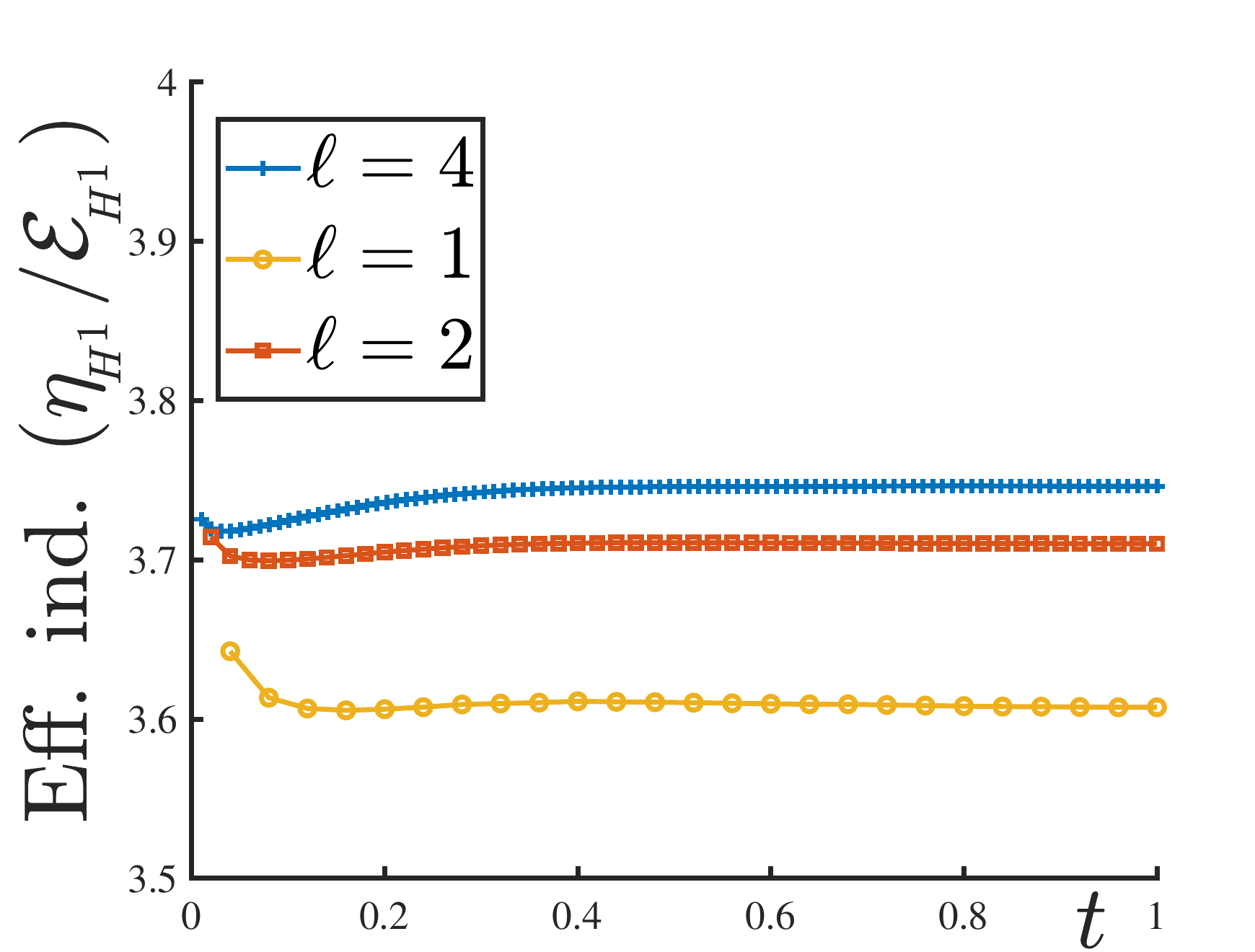}
\end{subfigure}
\caption{[\Cref{sec:NumTestDeg}] Reliability estimator $\eta_{H^1}$ (left) and its effectivity index (right). }\label{fig:DegEffH1}
\end{figure}
\begin{figure}[h!]
\begin{subfigure}{.48\textwidth}
\includegraphics[width=.8\textwidth]{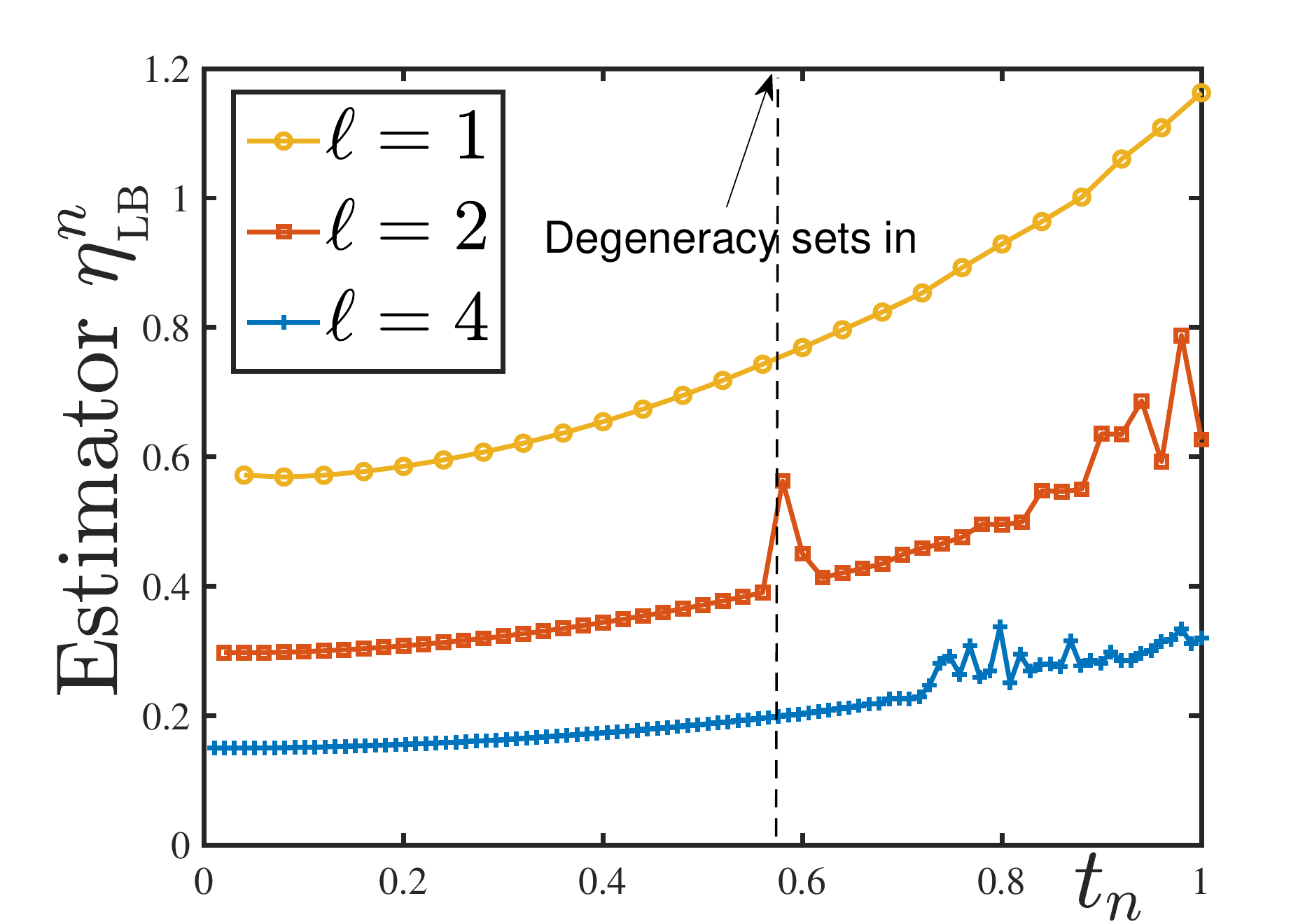}
\end{subfigure}
\begin{subfigure}{.48\textwidth}
\includegraphics[width=.8\textwidth]{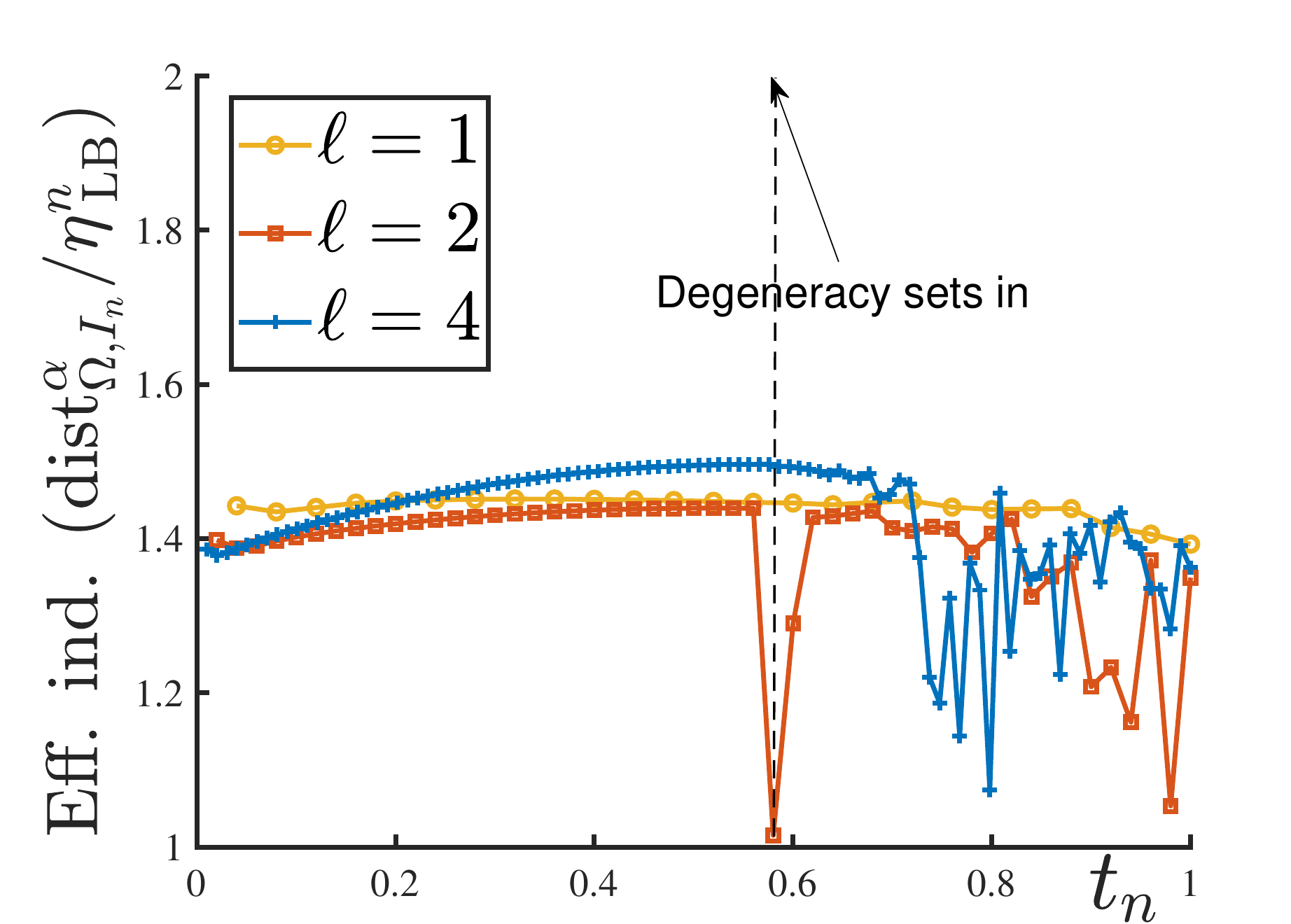}
\end{subfigure}
\caption{[\Cref{sec:NumTestDeg}] Lower bound $\eta^n_{\mathrm{LB}}$ (left) and its effectivity index (right). }\label{fig:DegEffLB}
\end{figure}

\begin{figure}[h!]
\begin{subfigure}{.32\textwidth}
\includegraphics[width=\textwidth]{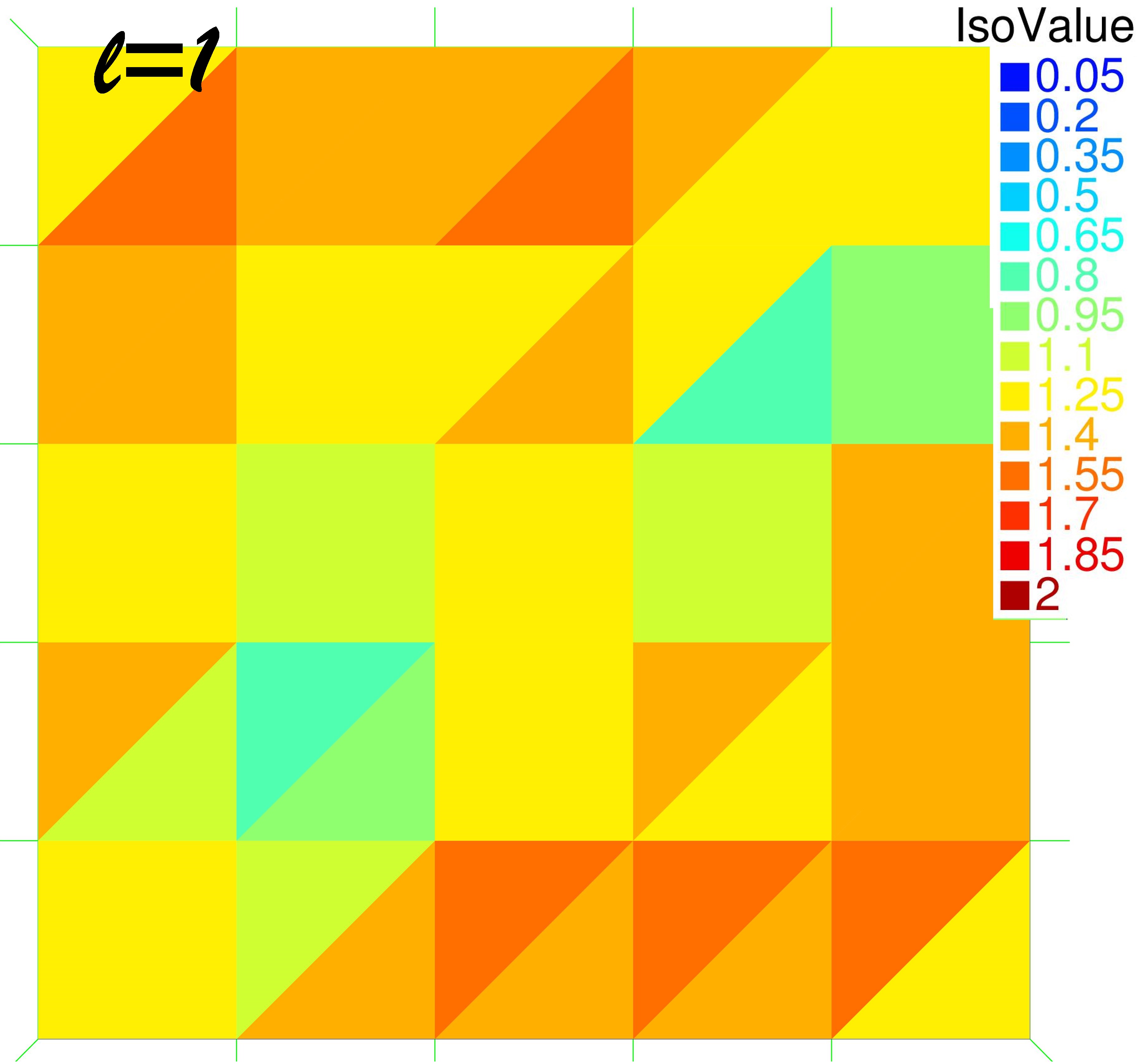}
\end{subfigure}
\begin{subfigure}{.32\textwidth}
\includegraphics[width=\textwidth]{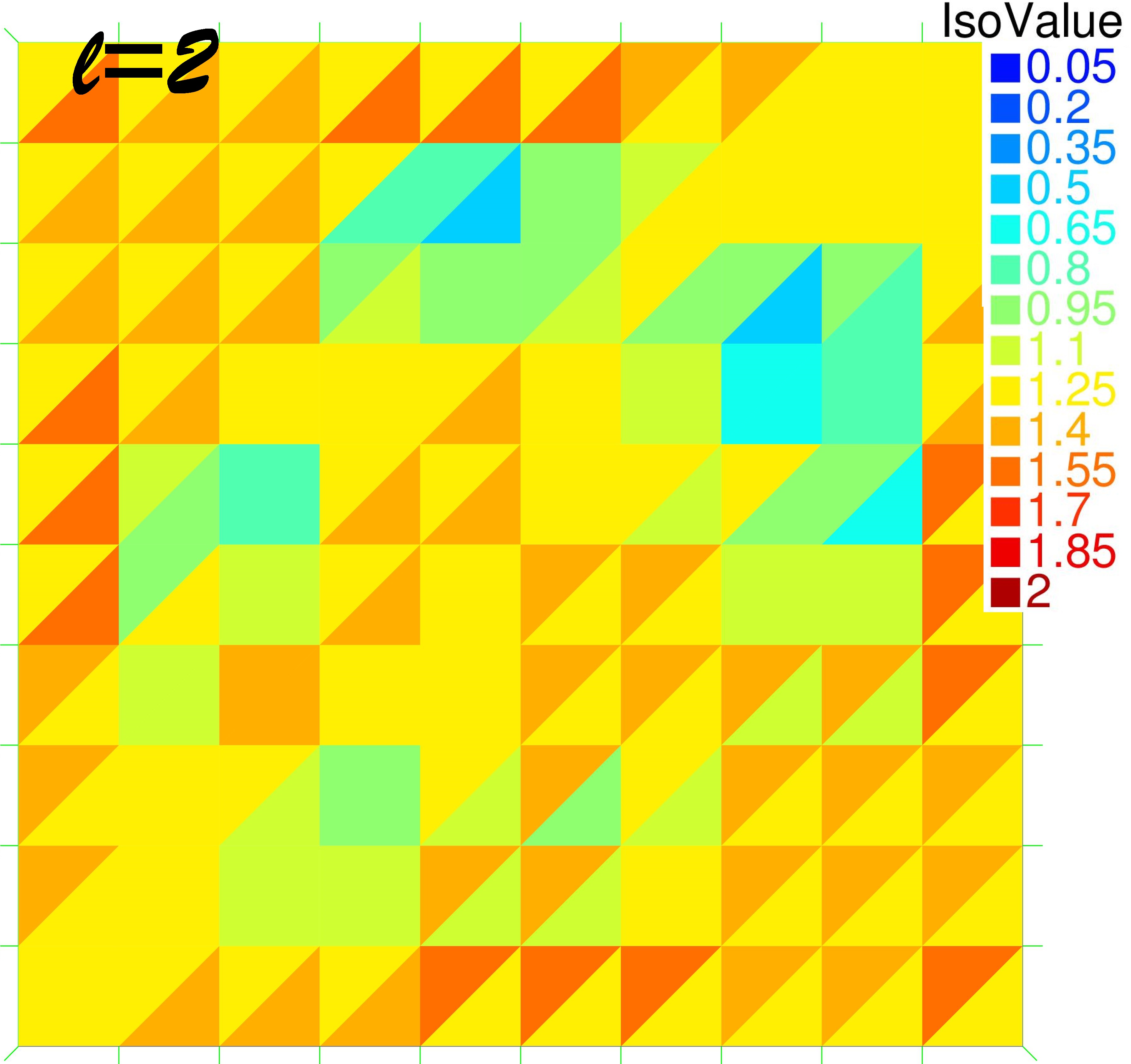}
\end{subfigure}
\begin{subfigure}{.32\textwidth}
\includegraphics[width=\textwidth]{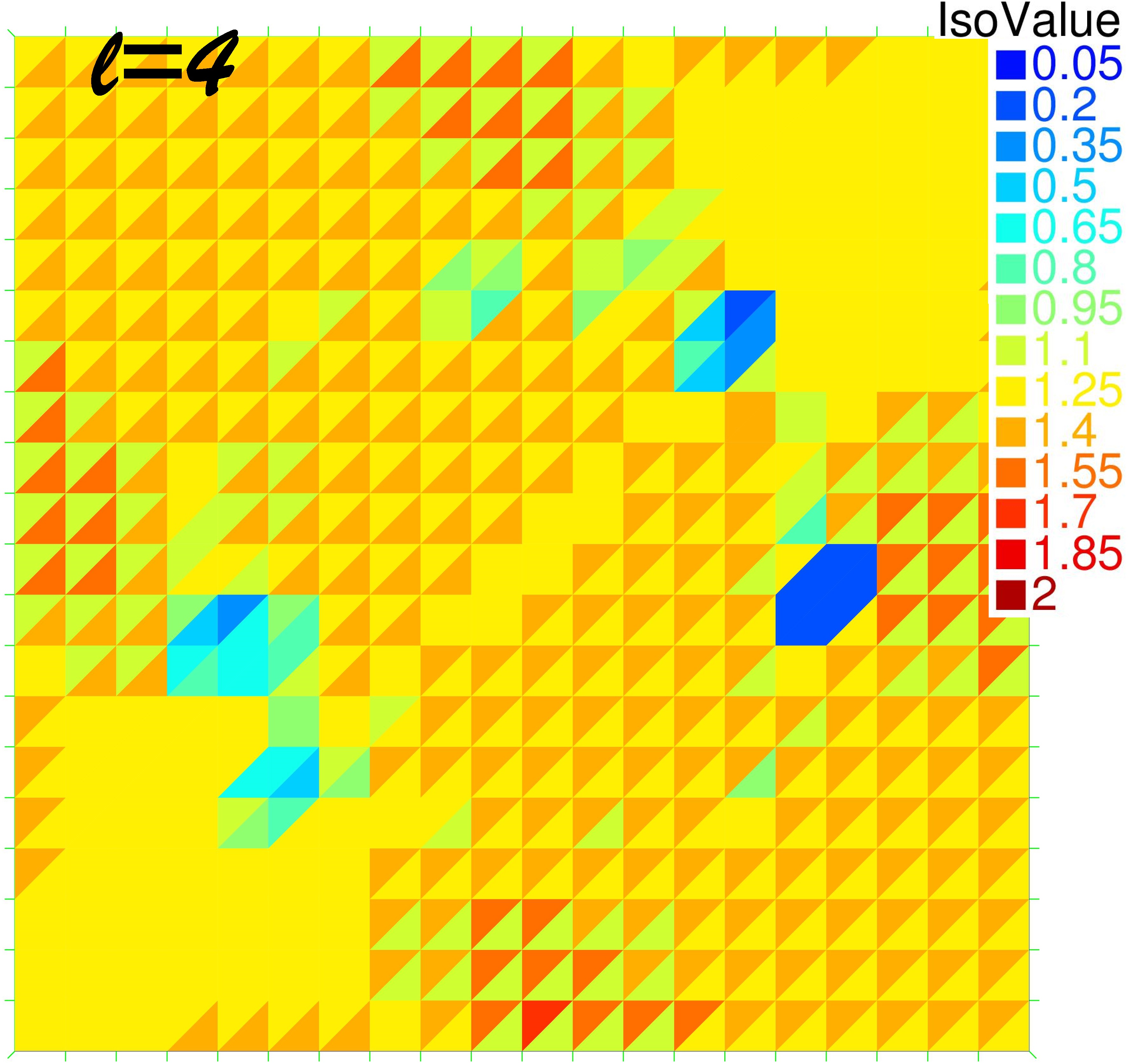}
\end{subfigure}
\caption{[\Cref{sec:NumTestDeg}] Local effectivity indices  \eqref{eq:EffInd_Loc} for $t_n=1$ and $\ell=1,2,4$.}\label{fig:LocalEffDeg}
\end{figure}

The effectivity indices for \Cref{theo:GlobalReliability,theo:Y_norm_guaranteed_efficiency} are defined as before. \Cref{fig:DegEffL2} shows the estimate $\eta_{L^2}$ and its effectivity. The effectivity index increases with $\ell$, despite $\eta_{L^2}$ decreasing monotonically, possibly since $s=s_{h\t}=1$ in a major portion of the domain towards the end of the simulation. \Cref{fig:DegEffH1} shows the results for $\eta_{H^1}$. The effectivity remains more stable in this case.  The effectivity indices for the lower bound are shown in \Cref{fig:DegEffLB} (right). An oscillation in the lower bound is observed for higher values of $\ell$. The reason for this behaviour is not clear. \Cref{fig:LocalEffDeg} shows the distribution of local space--time effectivity indices. They are close to 1 in most regions and only take a lower value close to the free-boundary $s=1$. Overall, we find these results satisfactory.

\subsection{Realistic case}\label{sec:NumTestReal}
\begin{figure}[h!]
\begin{subfigure}{.33\textwidth}
\includegraphics[width=\textwidth]{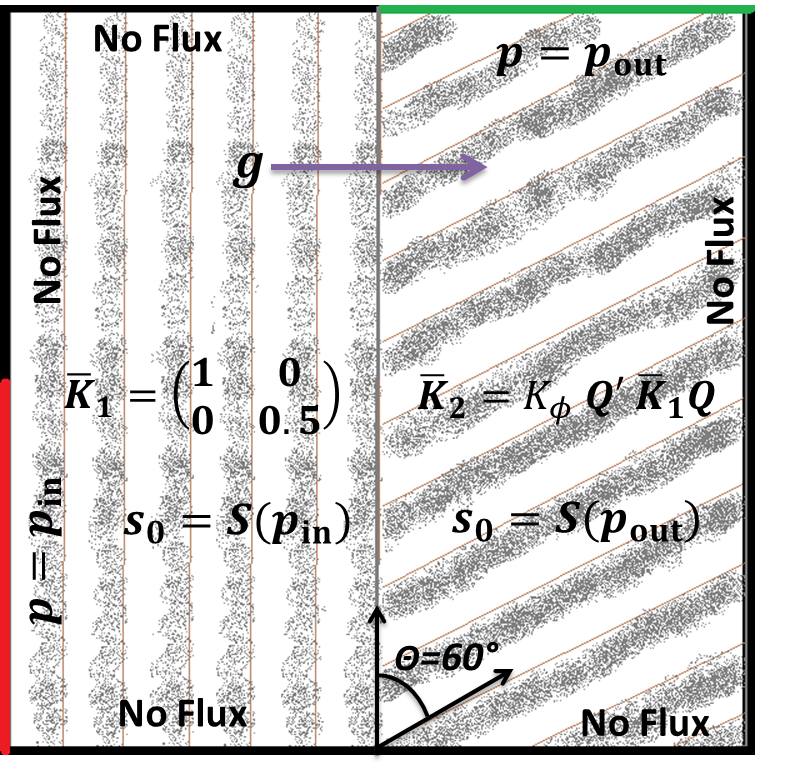}
\end{subfigure}
\begin{subfigure}{.32\textwidth}
\includegraphics[width=\textwidth]{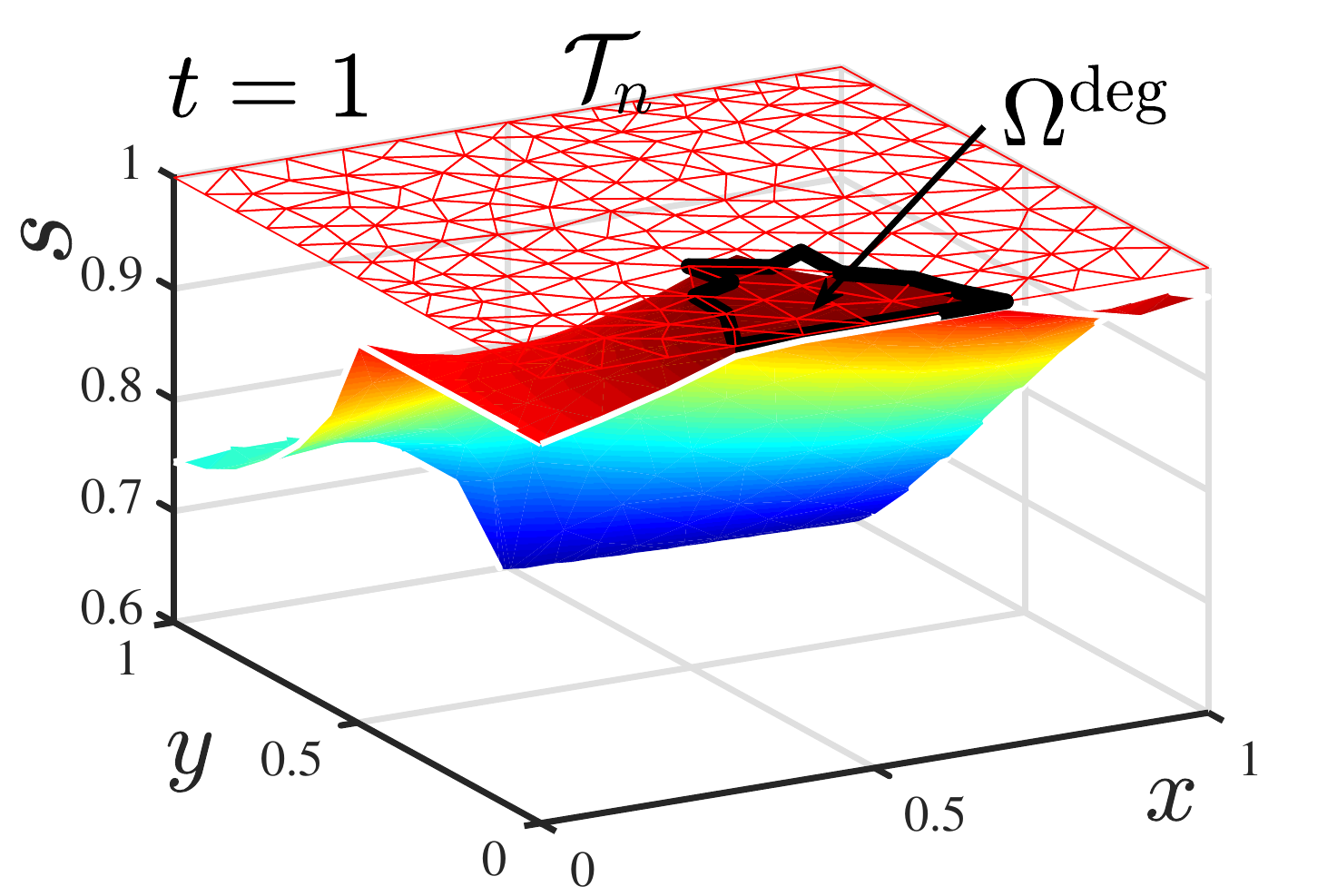}
\end{subfigure}
\begin{subfigure}{.32\textwidth}
\includegraphics[width=\textwidth]{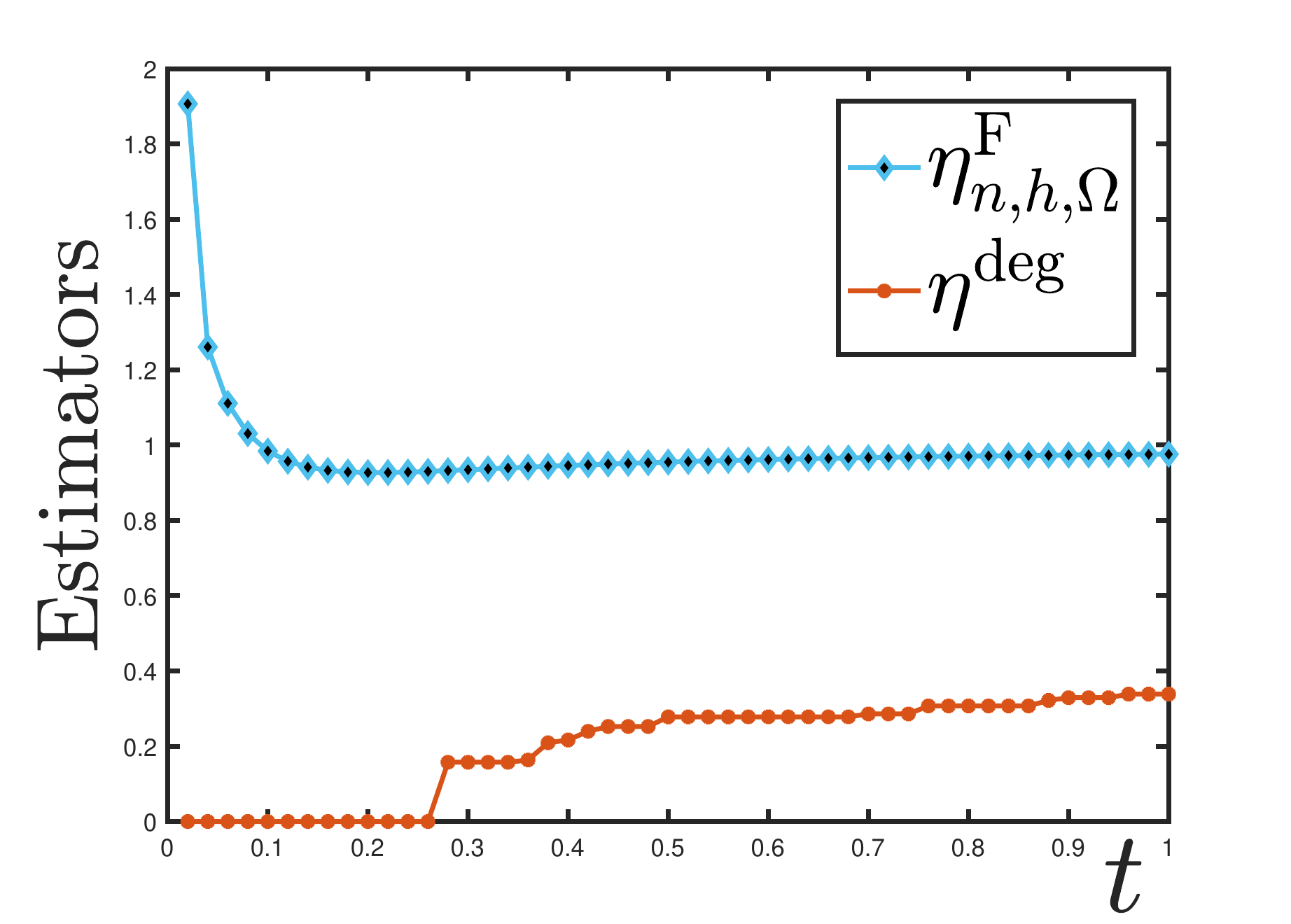}
\end{subfigure}
\caption{[\Cref{sec:NumTestReal}] Computational domain showing heterogeneities, initial, and boundary conditions (left). Saturation of the numerical solution ($S_{n,h}$) for $\ell=2$ at $t=1$. The mesh and the domain $\Om^{\mathrm{deg}}(1)$ containing the degenerate region is shown (center). The main estimators for the same simulation (right). The flux estimator $\eta^{\mathrm{F}}_{n,h,\Om}(t)$ contributes the most to the error, along with the degeneracy estimator $\etaDeg(t)$ which becomes non-zero only after the onset of degeneracy. \label{fig:Dom_Test3}}
\end{figure}

In this case, the domain and the functions given in \eqref{eq:NonlinearityTest1} are kept unchanged. The source term $f$ is made 0. However, the medium used is  heterogeneous and anisotropic with 
\begin{align}
\K=\begin{cases}
\K_1 &\text{ for } x<0.5,\\
K_{\phi}\, \bm{Q}^{\mathrm{T}}\K_1 \bm{Q} &\text{ for } x\geq 0.5,
\end{cases}
\text{ where }
\K_1:=
\begin{bmatrix}
1 &0\\
0 &0.5
\end{bmatrix}, \; \bm{Q}:=\begin{bmatrix}
\cos\theta &-\sin\theta\\
\sin\theta &\cos\theta
\end{bmatrix}.
\end{align}
Here, $\theta$ represents a tilted alignment of the principle axes of $\K$, and $K_{\phi}$ represents a factor stemming from the change in porosity. The choice of $\theta=\pi/3$ and $K_{\phi}=0.1$ is fixed. 
Both Neumann and inhomogeneous Dirichlet boundary conditions are used for the computation. The initial condition used is discontinuous. The details are shown in \Cref{fig:Dom_Test3} (left). The input and output pressures are 
$$
p_{\mathrm{in}}=0.8,\quad p_{\mathrm{out}}=-3.
$$
 A nonuniform mesh is used for the computation. No exact solution is known for this system.

Degeneracy occurs in the system close to the interface $x=0.5$ at $ y=0$. This is caused by the jump in $\K$, but also partly by the no-flux boundary condition. The error caused by this additional component is estimated by adding to $[\etaDeg(t)]^2$,
$$
\frac{2}{D(1)|\Om^{\mathrm{deg}}|}\int_{\p\Om} \hat{\bm{n}}^{\mathrm{T}}\left (\int_{\Om^{\mathrm{deg}}} \K \vg \right ) [\Psi_{h\t}(t) -P_{\mathrm{M}}]_+.
$$

\begin{figure}[h!]
\begin{subfigure}{.48\textwidth}
\includegraphics[width=.8\textwidth]{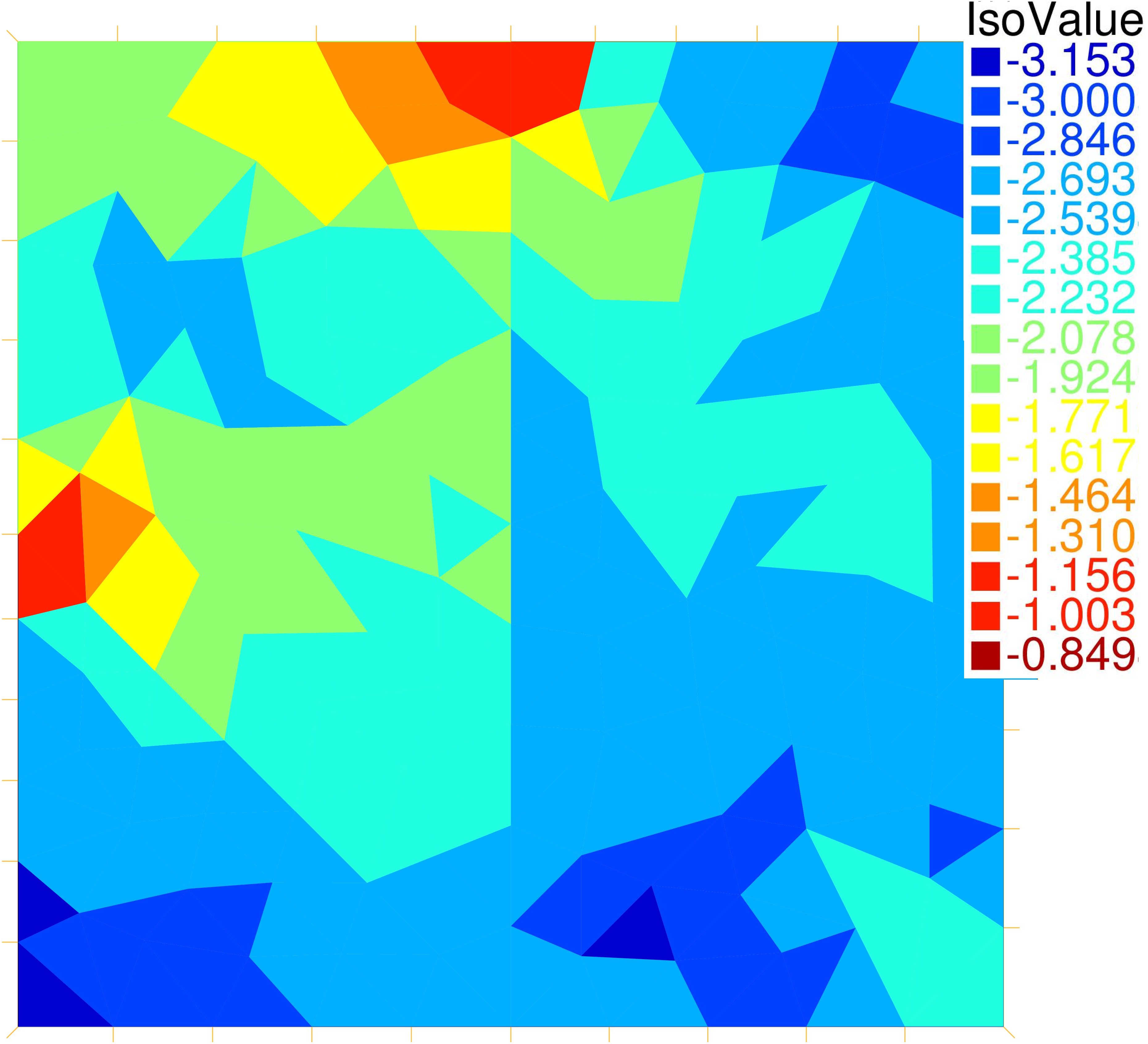}
\end{subfigure}
\begin{subfigure}{.48\textwidth}
\includegraphics[width=.8\textwidth]{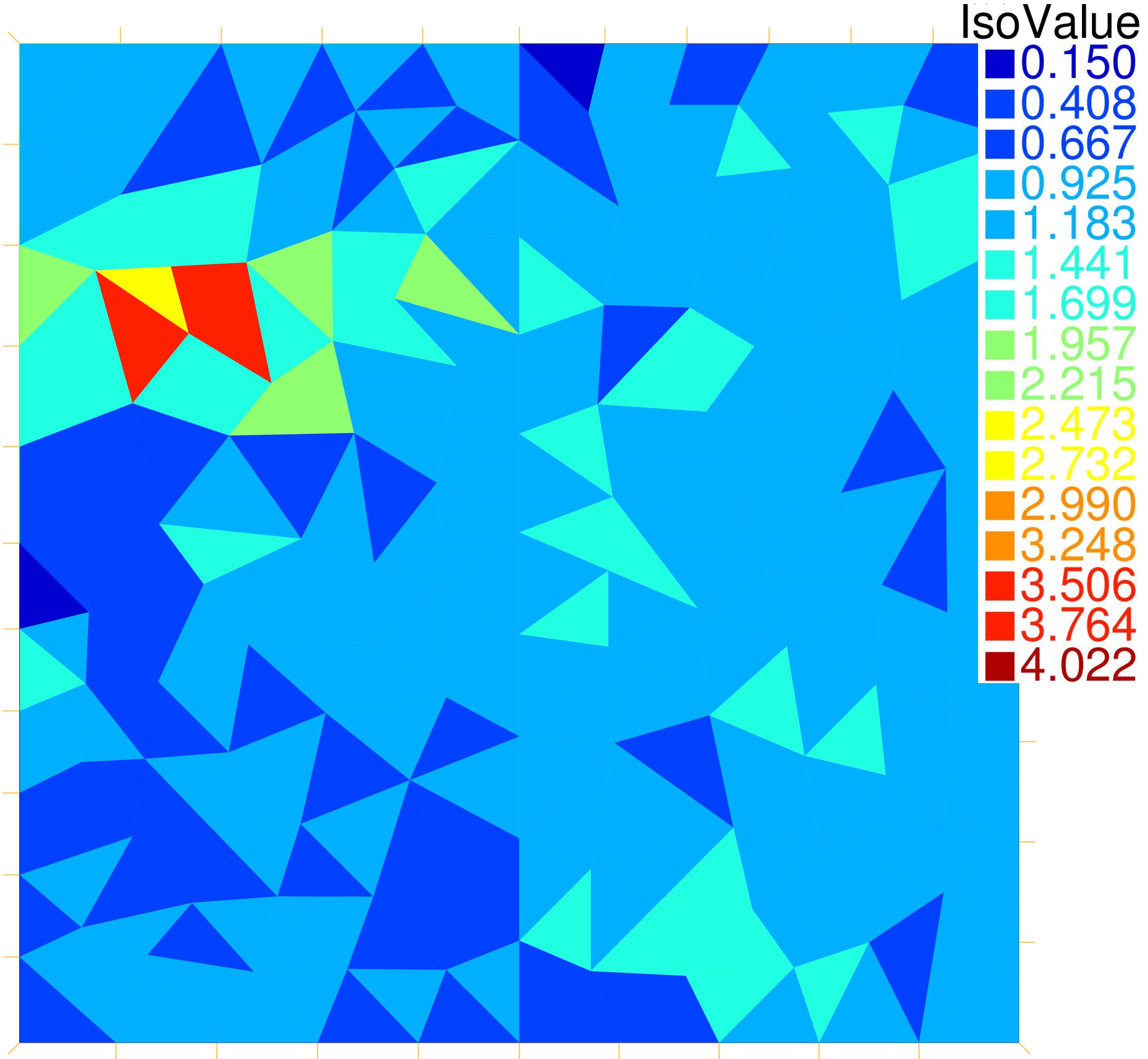}
\end{subfigure}
\caption{[\Cref{sec:NumTestReal}] Elementwise distribution of $\etaEq$ at time $t=1$ in $\log_{10}$-scale for $\ell=2$ (left). Estimators take larger values near the  inlet and the outlet. Local effectivity indices  defined in \eqref{eq:EffInd_Loc} (right).  Here,  the numerical solution for $\ell=4$ is used in the place of the exact solution. }\label{fig:RealCase}
\end{figure}
\Cref{fig:Dom_Test3} (center) shows the saturation distribution and degenerate zone for this problem at $t=1$. 
\Cref{fig:Dom_Test3} (right) shows the  main estimators for $\ell=2$.  The flux estimator is still the largest component, followed by $\etaDeg$. 
\Cref{fig:RealCase} (left) plot shows the spatial distribution of $\etaEq$ indicating high error concentrations  located around the inlet and the outlet. The right plot shows the local effectivity indices,  where the numerical solution for $\ell=4$ is used as the reference solution. Although the estimators vary by almost 3 orders of magnitude, the effectivity varies between 0.15--4 with most of the region having effectivity close to 1. We find this again quite satisfactory.

\begin{appendices}
\section{Iterative linearization}\label{App:linear}
In practice, since the problem \eqref{eq:TimeDiscrete} is nonlinear, its solution $p_{n,h}$ cannot be directly enumerated, and linearization iterations have to be used. We address this issue here.

\subsection{Linearization}\label{sec:AppLinScheme}
Set $S^{\bar{i}}_{0,h}=\Pi_{1,h} s_0$. For $n\in \{1,\dots,N\}$, let $\pinh{S}{n-1}$ be an approximation of $S_{n-1,h}$. Let $p_{n,h}^i\in \Vn$ denote the pressure at iteration $i\in \N$. Then, for a given $p^{i-1}_{n,h}\in \Vn$, we look for $\eth p^{i}_{n,h}:=p^{i}_{n,h}-p^{i-1}_{n,h}\in \Vn$ satisfying for all  $\f_h\in \Vn$,
\begin{align}\label{eq:MainNumScheme}
&  \tfrac{1}{\t_n}(L\,\eth p^{i}_{n,h},\f_h)+ (\K[\k(S(p_{n,h}^{i-1}))\nabla p^{i}_{n,h} + \bm{\xi}\, \eth p^{i}_{n,h}], \nabla \f_{h}) \nonumber \\
= &- \tfrac{1}{\t_n}(S(p_{n,h}^{i-1})-\pinh{S}{n-1}, \f_{h}) - (\K\vg \k(S(p_{n,h}^{i-1})),\del \f_h) + (f(S(p_{n,h}^{i-1}),\bm{x},t_n),\f_h),
\end{align}
Here, $(L,\bm{\xi}) \in \mathbf{L}^\infty(\Om;\R^{d+1})$ with $L\geq 0$, depends on the specific scheme used.  Since \eqref{eq:MainNumScheme} is linear with respect to $p^{i}_{n,h}$, it is directly computable. Observe that $\eth p^{i}_{n,h}=0$ if and only if $p^{i-1}_{n,h}$ solves \eqref{eq:TimeDiscrete} (provided $\pinh{S}{n-1}=S_{n-1,h}$) which shows that the schemes are consistent.

\begin{table}[H]
\centering
\begin{tabular}{|l|c|c|}
 \hline Scheme &$L,\;\bm{\xi}$ &Convergence\\
   \hline Picard  &$\scriptscriptstyle{0,\; \bm{0}}$ &\\
   \hline modified Picard \cite{celia_Picard} &$\scriptscriptstyle{S'(p^{i-1}_{n,h}) - \t_n \p_s f(S(p^{i-1}_{n,h}),\bm{x},t_n),\; \bm{0}}$ &Linear\\
 \hline
Newton \cite{bergamaschi1999mixed} &$\scriptscriptstyle{S'(p^{i-1}_{n,h}) - \t_n \p_s f(S(p^{i-1}_{n,h}),\bm{x},t_n)}$, $\scriptscriptstyle{\k'(S(p_{n,h}^{i-1}))[\nabla p^{i-1}_{n,h}+\vg]}$ &Quadratic\\
 \hline J{\"a}ger--Ka{\v{c}}ur \cite{Jager1991} &$\scriptscriptstyle{\sup_{p\in \R}\{(S(p)-S(p^{i-1}_{n,h})-\t_n (f(S(p))-f(S(p^{i-1}_{n,h}))))/(p-p^{i-1}_{n,h})\},\;\bm{0}}$ &Linear\\
 \hline L--scheme \cite{list2016study} &$\scriptscriptstyle{ L \text{ (constant)}\geq \frac{1}{2}\sup\{S'-\t_n\p_s f\},\; \bm{0}}$ &Linear\\\hline modified L--scheme \cite{MITRA2018} &$\scriptscriptstyle{S'(p^{i-1}_{n,h}) - \t_n \p_s f(S(p^{i-1}_{n,h}),\bm{x},t_n) + M\t_n \; (M>0 \text{ constant}), \; \bm{0}}$ &Linear\\\hline
\end{tabular}
\caption{Different iterative linearization schemes commonly used for Richards equation \eqref{eq:Richards1}. They fit into the common framework   \eqref{eq:MainNumScheme}. The corresponding $L$ and $\bm{\xi}$ quantities are displayed along with the convergence characteristics of the schemes. }\label{table:LinearIterative}
\end{table}

\Cref{table:LinearIterative} lists the commonly used schemes and the $L$ and $\bm{\xi}$ associated with them. The Picard scheme is generally unstable for the Richards equation. The modified Picard scheme \cite{celia_Picard} is linearly converging and the Newton method \cite{bergamaschi1999mixed} is quadratically converging  for an initial guess close to the solution of the nonlinear problem. However, the convergence is not guaranteed for degenerate cases.  The J{\"a}ger--Ka{\v{c}}ur scheme \cite{Jager1991} and the L--scheme \cite{list2016study} are unconditionally stable, meaning that they converge linearly, independent of the initial guess even in degenerate cases and for discontinuous initial conditions. However, a global supremum has to be computed for the J{\"a}ger--Ka{\v{c}}ur scheme, whereas, L--scheme converges slowly compared to the schemes mentioned above. The modified L--scheme \cite{MITRA2018} preserves the stability  of the L--scheme while being faster than the modified Picard scheme.

 For a given $n\in \{1,\dots,N\}$, let the linear iteration \eqref{eq:MainNumScheme} be terminated at $i=\bar{i}$. From the sequence $\{\pinh{p}{n}\}_{n=1}^N$, the space--time discrete total pressure and saturation are defined as
\begin{align}\label{eq:PsiDiscreteIter}
\pinh{\Psi}{n}:=\calK(\pinh{p}{n}), \text{ and } \pinh{S}{n}:= \Sf(\pinh{\Psi}{n})\overset{\eqref{eq:ConnectionVariables}}=S(\pinh{p}{n})\quad  \forall n\in \{1,\dots,N\},
\end{align}
analogous to \eqref{eq:PsiDiscrete}. Replacing $(\Psi_{n,h}, S_{n,h})$ by $(\pinh{\Psi}{n},\pinh{S}{n})$, we compute the time continuous solutions $\Psi_{h\t}$ and $s_{h\t}$ by following the steps of \Cref{sec:time_continuous}.

\subsection{Equilibrated flux} 
The terms $\Snt$ and $\Fnt$ from \eqref{eq:SourceFluxTerms}  are redefined as
\begin{subequations}\label{eq:SourceIterFluxTerms}
\begin{align}
\Snt &:=\left(f( S_{n,h}^{\bar{i}-1},\bm{x},t_n) - \tfrac{1}{\t_n}( S_{n,h}^{\bar{i}-1}- \pinh{S}{n-1}) -  L\,\eth \pinh{p}{n}\right ) \big|_{I_n},\\
\Fnt &:= \left [\k( S_{n,h}^{\bar{i}-1})\nabla p_{n,h}^{\bar{i}} +\vg \k( S_{n,h}^{\bar{i}-1}) + \bm{\xi}\,\eth \pinh{p}{n} \right ]\big|_{I_n}.
\end{align}
\end{subequations}
Observe that upon rearranging \eqref{eq:MainNumScheme}, $\Snt$ and $\Fnt$ play the role of the source-like and flux-like terms, just as in \eqref{eq:SourceFluxTerms}. Moreover, $\Snt$ and $\Fnt$ converge to their definitions in  \eqref{eq:SourceFluxTerms} if the iterate $\pinh{p}{n}$ converges to $p_{n,h}$.

The equilibrated flux $\sith\in \Hdiv$ is then constructed as stated in \Cref{def:flux_construction_1}.

\subsection{Estimators} The estimators in \eqref{eq:etaEq_def}--\eqref{eq:etaIni_def} and \eqref{eq:estimators2} are defined exactly the same way replacing $(\Psi_{n,h}, S_{n,h})$ by $(\pinh{\Psi}{n},\pinh{S}{n})$. The linearization estimators for the source-like and flux-like terms are introduced for $\om\subseteq \Om$ and $n\in\{1,\dots,N\}$ as
\begin{subequations}\label{eq:etaLin_def}
\begin{align}
&\eta^{\mathrm{lin},1}_{n,\om}:= C_{\mathrm{P},\om}\,h_{\om}\,\left \|\tfrac{1}{\t_n} (\pinh{S}{n} - S_{n,h}^{\bar{i}-1}- L\, \eth p_{n,h}^{\bar{i}})- (f(\pinh{S}{n})-f(S_{n,h}^{\bar{i}-1}))\right \|_{\om},\\
& \eta^{\mathrm{lin},2}_{n,\om} := \left \|\K^{\frac{1}{2}} \left  ((\k(\pinh{S}{n})-\k(S_{n,h}^{\bar{i}-1}))\, [\nabla \pinh{p}{n} +\vg]+ \bm{\xi}\,\eth p_{n,h}^{\bar{i}}\right )\right \|_\om,
\end{align}
\end{subequations}
($C_{\mathrm{P},\om}>0$ is the Poincar\'e constant). The new total estimator becomes, for $t\in I_n$,
\begin{align}\label{eq:NewEtaR}
\eta_{\mathcal{R}}(t):= \left [ \sum_{K\in \calT_n} [\etaEq(t)  + \etaOscS]^2 \right ]^{\frac{1}{2}}+ \etaOscTom(t) +\eta^{\mathrm{osc}}_{n,\Om}(t)  +\eta^{\mathrm{lin},1}_{n,\Om}.
 \end{align} 

 \begin{remark}
Observe that in \eqref{eq:etaLin_def} to define $\eta^{\mathrm{lin},1}_{n,\Om}$, we have used the $L^2$ norm instead of the $\Kn$-norm, which is costly to evaluate at every iteration. Observe also that only $\eta^{\mathrm{lin},1}_{n,\Om}$ appears in \eqref{eq:NewEtaR}.
\end{remark}

\subsection{Adaptive linearization}
Inspired by \cite{ern2013adaptive}, we propose the following adaptive algorithm for the linearization:
 
 \begin{algo}[Adaptive linearization]\label{algo:1}
 For a fixed $\gamma\in (0,1)$ and $n\in \N$, let $S^{\bar{i}}_{n-1,h}\in L^2(\Om)$ and $p^0_{n,h}\in H^1(\Om)$ be given. Then, for each $i\in\N$, solve  \eqref{eq:MainNumScheme}  until for some $i=\bar{i}$, upon computation of $\eta^{\mathrm{F}}_{n,h,\Omega}$ from \eqref{eq:etaEq_def} and $\eta^{\mathrm{lin},1}_{n,\Om}$, $\eta^{\mathrm{lin},2}_{n,\Om}$ from \eqref{eq:etaLin_def}, the following holds
  \begin{align}\label{eq:LinearAdaptiveCriteria}
\eta^{\mathrm{lin},1}_{n,\Om} + \eta^{\mathrm{lin},2}_{n,\Om}\leq \g  \, \eta^{\mathrm{F}}_{n,h,\Omega}.
\end{align}
 \end{algo}

With Algorithm \ref{algo:1},  \Cref{theo:GlobalReliability,theo:Y_norm_guaranteed_efficiency}
are restated as
\begin{proposition}[Reliability and efficiency with linearization]
 Let $\{\pinh{\Psi}{n}\}_{n=1}^{N}\subset H^1_0(\Om)$ and $\{\pinh{S}{n}\}_{n=1}^{N}\subset H^1(\Om)$ be defined using the numerical scheme \eqref{eq:MainNumScheme}--\eqref{eq:PsiDiscreteIter} with stopping criteria set by Algorithm \ref{algo:1}. Let $\Psi_{h\t}\in C(0,T;H^1_0(\Om))$ with $s_{h\t}=\Sf(\Psi_{h\t})\in W^{1,\infty}(0,T;H^1(\Om))$,  be their time-continuous interpolates as defined  in \eqref{eq:timeContSol}, with $(S_{n,h},\Psi_{n,h})$ replaced by $(\pinh{S}{n},\pinh{\Psi}{n})$. Let the estimators $\eta^{\mathrm{lin},j}_{n,\Om}$, $j=1,2$, be defined in \eqref{eq:etaLin_def} and $\eta_{\res}$ in \eqref{eq:NewEtaR}. Then
\begin{enumerate}[label=(\alph*)]
\item \textbf{Reliability:}  Under the assumptions of \Cref{theo:GlobalReliability}, the  estimates \eqref{eq:GlobalReliabilityEstimates} hold.
\item \textbf{Efficiency:}  Under the assumptions of \Cref{theo:Y_norm_guaranteed_efficiency}  and given that $\g$ is smaller than a threshold independent of the discretization, the estimate \eqref{eq:Y_norm_guaranteed_efficiency_lower_global}
holds.
\end{enumerate} 
\end{proposition}

\noindent
 The proofs are simple extensions to the proofs of \Cref{theo:GlobalReliability,theo:Y_norm_guaranteed_efficiency}.

\subsection{Numerical study}\label{app:num_linear}
\begin{figure}[h!]
\begin{subfigure}{.48\textwidth}
\includegraphics[width=.8\textwidth]{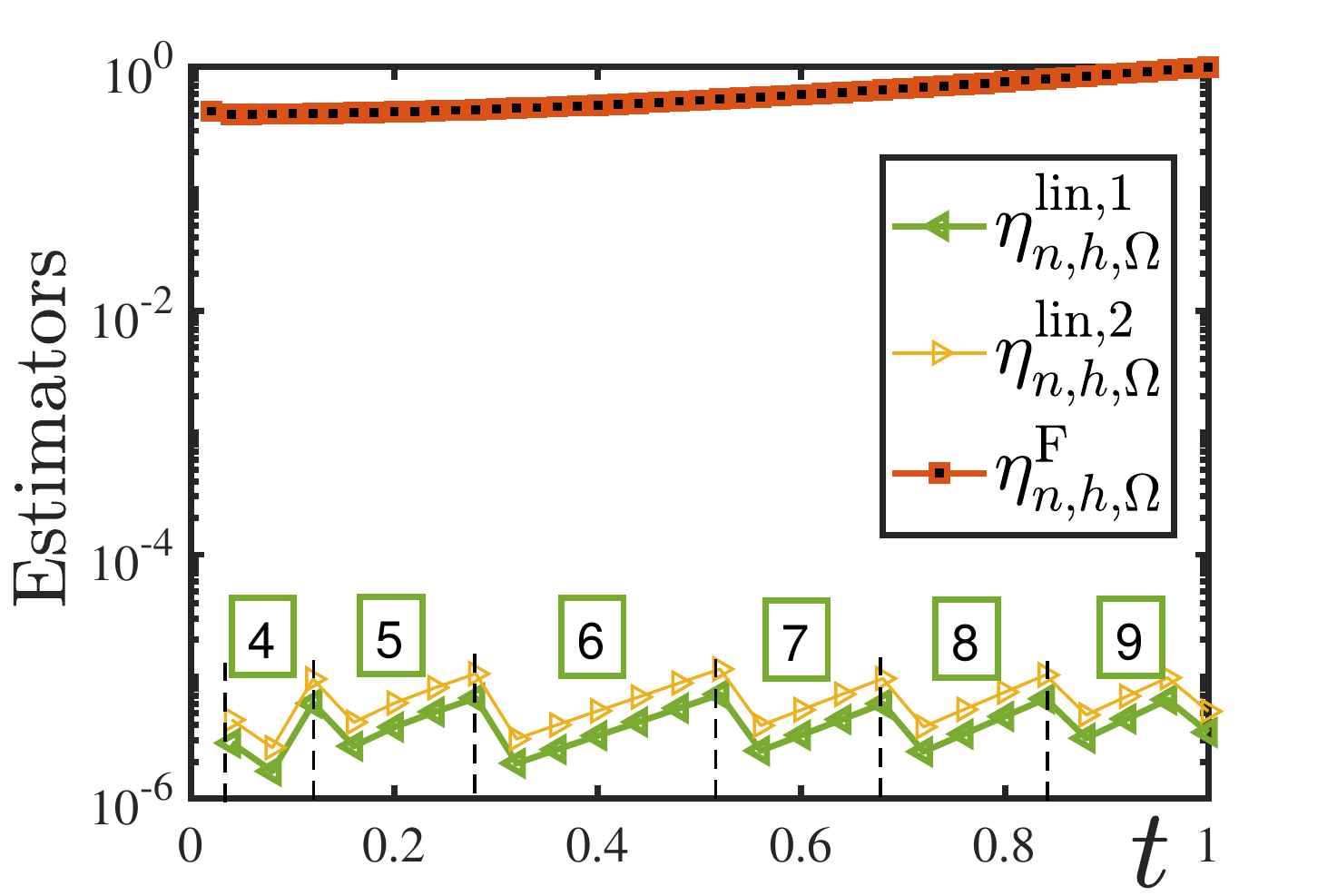}
\end{subfigure}
\begin{subfigure}{.48\textwidth}
\includegraphics[width=.8\textwidth]{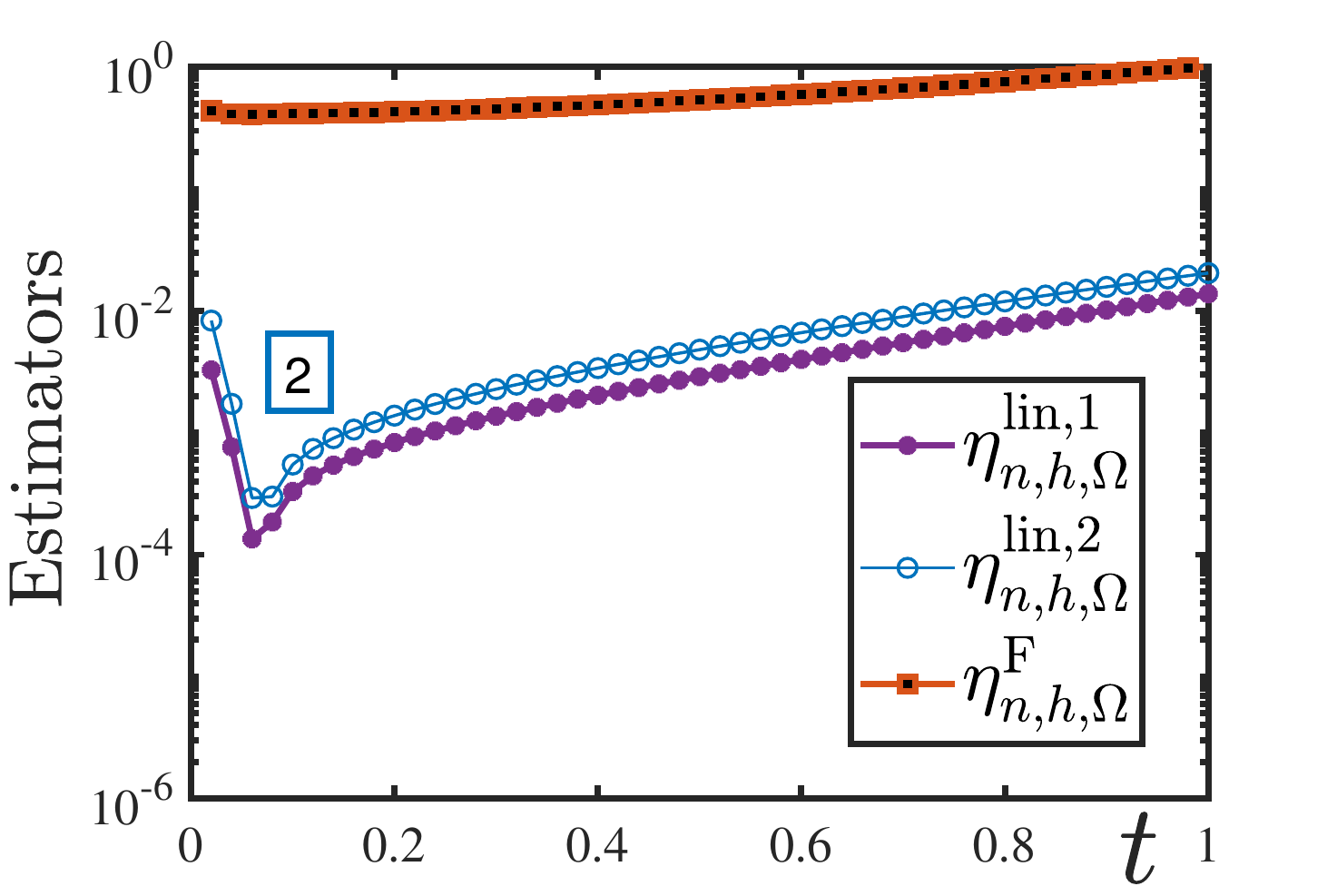}
\end{subfigure}
\caption{[\Cref{sec:NumTestNonDeg}, with adaptive linearization] The usual fixed error vs. the adaptive approach for linearization.  Here $\ell=2$.  The linearization estimators from \eqref{eq:etaLin_def} are plotted along with $\eta^{\mathrm{F}}_{n,h,\Om}$. Fixed error approach using $\|p^{\bar{i}}_{n,h}-p^{\bar{i}-1}_{n,h}\|_{\Kh(\Om)}\leq 10^{-4}$ as the stopping criterion (left).  Adaptive approach using Algorithm \ref{algo:1} with $\g=0.1$ (right).  The iterations required per time step are mentioned in the square boxes. They increase with time for the fixed approach and remain constant at 2 for the adaptive approach. \label{fig:AdapLin}}
\end{figure}

We present the results for the test cases from \Cref{sec:NumTestNonDeg}.
For this purpose, we use the modified L--scheme because of its stability and speed as discussed in Section \ref{sec:AppLinScheme}. The expression taken from \Cref{table:LinearIterative} becomes $L=S'(p^{i-1}_{n,h})+M\t_n$, $\bm{\xi}=\bm{0}$ with $M=1$ fixed throughout. Linear iterations are stopped when $\|p^{\bar{i}}_{n,h}-p^{\bar{i}-1}_{n,h}\|_{\Kh(\Om)}\leq 10^{-4}$. This  \textit{fixed error approach} is compared with the \textit{adaptive approach} which follows Algorithm \ref{algo:1}. For this purpose, $\g=0.1$ is chosen. \Cref{fig:AdapLin}  and \Cref{table:AdapLin} show that the adaptive approach requires much fewer iterations while having negligible impact on the quality of solutions. Moreover, the number of iterations required is stable as opposed to the fixed error approach.

\begin{table}[h!]
\centering
\begin{tabular}{|l|l|l|l|l|l|l|l|l|}
 \hline &\multicolumn{4}{c}{Fixed error approach}\vline &\multicolumn{4}{c}{Adaptive approach}\vline\\
 \hline
$\ell$ &\begin{tabular}{l}
avg.\\ iter.
\end{tabular} &$\eta^{\mathrm{lin},1}_{n,h,\Om}$ &$\eta^{\mathrm{lin},2}_{n,h,\Om}$ &$\eta_{\mathcal{R}}$ &\begin{tabular}{l}
avg.\\ iter.
\end{tabular} &$\eta^{\mathrm{lin},1}_{n,h,\Om}$ &$\eta^{\mathrm{lin},2}_{n,h,\Om}$ &$\eta_{\mathcal{R}}$\\
\hline
1 &7.72 &3.4e-6 &5.2e-6 &1.859 &2.00 &0.021 &0.038 &1.869\\
%\hline 1.6 &7.07 &3.4e-6 &3.8e-6 &1.230 &2 &0.016 &0.0253 &1.242\\
 \hline 2 &6.74 &5.7e-6 &5.6e-6 &0.998 &2.00 &0.014 &0.020 &1.088\\
%\hline 3 &6.29 &1.6e-6 &1.3e-6 &0.687 &2.03 &.0098 &0.137 &0.696\\
\hline 4 &5.72 &1.4e-6 &9.6e-7 &0.497 &1.98 &0.007 &0.009 &0.506\\\hline
\end{tabular}
\caption{[\Cref{sec:NumTestNonDeg}, with adaptive linearization] Average iterations required per time step together with $\eta^{\mathrm{lin},1}_{n,h,\Om}$, $\eta^{\mathrm{lin},2}_{n,h,\Om}$, and $\eta_{\mathcal{R}}$ at $t_n=1$ for the usual fixed error (left) and the adaptive linearization (right) approaches.}\label{table:AdapLin}
\end{table}

\section{Proofs of \Cref{sec:MaxPrin}}\label{App:proof}
We collect here the proofs of the statements of \Cref{sec:MaxPrin}.

\begin{proof}[\textbf{Proof of \Cref{pros:MaximumPrinciple}}]
We claim that $(S_{\mathrm{m}},p_{\mathrm{c}}(S_{\mathrm{m}}))$  serves as a subsolution of $(s,p)$ for a constant $\K$. From \eqref{eq:DefSMSm} we have that $\bar{S}_\mathrm{m}(t)> S(0)$ for some $t> 0$ only if $f_{\mathrm{m}}(S(0))\geq 0$ or $f(S(0),\bm{x},t)\geq 0$ a.e. in $(\bm{x},t)\in \Om\times \R^+$. Hence, if $S_\mathrm{m}(t)=\min(\bar{S}_\mathrm{m}(t),S(0))= S(0)$ in some interval $I\subseteq \R^+$ then $\p_t S_{\mathrm{m}}- f( S_{\mathrm{m}},\bm{x},t)\leq 0- f_{\mathrm{m}}(S(0)) \leq 0$. On the other hand, if $S_\mathrm{m}(t)= \bar{S}_{\mathrm{m}}(t)$ then $\p_t S_{\mathrm{m}}- f( S_{\mathrm{m}},\bm{x},t)\leq \p_t \bar{S}_{\mathrm{m}}- f_{\mathrm{m}}(\bar{S}_{\mathrm{m}})= 0$. Hence, 
\begin{align*}
&\p_t S_{\mathrm{m}}-\del\cdot[\K \, \k( S_{\mathrm{m}})\, (\del p_\mathrm{c}( S_{\mathrm{m}}) + \vg)]- f( S_{\mathrm{m}},\bm{x},t)=\p_t S_{\mathrm{m}}- f( S_{\mathrm{m}},\bm{x},t)\leq 0.
\end{align*}
Moreover, $p_\mathrm{c}(S_\mathrm{m})\leq p_\mathrm{c}(S(0))=0$ in relation to the boundary.
Thus, invoking  the comparison principle \cite{otto1996l1}, we conclude that $(S_{\mathrm{m}},p_{\mathrm{c}}(S_{\mathrm{m}}))$  is  a subsolution of $(s,p)$.
\end{proof}

\begin{proof}[\textbf{Proof of \Cref{pros:IntroVS}}] Let $J_1:=\min(J,p_l)<0$. For the sake of simplicity, let the space coordinate be translated such that  $\min\{\vg\cdot\bm{x}\}=0$. To show the lower bound of $\vs$ we use $\del(\vg\cdot\bm{x})=\vg$, and rewrite \eqref{eq:VSweak} as
\begin{align}\label{eq:VSinTranslatedCor}
(\K \k(S(\vs))\del [\vs +\vg\cdot\bm{x}],\del \f)=\left (\inf_{t\in \R^+} [f(S(\vs),\bm{x},t)]_-,\f\right ).
\end{align}
Selecting the test function $\f=[\vs -J_1 +\vg\cdot\bm{x}]_-\in H^1_0(\Om)$  (observe that $\f=0$ on $\p\Om$ since $\vg\cdot\bm{x} -J_1\geq 0$ for all $\bm{x}\in \Om$) one then obtains in the left hand side of \eqref{eq:VSinTranslatedCor},
\begin{align*}
&(\K \k(S(\vs))\del [\vs +\vg\cdot\bm{x}],\del [\vs-J_1 +\vg\cdot\bm{x} ]_-) \geq \int_\Om \k(S(\vs))\left |\K^{\frac{1}{2}} \del [\vs-J_1 +\vg\cdot\bm{x} ]_-\right |^2.
\end{align*}
Observe that $\f$ is nonzero only when $\vs\leq J_1-\vg\cdot\bm{x} \leq p_l$, implying $f(S(\vs),\bm{x},t)\geq 0$. Hence, the right hand side of \eqref{eq:VSinTranslatedCor} yields
\begin{align*}
\left (\inf_{t\in \R^+} [f(S(\vs),\bm{x},t)]_-,[\vs-J_1 +\vg\cdot\bm{x} ]_-\right )= 0.
\end{align*}
Hence, from \eqref{eq:VSinTranslatedCor}, one obtains $\vs\geq J_1 - \vg\cdot\bm{x}$. We obtain the upper bound by testing with $\f=[\vs -J+\vg\cdot\bm{x}-\max\{\vg\cdot\bm{x}\}]_+$ and following the arguments as before.
\end{proof}
\begin{proof}[\textbf{Proof of \Cref{pros:MaximumPrincipleB}}]
Observe that the choice of $J$ implies from \Cref{pros:IntroVS} that $\vs\leq 0$ and $S(\vs)\leq s_0$ a.e. in $\Om$. Moreover, from \eqref{eq:VSweak},
\begin{align*}
&\p_t \vs- \del\cdot(\K \k(S(\vs))[\del \vs +\vg])-f(S(\vs),\bm{x},t)\nonumber\\
\leq & 0- \del\cdot(\K \k(S(\vs))[\del \vs +\vg])-\inf_{\z\in \R^+} [f(S(\vs),\bm{x},\z)]_-= 0,
\end{align*}
since $f(S(\vs),\bm{x},t) \geq \inf_{\z\in \R^+} [f(S(\vs),\bm{x},\z)]_-$. Hence, similar to the proof of \Cref{pros:MaximumPrinciple}, the result follows from applying the comparison principle.
\end{proof}
\end{appendices}

\end{document}